\documentclass[11pt,reqno]{amsart}
\usepackage{amsmath}
\usepackage{cases}
\usepackage{mathrsfs}
\usepackage{amssymb}
\usepackage{amscd}
\usepackage{amsfonts,latexsym,amsmath,amsthm,amsxtra,amssymb,latexsym}
\usepackage[all,cmtip]{xy}
\RequirePackage{amsmath} \RequirePackage{amssymb}
\usepackage{color}
\usepackage{colordvi}
\usepackage{multicol}
\usepackage{hyperref}
\usepackage{mathtools}
\allowdisplaybreaks[3]
\usepackage[margin=1in]{geometry}
\usepackage{xcolor}
\hypersetup{
    colorlinks,
    linkcolor={red!50!black},
    citecolor={blue!100!black},
    urlcolor={blue!100!black},
    bookmarksdepth=3,
    bookmarksnumbered=true,
    bookmarksopen=true,
    bookmarksopenlevel=2
}
\usepackage{cite}

\newcommand{\bea}{\begin{eqnarray}}
\newcommand{\eea}{\end{eqnarray}}
\newcommand{\bna}{\begin{eqnarray*}}
\newcommand{\ena}{\end{eqnarray*}}

\numberwithin{equation}{section}

\theoremstyle{plain}
\newtheorem{lemma}{Lemma}[section]
\newtheorem{theorem}[lemma]{Theorem}

\newtheorem{proposition}[lemma]{Proposition}

\theoremstyle{definition}

\newtheorem{remark}{Remark}

\renewcommand{\Re}{\operatorname{Re}}
\renewcommand{\Im}{\operatorname{Im}}

\newcommand{\GL}{\operatorname{GL}}
\newcommand{\SL}{\operatorname{SL}}

\newcommand{\dd}{\mathrm{d}}

\newcommand{\eps}{\varepsilon}
\newcommand{\sgn}{\operatorname{sgn}}

\renewcommand{\bmod}[1]{\,(\mathrm{mod\,}{#1})}
\renewcommand{\i}{\mathrm{i}}
\makeatletter
\@namedef{subjclassname@2020}{%
  \textup{2020} Mathematics Subject Classification}

\title[]
{Non-vanishing of central $\mathrm{GL}_3$ $L$-values
in the level aspect}

\author{Hui Wang}
\address{The University of Hong Kong Shenzhen Institute of Research and Innovation,
Shenzhen, Guangdong 518057, China}
\email{wh0315@mail.sdu.edu.cn}

 \author{Shilun Wang}
\address{Tongji University, Shanghai 200092, China}
\email{slwang@tongji.edu.cn}
\date{}

\begin{document}
\raggedbottom
\begin{abstract}
We study non-vanishing of central standard $L$-values of spherical
cuspidal $\mathrm{GL}_3$ newforms of prime level.
A logarithmic M\"obius mollifier with finite inverse Euler factors,
combined with a conjugate term from the functional equation, gives
a non-vanishing harmonic proportion of at least $1/2-\eps$ for
every $\eps>0$ and all sufficiently large prime levels.
The required mixed moment is evaluated by an asymmetric
approximate functional equation.
\end{abstract}
\maketitle

\tableofcontents

\section{Introduction}
\setcounter{equation}{0}

Non-vanishing of central \(L\)-values is a basic problem in analytic number theory. A standard approach is to approximate the reciprocal by a short Dirichlet polynomial and apply Cauchy--Schwarz to the resulting mollified moments. This method originates in Selberg's work \cite{Selberg} and was developed further by Levinson and Conrey \cite{Levinson,Conrey}. 

For Dirichlet $L$-functions, positive-proportion non-vanishing has been established in families of primitive characters of large modulus and in quadratic families \cite{BM,ISD,Sound}. Further developments of mollification and improvements of the non-vanishing proportions can be found in \cite{MV,Bui,KMN,QW}.

For automorphic $L$-functions on $\GL_2$, 
we have non-vanishing results in the spectral aspect for Maass forms. See \cite{Liu,BHS,WW,LQ}. For the weight aspect for holomorphic forms, we can see \cite{Luo,BF}. In the level aspect for newforms, we refer to the work \cite{Duke,KM,IS,PWZ,RouymiTrace,RouymiMoll,BFLevel}. Results for twists of a fixed form are given in \cite{Rohrlich,FH,OS,BFKMMS}, with further moment estimates for quadratic twists in \cite{SY,LiQuad}. Related results for Rankin--Selberg convolutions and for non-vanishing at spectrally varying special points are given in \cite{KMV,LuoWeyl,LuoRefined,XuSpecial,QiMaass,QiRS} and the references therein.

For $\GL_3$, earlier work treats character twists of fixed forms \cite{HK,LuoGLn} and symmetric square families \cite{BlomerSym,KhanSym,BFMaassSym}. The Kuznetsov formulas of Blomer \cite{Blomer} and Buttcane \cite{ButtcaneSpectral} provide tools for general \(\mathrm{GL}_3\) families. Building on the kernel estimates of Blomer--Buttcane
\cite{BBSubconvex}, Huang--Liu--Xu \cite{HLX} used a suitable
mollifier to prove a positive harmonic non-vanishing proportion
for forms of level one in generic short spectral intervals.

In the level aspect, Blomer--Buttcane--Maga \cite{BBM} developed the Kuznetsov trace formula of the congruence subgroup and its applications, while Blomer--Comtat \cite{BC} supplied local and spectral input for prime level averages. Tsuzuki \cite{Tsuzuki} proved simultaneous non-vanishing in families whose conductors have a fixed ramified part and an additional varying prime.

Motivated by these developments, we study non-vanishing of central standard $L$-values of degree three for some families of cuspidal automorphics of  $\GL_3$. We state the following non-vanishing theorem in \S\,\ref{MT} and the underlying mollified moment formulas in \S\,\ref{subsec:intro-mollified-moments}.

\subsection{The family and the non-vanishing theorem}
\label{MT}
Let $p$ be an odd prime and define
\begin{equation*}
\Gamma_0(p)\colon \!\!= \{\gamma=(\gamma_{ij})_{1\leq i,j\leq3}\in\SL_3(\mathbb{Z}) \colon  \gamma_{31}\equiv\gamma_{32}\equiv0\bmod p\}.
\end{equation*}
We denote its index in $\SL_3(\mathbb{Z})$ by
\begin{equation*}
P\colon \!\!=[\SL_3(\mathbb{Z}): \Gamma_0(p)]=p^2+p+1.
\end{equation*}
Let $\mathcal{B}(p)$ denote the family of cuspidal automorphic representations $\pi$ on $\GL_3$ of conductor $p$ with trivial central character and spherical archimedean component. We normalize the Fourier coefficients by $A_\pi(1,1)=1$. Write $\mu_\pi$ for the archimedean spectral parameters and $L(s,\pi)$ for the standard $L$-function. In \S\,\ref{sec:preliminaries}, we specify the normalization and the local description at $p$.

Fix the test function $h$ as in \eqref{eq:positive-moment-test}. We define
\begin{equation*}
\mathcal H_h\colon \!\!= \int_{\Re(\mu)=0}h(\mu)  \operatorname{spec}(\mu)\,\dd\mu>0.
\end{equation*}
The density $\operatorname{spec}(\mu)\,\dd\mu$ is specified in \S\,\ref{subsec:kuz-transforms}. Here the finite adjoint $L$-function is defined for $\Re(s)>1$ as in \cite[\S\,2.4]{BZ}.
\begin{equation*}
L(s,\pi,\mathrm{Ad})\colon \!\!=\dfrac{L(s,\pi\times\widetilde\pi)}{\zeta(s)} =\prod_q(1-q^{-s})L_q(s,\pi\times\widetilde\pi).
\end{equation*}

The total harmonic mass satisfies, as $p\to\infty$ through primes
(see \S\,\ref{subsec:total-harmonic-mass}),
\begin{equation*}
c_{\rm cusp}\sum_{\pi\in\mathcal B(p)}
 \frac{h(\mu_\pi)}{L(1,\pi,\mathrm{Ad})}
=P\mathcal H_h+O_{h,\varepsilon}(p^{1+\varepsilon}).
\end{equation*}
The positive constant $c_{\rm cusp}$ is fixed in \S\,\ref{subsec:kuz-normalization}.

\begin{theorem}
\label{thm:intro-weighted-nonvanishing}
Let $h$ be the fixed test function. For every $\eps>0$ and all sufficiently large primes $p$, we have
\begin{equation*}
\sum_{\substack{\pi\in\mathcal B(p)\\L(1/2,\pi)\ne0}}
\dfrac{h(\mu_\pi)}{L(1,\pi,\mathrm{Ad})}
\geq \left(\frac{1}{2}-\eps\right)\sum_{\pi\in\mathcal B(p)}
\dfrac{h(\mu_\pi)}{L(1,\pi,\mathrm{Ad})}.
\end{equation*}
\end{theorem}

The asymptotic for the total harmonic mass shows that the non-vanishing subfamily has harmonic mass of order $p^2$. The theorem provides a counterpart in the level aspect to the spectral result of Huang--Liu--Xu \cite[Theorem~1.1 and Remark~1.2]{HLX}.

The factor $L(1,\pi,\mathrm{Ad})^{-1}$ is the harmonic weight arising from the Kuznetsov formula. Removing it would require uniform control of an adjoint $L$-value approximation inside the mollified moments, together with the additional Hecke indices it introduces. For a treatment at level one, see \cite[\S\S\,2.4 and 3]{BZ}. Pointwise bounds for the adjoint values alone do not preserve the proportion $1/2$. We therefore retain the harmonic weight throughout.

\begin{remark}
The present paper treats the level aspect with the archimedean test function fixed. A natural continuation is to establish a hybrid non-vanishing theorem in which the prime level $p$ and the archimedean spectral parameter $T$ tend to infinity simultaneously. More precisely, one would seek a positive lower bound for the harmonic proportion of non-vanishing in the corresponding family, uniformly in a suitable joint range of
$p$ and $T$. Such a result would combine the level aspect studied here with the spectral aspect considered by Huang--Liu--Xu. Its proof would require hybrid estimates for the prime-level Kuznetsov formula and the associated mollified moments, with uniform control of the oldform and Eisenstein contributions as both parameters grow. We leave this problem for future work.
\end{remark}

\subsection{Mollified moments}
\label{subsec:intro-mollified-moments}
Fix a square-free positive integer $Q$ and $0<\eta<1/2$, independently of $p$. For primes $p\nmid Q$, put $L=p^\eta$ and define
\begin{equation}
x_\ell\colon \!\!= \mu(\ell)\prod_{q\mid\ell}\dfrac q{q+1}
 \dfrac{\log(L/\ell)}{\log L}\,
 \mathbf1_{\ell\leq L,\;(\ell,Q)=1}.
\label{eq:mollifier-coefficients}
\end{equation}
Here $\mu(n)$ denotes the M\"obius function. These coefficients satisfy $x_1=1$ and $|x_\ell|\leq1$. These $x_\ell$ vanish unless $\ell$ is square-free and $(\ell,Q)=1$. We use them to define the mollifier
\begin{equation}
M_\pi\colon \!\!= \prod_{q\mid Q}
 \left(1-\dfrac{A_\pi(1,q)}{\sqrt q}
 +\dfrac{A_\pi(q,1)}q-\dfrac1{q^{3/2}}\right)
 \sum_{\ell\leq L}\dfrac{x_\ell A_\pi(1,\ell)}{\sqrt\ell}.
\label{eq:intro-mollifier}
\end{equation}
The finite product consists of the inverse of Euler factors at primes dividing $Q$. Since $L<p$, all mollifier indices are prime to $p$. The mollified harmonic first and second moments are
\begin{equation}
\mathcal M_1(p)\colon \!\!= c_{\rm cusp}\sum_{\pi\in\mathcal B(p)}  \dfrac{h(\mu_\pi)L(1/2,\pi)}{L(1,\pi,\mathrm{Ad})}M_\pi,
\label{eq:intro-mollified-first-moment}
\end{equation}
and
\begin{equation}
\mathcal M_2(p)\colon \!\!= c_{\rm cusp}\sum_{\pi\in\mathcal B(p)}  \dfrac{h(\mu_\pi)|L(1/2,\pi)|^2}  {L(1,\pi,\mathrm{Ad})}|M_\pi|^2.
\label{eq:intro-mollified-second-moment}
\end{equation}

\begin{theorem}
\label{thm:intro-mollified-moments}
For every fixed square-free $Q$ and $0<\eta<1/2$, the mollifier in \eqref{eq:intro-mollifier} satisfies, as $p\to\infty$ through primes with $p\nmid Q$,
\begin{equation*}
\begin{split}
\mathcal M_1(p)=P\mathcal H_h+O_{h,Q,\eta}(p^{3/2}),
\end{split}
\end{equation*}
and
\begin{equation*}
\begin{split}
\mathcal M_2(p)=\zeta(2)\prod_{q\mid Q}(1-q^{-2})  (1+\eta^{-1})P\mathcal H_h  +O_{h,Q,\eta}\!\left(
 \dfrac{p^2\log(2+\log p)}{\log p}\right).
\end{split}
\end{equation*}
The test function and the parameters $Q,\eta$ are fixed before $p$ tends to infinity.
The finite inverse Euler factors replace the arithmetic factor $\zeta(2)$ by $\zeta(2)\prod_{q\mid Q}(1-q^{-2})$, which approaches one as $Q$ includes sufficiently many fixed small primes.
\end{theorem}

To improve the lower bound obtained from these two moments alone, we use a construction with two pieces inspired by Michel and VanderKam \cite{MV} and define
\begin{equation}
\widehat M_\pi:= M_\pi+
 \overline{\epsilon(\pi)}
 \frac{\gamma_{\mu_\pi}(1/2)}{\gamma_{-\mu_\pi}(1/2)}
 \overline{M_\pi}.
\label{eq:two-piece-mollifier}
\end{equation}
Here $\epsilon(\pi)$ and $\gamma_\mu$ are normalized as in \eqref{eq:functional-equation} and \eqref{eq:standard-gamma-factor}.

The functional equation gives $L(1/2,\pi)\widehat M_\pi=2\Re(L(1/2,\pi)M_\pi)$. The second moment for the mollifier with two pieces also involves the second moment without complex conjugation, which supplies its mixed term,
\begin{equation}
\mathcal M_{1,1}(p)\colon \!\!= c_{\rm cusp}\sum_{\pi\in\mathcal B(p)}  \dfrac{h(\mu_\pi)L(1/2,\pi)^2}{L(1,\pi,\mathrm{Ad})}M_\pi^2.
\label{eq:intro-mollified-mixed-moment}
\end{equation}
We refer to $\mathcal M_{1,1}(p)$ as the mixed moment below.

\begin{theorem}
\label{thm:intro-mollified-holomorphic-square}
For every fixed square-free $Q$ and $0<\eta<1/2$, as $p\to\infty$ through primes, we have
\begin{equation}
\mathcal M_{1,1}(p)= P\mathcal H_h \prod_{q\nmid Q}
 \left(1-\dfrac{2}{\sqrt q(q+1)^2}
             +\dfrac{1}{q(q+1)^2}\right)
 +O_{h,Q,\eta}\left(\dfrac{p^2}{\log p}\right).
\label{eq:intro-mollified-holomorphic-square}
\end{equation}
The Euler product converges absolutely and lies in $(0,1]$.
\end{theorem}

The role of the construction with two pieces is visible in the moment constants. As $Q$ ranges through products of increasingly many small primes, the arithmetic factor in Theorem~\ref{thm:intro-mollified-moments} and the Euler product in Theorem~\ref{thm:intro-mollified-holomorphic-square} both tend to one. In this limit, the constants in the lower bounds from the mollifiers with one and two pieces are governed, respectively, by
\[
\frac{1}{1+\eta^{-1}}\longrightarrow\frac13,
\qquad
\frac{2}{(1+\eta^{-1})+1}\longrightarrow\frac12,
\qquad \eta\to(1/2)^-.
\]
The second expression uses the mixed moment through \eqref{eq:two-piece-moment-identities}. For each prescribed error, we first fix $Q$, then fix $\eta<1/2$, and finally take $p$ sufficiently large. Together with the total harmonic mass, these three moments prove Theorem~\ref{thm:intro-weighted-nonvanishing}, as shown in \S\,\ref{sec:harmonic-nonvanishing}.

\subsection{Difficulties in the level aspect}
First, ramification at $p$ changes the dual Fourier indices produced by the functional equation. These terms require local coefficient identities at $p$, beyond the unramified Hecke relations. We use the local description in \cite[Lemma~1]{BC} to put them in a form suitable for the Kuznetsov formula.

Second, the Kuznetsov formula also contains oldforms and Eisenstein series, whose contributions must be controlled uniformly in $p$. We sum the three oldvectors before taking absolute values to
preserve cancellation in the inverse Gram matrix. The maximal Eisenstein new contribution is controlled by harmonic coefficient averages on $\GL_2$. 

Finally, the mollifier length $L=p^\eta$ makes the Fourier twists grow with $p$. We need estimates for the first moment that are uniform for $r\leq Qp^\eta$ and $s\leq Q$, and estimates for the second moment that are uniform for $ab\leq Q^2p^{2\eta}$, with all indices prime to $p$. These ranges must accommodate every fixed $0<\eta<1/2$ to obtain the proportion $1/2-\eps$.
\subsection{Outline of the proof}
We first prove the uniform twisted moment formulas in Theorems~\ref{thm:harmonic-first-moment} and \ref{thm:harmonic-second-moment} by the Kuznetsov formula at prime level.

We then insert the Hecke expansion of the mollifier to obtain $\mathcal M_1(p)$ and $\mathcal M_2(p)$. Retaining the full factor $\log(p/(ab))$ in the diagonal contribution to the second moment produces both quadratic forms in \eqref{eq:quadratic-forms}. Their evaluation gives the main term in Theorem~\ref{thm:intro-mollified-moments}.

The mixed moment is evaluated separately using the asymmetric approximate functional equation of Proposition~\ref{prop:holomorphic-square-afe}. After combining the Fourier coefficients, we apply the Kuznetsov formula. The first diagonal gives the Euler product in Theorem~\ref{thm:intro-mollified-holomorphic-square}, while the dual diagonal is absent. The kernel bounds and spectral estimates control the remaining terms.

Finally, the central functional equation combines the three moments through \eqref{eq:two-piece-moment-identities}. The Cauchy--Schwarz inequality and the total harmonic mass asymptotic give the required non-vanishing proportion.

\subsection{Plan of this paper}
\S\,\ref{prelim} fixes the coefficient conventions and proves the approximate functional equations and weight estimates. \S\,\ref{sec:kuznetsov} gives the Kuznetsov formula at prime level and its normalizations. \S\S\,\ref{sec:harmonic-twisted-moments} and \ref{sec:mollified-moments} evaluate the twisted and mollified moments, respectively. \S\,\ref{sec:harmonic-nonvanishing} completes the proof of Theorem \ref{thm:intro-weighted-nonvanishing}.

\subsection{Notation}

Throughout the paper, $\varepsilon$ and $A$ denote arbitrarily
small and arbitrarily large positive real numbers, respectively.
Their values may differ at different occurrences.
The Mellin contours and the remaining auxiliary parameters are then fixed before $p$ varies. 
In a truncation
estimate, the positive exponent is fixed before the level varies
and before the Mellin line for the tail is selected.
We write $e(z)=\exp(2\pi\i z)$, $f\ll g$ or $f=O(g)$ for
$|f|\leq Cg$ with an implied positive constant, and $f\asymp g$
when both $f\ll g$ and $g\ll f$ hold. The relation $f\sim g$
means $f-g=o(g)$. 
The symbols $(m,n)$, $\mu(n)$ and $\tau(n)$ denote the greatest
common divisor, the M\"obius function and the divisor function,
respectively. In products indexed by $q$, the index runs over
primes. For a prime $q$, $\operatorname{ord}_q n$
denotes the exponent of $q$ in $n$.
A superscript $(p)$ on an $L$-function indicates that its
Euler factor at $p$ is omitted.
We use the Euclidean norm $\|\mu\|=(\sum_j|\mu_j|^2)^{1/2}$.
\section{Preliminaries}
\label{prelim}
\subsection{$\rm GL_3$ $L$-functions}
\label{sec:preliminaries}

We work with the notations $\Gamma_0(p)$, $P=p^2+p+1$ and $\mathcal{B}(p)$ from the introduction, with $p$ an odd prime.
The finite Fourier coefficients and local bases follow
Blomer--Comtat \cite{BC}. The archimedean parameters use the convention
of \cite[\S\,1.3]{BZ}, with standard gamma factor
$\prod_j\Gamma_{\mathbb R}(s-\mu_j)$.

\subsubsection{Fourier coefficients}
 For every prime $q\neq p$, the Fourier coefficient $A_\pi(q^k,q^\ell)$ is given by
\begin{equation}
\det\begin{pmatrix}
\alpha_q^{k+\ell+2}&\beta_q^{k+\ell+2}&\gamma_q^{k+\ell+2}\\
\alpha_q^{k+1}&\beta_q^{k+1}&\gamma_q^{k+1}\\1&1&1
\end{pmatrix}/
\det\begin{pmatrix}\alpha_q^2&\beta_q^2&\gamma_q^2\\
\alpha_q&\beta_q&\gamma_q\\1&1&1\end{pmatrix},
\qquad k,\ell\geq0,
\label{eq:local-Schur}
\end{equation}
in terms of the Satake parameters $\alpha_q$, $\beta_q$, $\gamma_q$ of $\pi_q$.
At coinciding parameters, the quotient is interpreted
by its Schur polynomial continuation.
Since $\pi$ is unitary with trivial central character, we have
\begin{equation*}
\alpha_q\beta_q\gamma_q=1 \text{ and } \{\alpha_q^{-1},\beta_q^{-1},\gamma_q^{-1}\} =\{\overline{\alpha_q},\overline{\beta_q},\overline{\gamma_q}\}.
\end{equation*}
The bounds towards the Ramanujan conjecture of
Kim--Sarnak \cite{KS} and Blomer--Brumley \cite[Theorem~1]{BB} give
\begin{equation*}
\max\{|\alpha_q|,|\beta_q|,|\gamma_q|\}\leq  q^{\theta_3},
\qquad \theta_3=\dfrac{5}{14}.
\end{equation*}

For $q=p$,  we know $\pi_p\simeq |\cdot|_p^{-2\rho}\times  |\cdot|_p^\rho\mathrm{St}_2$ from the proof of \cite[Lemma 1]{BC}. Furthermore, we have the following lemma.
\begin{lemma}\label{lem:local-prime-level}
Let $\pi\in\mathcal{B}(p)$. Then there is $\rho\in \i\mathbb{R}$ such that
the parameters of its standard local $L$-factor are
\begin{equation*}
(\alpha_p,\beta_p,\gamma_p)=  (p^{-1/2-\rho},p^{2\rho},0).
\end{equation*}
\end{lemma}
\begin{proof}
See \cite[Lemma 1]{BC}. 
\end{proof}

The Fourier coefficients $A_\pi(p^k,p^\ell)$ are still given by
\eqref{eq:local-Schur} with these local $L$-parameters.
\begin{lemma}\label{ramified-coefficient}
For $k,\ell\geq0$ we have
\begin{equation}
A_\pi(p^k,p^\ell)= 
(\alpha_p\beta_p)^k\sum_{j=0}^{\ell}\alpha_p^j\beta_p^{\ell-j}
=p^{k(-1/2+\rho)}A_\pi(1,p^\ell).
\label{eq:ramified-Hecke}
\end{equation}
In particular, $A_\pi(p,1)=p^{-1/2+\rho}$.

\end{lemma}
\begin{proof}
Setting $\gamma_p=0$ in \eqref{eq:local-Schur} gives $A_\pi(p^k,p^\ell)=(\alpha_p\beta_p)^k (\alpha_p^{\ell+1}-\beta_p^{\ell+1})/(\alpha_p-\beta_p)$.
This proves \eqref{eq:ramified-Hecke}. 
\end{proof}
We extend these local coefficients multiplicatively to all positive integers $m,n$.
\begin{equation}
A_{\widetilde{\pi}}(m,n)= \overline{A_\pi(m,n)}.
\label{eq:contragredient-coeff}
\end{equation}
In particular,
\begin{equation}
A_\pi(m,n)=\overline{A_\pi(n,m)},\qquad (mn,p)=1.
\label{eq:unramified-unitarity}
\end{equation}
We will use, uniformly in $\pi$,
\begin{equation*}
|A_\pi(m,n)|\ll_\varepsilon (mn)^{\theta_3+\varepsilon}.
\end{equation*}
At $p$ this follows with exponent zero from \eqref{eq:ramified-Hecke}.
At unramified primes, apply the preceding Ramanujan bound
to the Schur polynomial in \eqref{eq:local-Schur} and use the divisor
estimate.

\subsubsection{The standard $L$-function and its root number}
For $\Re(s)>1$,
\begin{equation*}
L(s,\pi)= \sum_{n\geq1}\dfrac{A_\pi(1,n)}{n^s}=\prod_q\prod_{\xi\in\{\alpha_q,\beta_q,\gamma_q\}} (1-\xi q^{-s})^{-1},
\end{equation*}
The contragredient $L$-function $L(s,\widetilde\pi)$ is defined analogously.
Let $\mu=\mu_\pi$ denote the parameters in the Kuznetsov convention, so $\sum_j\mu_j=0$, $|\Re(\mu_j)|\leq\theta_3$ and $\{\overline\mu_j\}=\{-\mu_j\}$. Define
\begin{equation}
\gamma_\mu(s)\colon = \prod_{j=1}^3\Gamma_{\mathbb{R}}(s-\mu_j),
\qquad \Gamma_{\mathbb{R}}(s)=\pi^{-s/2}\Gamma(s/2).
\label{eq:standard-gamma-factor}
\end{equation}
Then the completed $L$-function is 
\begin{equation*}
\Lambda(s,\pi):= p^{s/2}\gamma_\mu(s)L(s,\pi),
\end{equation*}
which is entire and satisfies the functional equation
\begin{equation}
\Lambda(s,\pi)= \epsilon(\pi)\Lambda(1-s,\widetilde\pi),
\qquad \mu_{\widetilde\pi}=-\mu_\pi.
\label{eq:functional-equation}
\end{equation}
The conductor at height $t$ is
\begin{equation*}
C(\pi,t)\asymp p\prod_{j=1}^3(1+|t-\Im(\mu_j)|).
\end{equation*}

Fix the standard additive character $\psi=\prod_v\psi_v$ of $\mathbb{Q}\backslash\mathbb{A}$ and the corresponding self-dual Haar measures. The global root number is
\begin{equation*}
\epsilon(\pi)=\prod_v\epsilon\!\left(\dfrac{1}{2},\pi_v,\psi_v\right).
\end{equation*}

\begin{lemma}[Root number at prime conductor]
\label{lem:root-number}
We have
\begin{gather}
\epsilon(s,\pi_p,\psi_p)     =-p^{-\rho-(s-1/2)},
\label{eq:local-epsilon}\\
\epsilon(\pi)     =-p^{-\rho}   =-p^{1/2}\overline{A_\pi(p,1)},
\label{eq:global-root}\\
\epsilon(\pi)A_{\widetilde\pi}(1,n)     =-p^{1/2}\overline{A_\pi(p,n)}       \qquad(n\geq1).
\label{eq:coefficient-identity}
\end{gather}
In particular, \eqref{eq:coefficient-identity} also holds when $p\mid n$.
\end{lemma}
\begin{proof}
Let $D:=\text{WD}(\pi_p)$ be the associated Weil--Deligne representation. Explicitly, $D=(r,N)$ acts on $V=\langle e_0,e_1,e_2\rangle$ by
\begin{equation*}
\begin{aligned}
r(w)e_0&=|w|^{-2\rho}e_0,   & Ne_0&=0,\\
  r(w)e_1&=|w|^{\rho-1/2}e_1,   & Ne_1&=e_2,\\
  r(w)e_2&=|w|^{\rho+1/2}e_2,  & Ne_2&=0.
\end{aligned}
\end{equation*}
Let $I_p$ be the inertia subgroup and choose a geometric Frobenius $\Phi_p$ with $|\Phi_p|=p^{-1}$. Since the Weil representation $r$ is unramified, we have
\begin{equation*}
V^{I_p}=V,   \qquad   (\ker N)^{I_p}=\langle e_0,e_2\rangle,   \qquad   V^{I_p}/(\ker N)^{I_p}=\langle e_1 \rangle.
\end{equation*}
The Frobenius $\Phi_p$ acts on $\langle e_1 \rangle$ by \begin{equation*}
r(\Phi_p)e_1=p^{1/2-\rho} e_1.
\end{equation*}
The conductor formula gives
\begin{equation*}
a(D)  =a(r)+\dim V^{I_p}-\dim(\ker N)^{I_p}  =0+3-2=1.
\end{equation*}

Hence Deligne's monodromy correction formula is
\begin{equation*}
  \epsilon(s,D,\psi_p)  =\epsilon(s,r,\psi_p)    \det\!\left(
        -p^{-s}r(\Phi_p)    \,\middle|\,
        V^{I_p}/(\ker N)^{I_p}
      \right),
\end{equation*}
see \cite[\S\,8.12]{Deligne} and \cite[\S\,4.1.6, pp. 20--21]{Tate}. Since $r$ is a direct sum of unramified characters and $\psi_p$ has conductor $\mathbb{Z}_p$, we have $\epsilon(s,r,\psi_p)=1$. Consequently,
\begin{equation*}
\epsilon(s,\pi_p,\psi_p)    =-p^{-\rho-(s-1/2)}
\end{equation*}
which proves \eqref{eq:local-epsilon}.

At each finite prime $\ell\neq p$, both $\pi_\ell$ and $\psi_\ell$ are unramified, and hence $\epsilon(1/2,\pi_\ell,\psi_\ell)=1$. At infinity, the trivial central character and $\text{O}_3=\text{SO}_3\cdot\{\pm I_3\}$ imply that an $\text{SO}_3$-fixed vector is also $\text{O}_3$-fixed. The archimedean Langlands parameter is therefore a sum of even characters of $\mathbb{R}^\times$, each of which has epsilon factor $1$. See \cite[Eqs.\,(3.7)--(3.8), p. 404]{Knapp}. It follows that
\begin{equation}
\label{eq:root-rho}
  \epsilon(\pi)
  =\epsilon\!\left(\dfrac12,\pi_p,\psi_p\right)
  =-p^{-\rho}.
\end{equation}

 The coefficients $A_\pi(1,n)$ are the Dirichlet coefficients of the standard $L$-function, so the compatible normalization for
$\widetilde\pi$ gives
\begin{equation*}
A_{\widetilde\pi}(1,n)
  =\overline{A_\pi(1,n)}
  \qquad(n\geq1).
\end{equation*}
The \eqref{eq:root-rho} and Lemma \ref{ramified-coefficient} imply
\begin{equation*}
-p^{1/2}\overline{A_\pi(p,1)}     =-p^{1/2}p^{-1/2-\rho}     =-p^{-\rho}
    =\epsilon(\pi),
\end{equation*}
which proves \eqref{eq:global-root}. 
Using \eqref{eq:ramified-Hecke} at $p$, multiplicativity at the other
primes and \eqref{eq:contragredient-coeff}, we obtain
\begin{equation*}
\begin{aligned}
-p^{1/2}\overline{A_\pi(p,n)}    &=-p^{1/2}\overline{A_\pi(p,1)} \,\overline{A_\pi(1,n)}\\
   &=\epsilon(\pi)A_{\widetilde{\pi}}(1,n).
\end{aligned}
\end{equation*}
This proves \eqref{eq:coefficient-identity} for every $n\geq1$.
\end{proof}

\subsection{Stirling's formula and bounds for the zeta function}
\label{subsec:stirling}

Stirling's formula \cite[\S\,8.4, Eq.\,(4.03)]{O} gives,
for $|\arg(s)|\leq\pi-\varepsilon$, $|s|\gg1$,
any $\varepsilon>0$ and any integer $K_1\geq0$,
\begin{equation*}
\log \Gamma(s)=\left(s-\frac{1}{2}\right)\log s-s+\frac{1}{2}\log(2\pi)
+\sum_{j=1}^{K_1}\frac{B_{2j}}{2j(2j-1)s^{2j-1}}+O_{K_1,\varepsilon}\left(\frac{1}{|s|^{2K_1+1}}\right),
\end{equation*}
where $B_j$ is the $j$-th Bernoulli number. Thus, for $s=\sigma+\i t$ with fixed 
$\sigma$ and $|t|\geq1$, we have
\begin{equation*}
|\Gamma(\sigma+\i t)|= \sqrt{2\pi}\,
|t|^{\sigma-1/2}e^{-\pi|t|/2}
\big(1+O(|t|^{-1})\big),
\end{equation*}
and\begin{equation*}
\frac{\Gamma'}{\Gamma}(\sigma+\i t)= 
\log(\sigma+\i t)-\frac{1}{2(\sigma+\i t)}+O(|t|^{-2}).
\end{equation*}
Consequently,
\begin{equation*}
|\Gamma_{\mathbb R}(\sigma+\i t)|\asymp 
(1+|t|)^{\sigma/2-1/2}e^{-\pi|t|/4}
\qquad (|t|\geq1).
\end{equation*}

Fix $a\geq0$, $c_0\leq c_1$ and $\delta>0$.
For $c_0\leq\Re(u)\leq c_1$, assume that each argument
$1/2+u\mp\mu_j$ stays at distance at least $\delta$ from
$\{0,-2,-4,\ldots\}$. Stirling's formula and the triangle inequality give
\begin{equation}
\left|\frac{\gamma_{\pm\mu}(1/2+u)}{\gamma_\mu(1/2)}\right|
\ll_{a,c_0,c_1,\delta} (1+\|\mu\|+|\Im(u)|)^C e^{3\pi|\Im(u)|/4},
\qquad |\Re(\mu_j)|\leq a,
\label{eq:gamma-ratio-bound}
\end{equation}
where $C=C(a,c_0,c_1)$, and $1/\gamma_\mu(1/2)$ is interpreted
by its entire continuation.

We also record two bounds for the zeta function.
The integral formula in \cite[\S\,4, Eq.\,(7), pp.~32--33]{Davenport} gives
\begin{equation}
|(s-1)\zeta(s)|\leq |s|+\frac{|s(s-1)|}{\Re(s)}
\ll(1+|s|)^2,\qquad \Re(s)\geq1/2,
\label{eq:zeta-pole-removed}
\end{equation}
with continuous extension at $s=1$.
The convexity estimate \cite[\S\,5.2, Eqs.\,(5.19)--(5.20), p.\,100]{IK} gives
\begin{equation}
\zeta(5/8+\i t)\ll_\eps (1+|t|)^{3/16+\eps},
\qquad t\in\mathbb R.
\label{eq:zeta-five-eighths-bound}
\end{equation}

\subsection{Approximate functional equations}
\label{subsec:afe}
Fix an integer $J\geq2$ and put
\begin{equation}
P_J(u,\mu):= \prod_{j=1}^3\prod_{\ell=0}^{J}\prod_{\epsilon=\pm1}
 \big((\tfrac12+2\ell+\epsilon\mu_j)^2-u^2\big),
\label{eq:afe-polynomial}
\end{equation}
and
\begin{equation}
G_J(u,\mu):= e^{u^2}\frac{P_J(u,\mu)}{P_J(0,\mu)}.
\label{eq:afe-auxiliary-function}
\end{equation}
On the generic spherical unitary spectrum occurring in Theorem~\ref{thm:kuznetsov}, the parameters satisfy $|\Re(\mu_j)|\leq5/14<1/2$. Hence $P_J(0,\mu)\neq0$ on this spectrum. The function $G_J$ is even in $u$,
invariant under $\mu\mapsto-\mu$ and under the Weyl group, and
$G_J(0,\mu)=1$.

\subsubsection{The standard $L$-function}

Write $U_\mu^+=V_\mu$ and $U_\mu^-=W_\mu$. For $y>0$ and $c>1$, we define
\begin{equation*}
U_\mu^\pm(y):= \frac1{2\pi \i}\int_{(c)}y^{-u}G_J(u,\mu)
 \frac{\gamma_{\pm\mu}(1/2+u)}{\gamma_\mu(1/2)}\frac{\dd u}{u}.
\end{equation*}
%The denominator in both integrals is $\gamma_\mu(1/2)$.
%In particular the dual weight includes the archimedean factor
%$\gamma_{-\mu}(1/2)/\gamma_\mu(1/2)$.

\begin{proposition}
\label{prop:AFE}
For $\pi\in\mathcal B(p)$, the central value
$L(1/2,\pi)$ satisfies the following approximate functional equation
\begin{equation}
L(1/2,\pi)= \sum_{n\geq1}\frac{A_\pi(1,n)}{\sqrt n}
 V_{\mu_\pi}(n/\sqrt p)
-\sqrt p\sum_{n\geq1}\frac{\overline{A_\pi(p,n)}}{\sqrt n}
 W_{\mu_\pi}(n/\sqrt p).
\label{eq:AFE-final}
\end{equation}
\end{proposition}

\begin{proof}
For $\mu=\mu_\pi$ and $c>1$, consider
\begin{equation*}
\frac1{2\pi\i}\int_{(c)}
\frac{\Lambda(1/2+u,\pi)}{p^{1/4}\gamma_\mu(1/2)}
G_J(u,\mu)\frac{\dd u}{u}.
\end{equation*}
Moving the contour to $\Re(u)=-c$ gives the residue $L(1/2,\pi)$.
The functional equation \eqref{eq:functional-equation} and
\eqref{eq:coefficient-identity} then yield \eqref{eq:AFE-final}
after expanding the Dirichlet series.
\end{proof}
We record the following estimates for $U_\mu^\pm(y)$.
\begin{lemma}
\label{lem:afe-weights}
Let $J\geq2$, $A>0$ and $j\in\mathbb Z_{\geq0}$.
On the unitary automorphic spectrum with $|\Re(\mu_j)|\leq5/14$, for every $y>0$,
\begin{equation}
\left|y^j\frac{\dd^j}{\dd y^j}U_\mu^\pm(y)\right|
\ll_{A,j,J}(1+\|\mu\|)^C(1+y)^{-A},
\label{eq:V-decay}
\end{equation}
where $C=C(A,j,J)$. For $0<y\leq1$,
\begin{equation}
U_\mu^\pm(y)=\frac{\gamma_{\pm\mu}(1/2)}{\gamma_\mu(1/2)}
 +O_J\bigl(y^2(1+\|\mu\|)^C\bigr),
\label{eq:V-small}
\end{equation}
with a suitable $C=C(J)$. The ratio in the main term has modulus one.
\end{lemma}

\begin{proof}
Since $|\Re(\mu_k)|\leq5/14$, $P_J(0,\mu)$ is bounded away
from zero, and the factors in \eqref{eq:afe-polynomial} cancel
the numerator gamma poles in $\Re(u)\geq-2$.
The gamma recurrence and \eqref{eq:gamma-ratio-bound} give,
for $u=\sigma+\i v$ with $-2\leq\sigma\leq A+2$,
\begin{equation*}
\left|G_J(u,\mu)
 \frac{\gamma_{\pm\mu}(1/2+u)}{\gamma_\mu(1/2)}\right|
\ll_{A,J}(1+\|\mu\|+|v|)^C e^{-v^2+3\pi|v|/4}.
\end{equation*}
This estimate justifies the contour shifts and differentiation
under the integral below.

For $0<y\leq1$, move to $\Re(u)=-2$.
The residue at zero is $\gamma_{\pm\mu}(1/2)/\gamma_\mu(1/2)$,
which gives \eqref{eq:V-small}.
For $j\geq1$, differentiating $y^{-u}$ introduces the factor
$(-1)^ju(u+1)\cdots(u+j-1)$ after multiplication by $y^j$.
This cancels the pole at zero, so the same shift gives
\begin{equation*}
\left|y^j\frac{\dd^j}{\dd y^j}U_\mu^\pm(y)\right|
\ll_{j,J}y^2(1+\|\mu\|)^C.
\end{equation*}
Unitarity gives $\gamma_{-\mu}(1/2)=\overline{\gamma_\mu(1/2)}$,
so the ratio in the main term has modulus one. This proves
\eqref{eq:V-decay} for $y\leq1$.
For $y\geq1$, moving to $\Re(u)=A+1$ bounds the differentiated
integral by $O_{A,j,J}(y^{-A-1}(1+\|\mu\|)^C)$,
completing the proof.
\end{proof}

\subsubsection{The product of the standard and dual $L$-functions}
For the second moment, all mollified moments, and non-vanishing,
we use the fixed test function
\begin{equation}
h(\mu):= \exp\!\left(\sum_{j=1}^3\mu_j^2\right)
 \prod_{j=1}^3\cos^4(\pi\mu_j)
 \prod_{1\leq i<j\leq3}
 \cos^2\!\left(\frac{\pi(\mu_i-\mu_j)}2\right).
\label{eq:positive-moment-test}
\end{equation}
The test $h$ is independent of $p$, $\eta$ and $J$.
It is even, Weyl invariant and rapidly decreasing on every finite
tube. Unitarity makes the multiset of parameters invariant under $\mu_j\mapsto-\overline{\mu_j}$, so $\sum_j\mu_j^2$ is real and its exponential is positive. The coordinate cosine factors pair into absolute-value powers. The difference factors occur in conjugate pairs or as squares of real cosines. Thus $h$ is non-negative on the unitary spectrum. On $\Re(\mu)=0$, all cosine factors are positive hyperbolic cosines, so $h>0$.
Its zeros on $\mu_j\in\tfrac12+\mathbb Z$ and
$\mu_i-\mu_j\in2\mathbb Z+1$ have orders at least four and two,
respectively.
For $a_0\geq10$ and an integer $J\geq\lceil a_0\rceil+4$,
the quotient
\begin{equation}
h_0(\mu):=\frac{h(\mu)}{P_J(0,\mu)}
\label{eq:moment-test-quotient}
\end{equation}
belongs to $\widetilde{\mathcal H}(a_0+1)$.
The test class, the tubes $\Lambda_a$ and the
coordinate derivatives used below are defined in
\S\,\ref{subsec:kuz-test}.

For $y>0$ and $c>1$, we define
\begin{equation}
\mathcal U_\mu(y):=\frac1{2\pi\i}\int_{(c)}
 y^{-u}G_J(u,\mu)
\frac{\gamma_\mu(1/2+u)\gamma_{-\mu}(1/2+u)}
 {\gamma_\mu(1/2)\gamma_{-\mu}(1/2)}\frac{\dd u}{u}.
\label{eq:product-weight-mellin}
\end{equation}

\begin{proposition}
\label{prop:product-afe}
For $\pi\in\mathcal B(p)$, the absolute square
$|L(1/2,\pi)|^2$ satisfies the following approximate functional equation
\begin{equation}
|L(1/2,\pi)|^2= 2\sum_{m,n\geq1}
 \frac{\overline{A_\pi(1,m)}A_\pi(1,n)}{\sqrt{mn}}
 \mathcal U_{\mu_\pi}(mn/p).
\label{eq:product-afe}
\end{equation}
\end{proposition}

\begin{proof}
For $\mu=\mu_\pi$ and $c>1$, consider
\begin{equation*}
\dfrac1{2\pi\i}\int_{(c)}  \dfrac{\Lambda(1/2+u,\pi)\Lambda(1/2+u,\widetilde\pi)}  {p^{1/2}\gamma_\mu(1/2)\gamma_{-\mu}(1/2)}  G_J(u,\mu)\frac{\dd u}{u}.
\end{equation*}
The same contour argument as in the proof of Proposition \ref{prop:AFE} gives \eqref{eq:product-afe}.
\end{proof}

We next estimate $\mathcal U_\mu(y)$ and the associated spectral test functions
\begin{equation*}
H_t(\mu) \colon\!\!=h(\mu)\mathcal U_\mu(t/p),\qquad\text{for}\,\,t\geq1.
\end{equation*}

\begin{lemma}
\label{lem:product-weights}
Fix $a_0\geq10$ and an integer $J\geq\lceil a_0\rceil+4$. For every fixed $B>0$, on the unitary automorphic spectrum with $|\Re(\mu_j)|\leq5/14$,
\begin{equation}
|\mathcal U_\mu(y)|\ll_{B,J} 
 (1+\|\mu\|)^C(1+y)^{-B}
\label{eq:product-weight-decay}
\end{equation}
for a suitable $C$. For every $c>0$, $M\geq0$ and integer $k\geq0$,
$H_t$ belongs to
$\widetilde{\mathcal H}(a_0)$ and satisfies
\begin{equation}
\sup_{\substack{\mu\in\Lambda_{a_0}\\|\alpha|\leq k}}
(1+\|\mu\|)^M|\partial_\mu^\alpha H_t(\mu)|
\ll_{c,M,k,h,J} (p/t)^c.
\label{eq:product-test-seminorm}
\end{equation}
At zeros of $P_J(0,\mu)$, $H_t$ is understood by
holomorphic continuation.
\end{lemma}

\begin{proof}
Write
\begin{equation*}
R_\mu(u)\colon\!\!=\dfrac{\gamma_\mu(1/2+u)\gamma_{-\mu}(1/2+u)}
 {\gamma_\mu(1/2)\gamma_{-\mu}(1/2)}.
\end{equation*}
The factors of $P_J$ in \eqref{eq:afe-polynomial} cancel the numerator gamma poles crossed below. The gamma recurrence and \eqref{eq:gamma-ratio-bound} justify the contour shifts in 
\eqref{eq:product-weight-mellin}.
For $y\geq1$, move to $\Re(u)=B$. For $0<y\leq1$, move to $\Re(u)=-1/4$, crossing only the pole at zero with residue $1$.
Thus
\begin{equation*}
\mathcal U_\mu(y)=1+O_J\bigl(y^{1/4}(1+\|\mu\|)^C\bigr),
\qquad 0<y\leq1,
\end{equation*}
and \eqref{eq:product-weight-decay} follows.

For the test estimate, \eqref{eq:afe-auxiliary-function} gives
\begin{equation*}
h(\mu)G_J(u,\mu)R_\mu(u) =\dfrac{h(\mu)}{P_J(0,\mu)}e^{u^2}P_J(u,\mu)R_\mu(u).
\end{equation*}
The quotient in \eqref{eq:moment-test-quotient} is holomorphic and rapidly decreasing on $\Lambda_{a_0+1}$. For fixed $\Re(u)=c>0$, the gamma recurrence makes the displayed multiplier holomorphic on this tube and preserves Weyl invariance and the zeros in \eqref{eq:test-zero-Kuz} below.

Applying \eqref{eq:gamma-ratio-bound} and Cauchy's integral formula on this larger tube gives, for $u=c+\i v$, 
\begin{equation*}
\sup_{\substack{\mu\in\Lambda_{a_0}\\|\alpha|\leq k}}
(1+\|\mu\|)^M
\left|\partial_\mu^\alpha
\bigl(h(\mu)G_J(u,\mu)R_\mu(u)\bigr)\right|
\ll_{a_0,c,M,k,h,J}(1+|v|)^C e^{-v^2+3\pi|v|/2}.
\end{equation*}
Integrating on $\Re(u)=c$ with $y=t/p$ now gives \eqref{eq:product-test-seminorm} and $H_t\in\widetilde{\mathcal H}(a_0)$.
\end{proof}

\subsubsection{The square of the standard $L$-function}

For $y>0$ and $c>1$, we define
\begin{equation}
\mathcal V_\mu^\pm(y) \colon \!\!= \dfrac1{2\pi\i}\int_{(c)}  y^{-u}G_J(u,\mu)^2  \left(\dfrac{\gamma_{\pm\mu}(1/2+u)}
 {\gamma_\mu(1/2)}\right)^2\dfrac{\dd u}{u}.
\label{eq:holomorphic-square-weight}
\end{equation}

\begin{proposition}
\label{prop:holomorphic-square-afe}
For $\pi\in\mathcal B(p)$ and $L>0$, the square $L(1/2,\pi)^2$ satisfies the following asymmetric approximate functional equation
\begin{equation}
\begin{split}
L(1/2,\pi)^2= &~\sum_{m,n\geq1}
 \dfrac{A_\pi(1,m)A_\pi(1,n)}{\sqrt{mn}}
 \mathcal V_{\mu_\pi}^+\left(\frac{mn}{pL^2}\right)\\
&~+\epsilon(\pi)^2\sum_{m,n\geq1}
 \dfrac{\overline{A_\pi(1,m)A_\pi(1,n)}}{\sqrt{mn}}
 \mathcal V_{\mu_\pi}^-\left(\frac{mnL^2}{p}\right).
\end{split}
\label{eq:asymmetric-holomorphic-square-afe}
\end{equation}
\end{proposition}

\begin{proof}
For $\mu=\mu_\pi$ and $c>1$, consider
\begin{equation*}
\dfrac1{2\pi\i}\int_{(c)}
 \dfrac{\Lambda(1/2+u,\pi)^2}
 {p^{1/2}\gamma_\mu(1/2)^2}
 G_J(u,\mu)^2L^{2u}\dfrac{\dd u}{u}.
\end{equation*}
The same contour argument as in the proof of Proposition \ref{prop:AFE} gives \eqref{eq:asymmetric-holomorphic-square-afe}.
\end{proof}

For $y>0$, define the corresponding spectral test functions by
\begin{equation*}
K_\mu^\pm(y)\colon \!\!=h(\mu)\mathcal V_\mu^\pm(y),
\end{equation*}
where $h$ is the fixed test function in \eqref{eq:positive-moment-test}. At zeros of $P_J(0,\mu)$, these products are defined by holomorphic continuation, as justified in Lemma \ref{lem:holomorphic-square-weights}. The following bounds apply to $\mathcal V_\mu^\pm(y)$ and $K_\mu^\pm(y)$.

\begin{lemma}
\label{lem:holomorphic-square-weights}
On the unitary automorphic spectrum with $|\Re(\mu_j)|\leq5/14$, for every $B>0$ and some $C>0$, uniformly for $y>0$,
\begin{equation}
|\mathcal V_\mu^\pm(y)|\ll_{B,J} 
 (1+\|\mu\|)^C(1+y)^{-B}.
\label{eq:holomorphic-square-weight-decay}
\end{equation}
For $0<y\leq1$ and every fixed $0<c<1/4$,
\begin{equation}
\mathcal V_\mu^\pm(y)= 
 \left(\dfrac{\gamma_{\pm\mu}(1/2)}{\gamma_\mu(1/2)}\right)^2
 +O_{J,c}\bigl(y^c(1+\|\mu\|)^C\bigr).
\label{eq:holomorphic-square-weight-small}
\end{equation}
Let $a_0\geq10$, $J\geq\lceil a_0\rceil+4$, and let $h$ be the fixed test function in \eqref{eq:positive-moment-test}. For every $y>0$, the functions $\mu\mapsto K_\mu^\pm(y)$ belong to $\widetilde{\mathcal H}(a_0)$. For every fixed $c>0$, $M\geq0$ and integer $k\geq0$, they satisfy
\begin{equation}
\sup_{\substack{\mu\in\Lambda_{a_0}\\|\alpha|\leq k}}
(1+\|\mu\|)^M
\left|\partial_\mu^\alpha
K_\mu^\pm(y)\right|
\ll_{h,J,a_0,c,M,k} y^{-c}.
\label{eq:holomorphic-square-test-seminorm}
\end{equation}
\end{lemma}

\begin{proof}
By \eqref{eq:afe-polynomial}, $P_J(u,\mu)^2$ cancels the numerator gamma poles with their full multiplicities in $\Re(u)\geq-2$. The contour argument in the proof of Lemma \ref{lem:afe-weights} applies to the squared multiplier in \eqref{eq:holomorphic-square-weight}. Move to $\Re(u)=B+1$ for $y\geq1$ and to $\Re(u)=-c$ for $0<y\leq1$, where $0<c<1/4$. The latter shift crosses only the pole at zero, with residue $(\gamma_{\pm\mu}(1/2)/\gamma_\mu(1/2))^2$ of modulus one. This proves \eqref{eq:holomorphic-square-weight-decay} and \eqref{eq:holomorphic-square-weight-small}.

For the test estimate, write
\begin{equation*}
h(\mu)G_J(u,\mu)^2
 \left(\dfrac{\gamma_{\pm\mu}(1/2+u)}{\gamma_\mu(1/2)}\right)^2
=\dfrac{h(\mu)}{P_J(0,\mu)^2}e^{2u^2}P_J(u,\mu)^2
 \dfrac{\gamma_{\pm\mu}(1/2+u)^2}{\gamma_\mu(1/2)^2}.
\end{equation*}
The denominator $P_J(0,\mu)^2$ has zeros of order four on each of its half-integer coordinate hyperplanes. The factors $\cos^4(\pi\mu_j)$ in $h$ cancel these zeros with their full multiplicities. Thus $h(\mu)/P_J(0,\mu)^2$ is holomorphic and rapidly decreasing on finite tubes.
The gamma recurrence and $J\geq\lceil a_0\rceil+4$ make the displayed multiplier holomorphic on $\Lambda_{a_0+1}$ for $\Re(u)>0$ and preserve Weyl invariance and the zeros in \eqref{eq:test-zero-Kuz}.
Stirling's formula and Cauchy estimates in the proof of Lemma \ref{lem:product-weights} now apply with the decay factor $e^{-2v^2+3\pi|v|/2}$, where $u=c+\i v$.
Integrating on $\Re(u)=c$ proves \eqref{eq:holomorphic-square-test-seminorm}.
\end{proof}

\section{\texorpdfstring{The $\mathrm{GL}_3$ Kuznetsov formula at prime level}{The GL(3) Kuznetsov formula at prime level}}
\label{sec:kuznetsov}
We use the finite Fourier coefficients and local bases of Blomer--Comtat \cite[\S\S\,1.4, 4 and 5]{BC}, the archimedean parameter and kernel conventions of \cite[\S\S\,1.3 and 4]{BZ}.
The harmonic weights use full finite adjoint values. The geometric side is given in \cite[Theorem~6]{BBM}. For the archimedean estimates we use \cite[Lemmas~2.4--2.6]{BL}, in arXiv:2401.04037v2.

\subsection{Kuznetsov test functions}
\label{subsec:kuz-test}

Let $\Lambda_a=\{\mu\in\mathbb C^3\colon \mu_1+\mu_2+\mu_3=0, \: |\Re(\mu_j)|\leq a\}$. Following \cite[Definition~1]{BC}, $\widetilde{\mathcal H}(a)$ consists of the test functions that are Weyl group invariant, holomorphic on a neighborhood of $\Lambda_a$, rapidly decreasing there and vanishing on 
\begin{equation}
\mu_i-\mu_j= n,\qquad n\in2\mathbb Z+1,\quad |n|\leq a,
\quad 1\leq i<j\leq3.
\label{eq:test-zero-Kuz}
\end{equation}
Rapid decrease means that
$\sup_{\Lambda_a}(1+\|\mu\|)^M|h(\mu)|<\infty$ for every $M\geq0$.
We use $\mu_1,\mu_2$ as independent coordinates and set $\mu_3=-\mu_1-\mu_2$.

For $\alpha=(\alpha_1,\alpha_2)\in\mathbb Z_{\geq0}^2$, write
\begin{equation*}
|\alpha|=\alpha_1+\alpha_2,\qquad
\partial_\mu^\alpha f(\mu)\colon =
\dfrac{\partial^{\alpha_1+\alpha_2}f(\mu_1,\mu_2,-\mu_1-\mu_2)}
 {\partial\mu_1^{\alpha_1}\partial\mu_2^{\alpha_2}}.
\end{equation*}

The fixed test function for the moments is given in \eqref{eq:positive-moment-test}.

\subsection{Archimedean kernels and integral transforms}
\label{subsec:kuz-transforms}

On $\Re(\mu)=0$ the spectral density is \cite[\S~4]{BZ}
\begin{equation*}
\operatorname{spec}(\mu)\colon \!\!= \prod_{(i,j)\in\{(1,2),(2,3),(3,1)\}}  (\mu_i-\mu_j)\tan\!\left(\dfrac\pi2(\mu_i-\mu_j)\right).
\end{equation*}
For the archimedean transforms, write $\mu_j=\i t_j$, $t_3=-t_1-t_2$, and $\dd\mu=\dd\mu_1\,\dd\mu_2=-\dd t_1\,\dd t_2$. Then $\operatorname{spec}(\mu)\dd\mu$ is non-negative and is positive away from the Weyl chamber walls.
Since $|\tan(\i x)|=|\tanh x|\leq1$ for $x\in\mathbb R$,
\begin{equation*}
|\operatorname{spec}(\mu)|\ll (1+\|\mu\|)^3
\qquad(\Re(\mu)=0).
\end{equation*}
Consequently, for every fixed $C\geq0$ and rapidly decreasing $h$,
\begin{equation}
\int_{\Re(\mu)=0}|h(\mu)|(1+\|\mu\|)^C
\operatorname{spec}(\mu)\,\dd\mu\ll_{h,C} 1.
\label{eq:integrable-spectral-density}
\end{equation}
The following kernels are $256\pi^5$ times those in \cite[Eqs.\,(10), (12), (16)--(18)]{BZ}, consistently with the spectral constants and the diagonal term \eqref{eq:Delta-exact}.
For $\epsilon\in\{\pm1\}$, we set 
\[
\widetilde G^\epsilon(z,\mu)\colon\!\!=
\dfrac{\pi^{-3z}}{2}
\left\{
\prod_{j=1}^{3}
\dfrac{\Gamma((z-\mu_j)/2)}     {\Gamma((1-z+\mu_j)/2)}+\epsilon i
\prod_{j=1}^{3}
\dfrac{\Gamma((1+z-\mu_j)/2)} {\Gamma((2-z+\mu_j)/2)}
\right\}.
\]
For $\epsilon_1,\epsilon_2\in\{\pm1\}$ and $z=(z_1,z_2)$, we set
\[
\begin{aligned}
G_{\mathrm{sym}}^{\epsilon_1,\epsilon_2}(z,\mu)
\colon ={}& \dfrac{\pi^{5/2}}{4} \sum_{\substack{d_1,d_2,d_3\in\{0,1\}\\
                 d_3\equiv d_1+d_2\!\!\!\pmod 2}}
\epsilon_1^{d_1}\epsilon_2^{d_2}(-1)^{d_1d_2}
\dfrac{\Gamma((1+d_3-z_1-z_2)/2)}
     {\Gamma((d_3+z_1+z_2)/2)}
\\
&\qquad\times
\prod_{j=1}^{3}
\dfrac{
\Gamma((d_1+z_1-\mu_j)/2)
\Gamma((d_2+z_2+\mu_j)/2)
}{
\Gamma((1+d_1-z_1+\mu_j)/2)
\Gamma((1+d_2-z_2-\mu_j)/2)
}.
\end{aligned}
\]

The kernels are defined by the Barnes integrals
\begin{equation*}
\begin{aligned}
K_{w_4}(y;\mu)\colon \!\!= &~\int_{\mathscr C}|y|^{-z}
 \widetilde G^{\sgn y}(z,\mu)\frac{\dd z}{2\pi \i},\\
K_{w_6}^{\mathrm{sym}}((y_1,y_2);\mu)\colon\!\!= &~
\int_{\mathscr C_1}\int_{\mathscr C_2}
(\pi^2|y_1|)^{-z_1}(\pi^2|y_2|)^{-z_2}
G_{\mathrm{sym}}^{\sgn(y_1),\sgn(y_2)}(z,\mu)
\frac{\dd z_1\,\dd z_2}{(2\pi \i)^2}.
\end{aligned}
\end{equation*}
We use the Barnes contours of \cite[\S\,2.4]{BL}, separating the poles, avoiding poles on the imaginary axis, and deforming the tails for convergence.
For $y,y_1,y_2\neq0$ define 
\begin{equation}
\begin{split}
\Phi_{w_4,h}(y)\colon \!\!= &~\int_{\Re(\mu)=0}h(\mu)K_{w_4}(y;\mu)  \operatorname{spec}(\mu)\,\dd\mu,\\
\Phi_{w_5,h}(y)\colon \!\!= &~\int_{\Re(\mu)=0}h(\mu)K_{w_4}(-y;-\mu)  \operatorname{spec}(\mu)\,\dd\mu,\\
\Phi_{w_6,h}(y_1,y_2)\colon \!\!= &~\int_{\Re(\mu)=0}h(\mu) K_{w_6}^{\mathrm{sym}}((y_1,y_2);\mu)
 \operatorname{spec}(\mu)\,\dd\mu.
\end{split}
\label{eq:kuz-transforms}
\end{equation}
We omit the subscript $h$ when the test function is fixed.

\subsection{Spectral decomposition}
\label{subsec:kuz-normalization}
For $q\in\{1,p\}$, let
\[
  X_q=\Gamma_0(q)\backslash
  \SL_3(\mathbb{R})/\mathrm{SO}_3.
\]
We equip $X_q$ with its invariant probability measure $d\nu_q$ and use the inner product \[
  \langle f,g\rangle_q
  =\int_{X_q}f(z)\overline{g(z)}\,d\nu_q(z).
\]
The spectral side of the non-degenerate Kuznetsov formula contains the cuspidal, minimal Eisenstein, and maximal Eisenstein contributions described in \cite[Lemma~16]{BC}.

\subsubsection{Cuspidal spectrum}
\label{subsec:kuz-oldforms}
Let $\mathcal B(1)$ denote the set of level-one spherical cuspidal automorphic representations of $\GL_3$ with trivial central character.
For each $\pi\in\mathcal B(1)$, choose the Hecke-normalized spherical eigenform $f\in\pi$ with $A_\pi(1,1)=1$.
The space of oldforms generated by $f$ has dimension three, with basis $T_0f,T_1f,T_2f$ given in \cite[Proposition~5]{BC}.
Let $A_\pi^{(i)}(m,n)$ be their Fourier coefficients respectively.  We have
\begin{equation}
\begin{split}
A_\pi^{(0)}(m,n)= &~A_\pi(m,n),\\
A_\pi^{(1)}(m,n)= &~\sqrt p\left(A_\pi(m/p,pn)+A_\pi(m,n/p)\right),\\
A_\pi^{(2)}(m,n)= &~pA_\pi(m/p,n).
\end{split}
\label{eq:kuz-old-coeff}
\end{equation}
Coefficients with nonintegral indices are zero. Put
$\lambda=A_\pi(1,p)$, so $A_\pi(p,1)=\overline\lambda$ and set
\begin{equation*}
a_p\colon = \dfrac{\sqrt{p}(p+1)}{P},\qquad b_p \colon=\dfrac{p}{P}.
\end{equation*}
By \cite[Proposition 5]{BC}, we have
\begin{equation}
Q_\pi\colon =
\left(\dfrac{\langle T_if,T_jf\rangle_p}{\|f\| ^2}\right)_{0\leq i,j\leq2}
= \begin{pmatrix}
1&a_p\lambda&b_p\overline\lambda\\
a_p\overline\lambda&1+b_p|\lambda|^2&a_p\lambda\\
b_p\lambda&a_p\overline\lambda&1
\end{pmatrix}.
\label{eq:kuz-gram}
\end{equation}
Gram--Schmidt orthogonalization gives an orthonormal basis $S_0f,S_1f,S_2f$ as in \cite[Corollary~4]{BC}.
Let $B_\pi^{(j)}\|f\|^{-1}$ be the Fourier coefficients of $S_0f$, $S_1f$ and $S_2f$ respectively. Then 
\begin{equation}
\begin{split}
B_\pi^{(0)}(m,n)= &~A_\pi^{(0)}(m,n),\\
B_\pi^{(j)}(m,n)= &~(1+c_{jj})A_\pi^{(j)}(m,n)
+\sum_{i=0}^{j-1}c_{ji}A_\pi^{(i)}(m,n),
\qquad j=1,2.
\end{split}
\label{eq:kuz-orth-coeff}
\end{equation}
Here $c_{ji}\ll p^{\theta_3-1/2}$ and $Q_\pi=I+O(p^{\theta_3-1/2})$, with $\theta_3=5/14$.

Let $\pi\in\mathcal{B}(p)$, and let $F_\pi$ be the corresponding normalized newform.
Define
\[
    D_\pi(u)  \colon=
    \sum_{m_1,m_2\geq1}
    \dfrac{|A_\pi(m_1,m_2)|^2}{m_1^{2u}m_2^u}.
\]
At primes $\ell\neq p$,  \eqref{eq:local-Schur} gives
\[
    D_{\pi,\ell}(u)  =  \dfrac{L_\ell(u,\pi_\ell\times\widetilde\pi_\ell)} {\zeta_\ell(3u)}.
\]

At $\ell=p$, Lemmas~\ref{lem:local-prime-level} and \ref{ramified-coefficient} give
\begin{equation}
    D_{\pi,p}(u) =   L_p(u,\pi_p,\mathrm{Ad}).
    \label{eq:cusp-new-local-Dirichlet}
\end{equation}
It follows that
\begin{equation}
    D_\pi(u)
    =
    \frac{L(u,\pi\times\widetilde\pi)}{\zeta(3u)}
    \frac{1-p^{-u}}{1-p^{-3u}}.
    \label{eq:cusp-new-global-Dirichlet}
\end{equation}
Since
\[
    L(u,\pi\times\widetilde\pi)
    =
    \zeta(u)L(u,\pi,\mathrm{Ad}),
\]
taking the residue at \(u=1\) gives
\begin{align}
    \operatorname*{Res}_{u=1}D_\pi(u)
    &=
    \frac{1-p^{-1}}{1-p^{-3}}\,
    \frac{L(1,\pi,\mathrm{Ad})}{\zeta(3)}
    \notag\\
    &=
    \chi_p^{-1}
    \frac{L(1,\pi,\mathrm{Ad})}{\zeta(3)}.
    \label{eq:cusp-new-Dirichlet-residue}
\end{align}

By the Rankin--Selberg norm formula in
\cite[Section~4, proof of Theorem~4.1]{LM},
together with the essential-vector pairing formula in
\cite[Theorem~3.5]{Jo}, we obtain
\begin{equation}
    \|F_\pi\|^2
    =
    C_{3,\infty}(\mu_\pi)\,
    \chi_p^{-1}L(1,\pi,\mathrm{Ad}).
    \label{eq:cusp-new-probability-norm}
\end{equation}
Here \(L(1,\pi,\mathrm{Ad})\) is the full finite adjoint value, \(C_{3,\infty}(\mu_\pi)\) is the factor occurring in the
level-one norm formula and  $\chi_p\colon \!\!= \dfrac{P}{p^2}=1+\dfrac{1}{p}+\dfrac{1}{p^2}$.
After applying the archimedean normalization fixed above, taking the
reciprocal of \eqref{eq:cusp-new-probability-norm} gives the arithmetic
weight
\[
    \dfrac{c_{\mathrm{cusp}}\chi_p}
         {L(1,\pi,\mathrm{Ad})}
\]
in the cuspidal new spectrum.

For \(\pi\in\mathcal B(1)\), after orthogonalizing the three oldvectors, their common arithmetic weight is
\[
    \dfrac{c_{\mathrm{cusp}}}{L(1,\pi,\mathrm{Ad})}.
\]

\subsubsection{The minimal Eisenstein spectrum}
\label{subsec:kuz-minimal}

Now we consider the minimal Eisenstein series. Let $\mathbf{s}=(s_1,s_2,s_3)\in(\i\mathbb R)^3$, $s_3=-s_1-s_2$. The local Langlands parameters are 
\begin{equation}
\alpha_q=q^{s_1}, \quad \beta_q=q^{s_2}, \quad \gamma_q=q^{s_3}. 
\end{equation}
Define $A_{\mathbf s}(m,n)$ by \eqref{eq:local-Schur} at each prime $q$ and by multiplicativity. In particular,
\begin{equation*}
A_{\mathbf s}(1,n)= \sum_{abc=n}a^{s_1}b^{s_2}c^{s_3}.
\end{equation*}
For each fixed inducing parameter $\mathbf s$, the local space of level-$p$ invariant sections has dimension three \cite[Proposition~6]{BC}. Thus the minimal Eisenstein spectrum is parameterized by $j\in\{0,1,2\}$ and $(s_1,s_2)\in(\i\mathbb R)^2$, with $s_3=-s_1-s_2$. We define $A_{\mathbf s}^{(j)}$ by \eqref{eq:kuz-old-coeff} and $B_{\mathbf s}^{(j)}$ by \eqref{eq:kuz-orth-coeff}. The orthogonalization coefficients satisfy $c_{ji}\ll p^{-1/2}$ by the same proposition.

In \cite[Lemma 16]{BC}, the arithmetic weight is
\[
\dfrac{c_{\min}}
{\displaystyle
 \prod_{1\le i<k\le3}
 |\zeta(1+s_i-s_k)|^2}.
\]

\subsubsection{The maximal Eisenstein spectrum}
\label{subsec:kuz-maximal}

For $q\in\{1,p\}$, let $\mathcal B_2(q)$ denote the set of Hecke-normalized even Maass newforms of level $q$ and trivial central character. For $\phi\in\mathcal B_2(q)$, write $\lambda_\phi(n)$ for its Hecke eigenvalues, with $\lambda_\phi(1)=1$, and
$(\tau_\phi,-\tau_\phi)$ for its archimedean parameters.
By \cite[Theorem~1]{BB}, we have
\begin{equation*}
\tau_\phi=\i t_\phi\in\i\mathbb R\cup[-\theta_2,\theta_2],
\qquad \theta_2=7/64.
\end{equation*}
Its gamma factor is
$\Gamma_{\mathbb R}(w+\tau_\phi)\Gamma_{\mathbb R}(w-\tau_\phi)$. For $s\in\i\mathbb{R}$, the induced spectral parameters are $\mu=(s+\tau_\phi,s-\tau_\phi,-2s)$.
At each prime $\ell$, let $(\alpha_{\phi,\ell},\beta_{\phi,\ell})$ be the inverse roots of the local standard $L$-factor of $\phi_\ell$, allowing a zero inverse root at a ramified prime.
For $\phi\in\mathcal B_2(p)$, we have
\begin{equation*}
(\alpha_{\phi,p},\beta_{\phi,p})=(\lambda_\phi(p),0),
\qquad \lambda_\phi(p)=\pm p^{-1/2}.
\end{equation*}
As in \cite[Eq.\,(1.9)]{BC}, we define the inverse roots of the induced local standard $L$-factor by
\begin{equation*}
(\alpha_\ell,\beta_\ell,\gamma_\ell)
=(\ell^s\alpha_{\phi,\ell},\ell^s\beta_{\phi,\ell},\ell^{-2s}).
\end{equation*}
Then we define $A_{\phi,s}(m,n)$ by \eqref{eq:local-Schur} and by multiplicativity. 
In particular,
\begin{equation*}
A_{\phi,s}(1,n)= \sum_{ab=n}\lambda_\phi(a)a^s b^{-2s}.
\end{equation*}
If $\phi\in\mathcal B_2(1)$, we define the Fourier coefficients $A_{\phi,s}^{(j)}$ and $B_{\phi,s}^{(j)}$ obtained from \eqref{eq:kuz-old-coeff} and \eqref{eq:kuz-orth-coeff}. 
By \cite[Proposition 7]{BC}, the space of maximal generic unitary Eisenstein series is parameterized by 
\begin{enumerate}
\item $(\phi,s)$ with $\phi\in \mathcal{B}_2(p)$. The associated Eisenstein series have Fourier coefficients $A_{\phi,s}$. 
\item $(j,\phi, s)$ with $j\in \{0,1,2\}$ and $\phi\in \mathcal{B}_2(1)$. The associated Eisenstein series have Fourier coefficients $B_{\phi,s}^{(j)}$ with the constant $c_{ji} \ll p^{\theta_2-1/2}$.
\end{enumerate}

For $\phi\in\mathcal B_2(p)$, the norm here is taken with respect to probability measure on $\Gamma_0(p)\backslash\mathbb H$. Its volume for the hyperbolic measure $\dd x\,\dd y/y^2$ is $\pi(p+1)/3$. In the same archimedean Fourier normalization at levels one and $p$, the Rankin--Selberg norm identity \cite[Lemma~4.6, Eq.\,(4.8)]{HumphriesKhan} therefore gives
\[
\begin{aligned}
\|\phi\|^2
&=\frac{3}{\pi(p+1)}
\int_{\Gamma_0(p)\backslash\mathbb H}
 |\phi(x+\i y)|^2\frac{\dd x\,\dd y}{y^2}\\
&=C_{2,\infty}(\tau_\phi)\frac{p}{p+1}
 L(1,\phi,\operatorname{Ad}),
\end{aligned}
\]
where $C_{2,\infty}(\tau_\phi)$ is the same archimedean
factor as in the level-one formula.
Let $\Phi_{\phi,s}$ be the corresponding inducing section,
normalized so that its component on the nontrivial
double coset is $\phi$.
Hence \cite[Eq.\,(3.6) and Lemma 10]{BC}
gives
\[
\begin{aligned}
\|\Phi_{\phi,s}\|^2
&=
\frac{p(p+1)}{P}\,
\|\phi\|^2 \\
&=
C_{2,\infty}(\tau_\phi)
\frac{p(p+1)}{P}\frac{p}{p+1}
L(1,\phi,\operatorname{Ad}) \\
&=
C_{2,\infty}(\tau_\phi)\chi_p^{-1}
L(1,\phi,\operatorname{Ad}).
\end{aligned}
\]
Consequently, the maximal Eisenstein new spectrum carries the arithmetic weight
\[
    \frac{c_{\max}\chi_p}
         {L(1,\phi,\mathrm{Ad})\,|L(1+3s,\phi)|^2}.
\]

\subsection{Kloosterman sums and level conditions}
\label{subsec:kuz-kloosterman}

For $D_1|D_2$, as in \cite[\S~5.2]{BC}, we define 
\begin{equation}
\widetilde S(n_1,n_2,m_1,D_1,D_2)\colon \!\!= 
\sum_{\substack{C_1\bmod{D_1},\ C_2\bmod{D_2}\\
(C_1,D_1)=(C_2,D_2/D_1)=1}}
e\!\left(\dfrac{n_2\overline C_1C_2+n_1C_1}{D_1}
 +\dfrac{m_1\overline C_2}{D_2/D_1}\right).
\label{eq:kuz-short-def}
\end{equation}
Here $\overline C_1$ is an inverse modulo $D_1$, and $\overline C_2$ is an inverse modulo $D_2/D_1$. For $p|D_1$, $D_2$,
the sum for the long Weyl element is
\begin{equation}
\begin{split}
&~S^{(p)}(n_1,n_2,m_1,m_2,D_1,D_2)\\
\colon \!\!= &~\sum_{\substack{B_1,C_1\bmod{D_1},\ B_2,C_2\bmod{D_2}\\
D_1C_2+B_1B_2+D_2C_1\equiv0\bmod{D_1D_2}\\
(B_j,C_j,D_j)=1\ (j=1,2),\ p\mid B_1}}
e\!\left(\dfrac{n_1B_1+m_1(Y_1D_2-Z_1B_2)}{D_1}
 +\dfrac{n_2B_2+m_2(Y_2D_1-Z_2B_1)}{D_2}\right),
\end{split}
\label{eq:kuz-long-def}
\end{equation}
where $Y_jB_j+Z_jC_j\equiv1\bmod{D_j}$. The sum is independent of the permissible choices of $Y_j,Z_j$.

For the sum in \eqref{eq:kuz-short-def}, the required bound is
\cite[Lemma~14(d)]{BC}:
\begin{equation}
\widetilde S(n_1,n_2,m_1,D_1,D_2)\ll_\varepsilon 
\!\big((m_1,D_2/D_1)D_1^2,
(n_1,n_2,D_1)D_2\big)(D_1D_2)^\varepsilon.
\label{eq:short-Kloo-bound}
\end{equation}
We recall the Ramanujan sum
\begin{equation*}
r_p(a)=\sideset{}{^\star}\sum_{x\bmod{p}}e(ax/p)
\end{equation*}
 and the classical Kloosterman  sum
\begin{equation*}
S(a,b,c)=\sideset{}{^\star}\sum_{x\bmod{c}}e((ax+b\overline x)/c).
\end{equation*}

The first two identities below are \cite[Lemmas~14(b) and 15(a)]{BC}.
The third corrects the parameter order in \cite[Lemma~15(b)]{BC}.
\begin{equation*}
\begin{aligned}
S^{(p)}(n_1,n_2,m_1,m_2,p,p)= &~p-1+r_p(m_1)r_p(n_2),\\
S^{(p)}(n_1,n_2,m_1,m_2,p,p^2)= &~r_p(m_1)S(n_2,pm_2,p^2),\\
S^{(p)}(n_1,n_2,m_1,m_2,p^2,p)= &~
r_p(n_2)S(m_1,pn_1,p^2).
\end{aligned}
\end{equation*}

To verify the last identity, we set $D_1=p^2$, $D_2=p$ and write $B_1=pb$. The defining congruence becomes
\begin{equation*}
pC_2+bB_2+C_1\equiv 0\pmod{p^2}.
\end{equation*}
Primitivity forces $C_1$, $b$ and $B_2$ to be units modulo $p$. Choose $Y_1=0$, $Z_1=C_1^{-1}\pmod{p^2}$, and set
$x=-B_2C_1^{-1}\pmod{p^2}$. Then
\begin{equation*}
x^{-1}\equiv b+pC_2B_2^{-1}\pmod{p^2},\qquad
B_1\equiv px^{-1}\pmod{p^2}.
\end{equation*}
For each fixed unit $B_2\bmod p$, the map $(b,C_2)\mapsto x$ is a bijection onto $(\mathbb Z/p^2\mathbb Z)^\times$. The entire $m_2$ part of the phase is an integer because $p\mid B_1$. The remaining sums therefore separate as
\begin{equation*}
\begin{aligned}
S^{(p)}(n_1,n_2,m_1,m_2,p^2,p)
= &~\sideset{}{^\star}\sum_{B_2\bmod p}e(n_2B_2/p)
 \sideset{}{^\star}\sum_{x\bmod{p^2}}e\!\left(\dfrac{m_1x+pn_1x^{-1}}{p^2}\right)\\
= &~r_p(n_2)S(m_1,pn_1,p^2).
\end{aligned}
\end{equation*}

\subsection{The complete Kuznetsov formula at prime level}
\label{subsec:kuznetsov-identity}

We use the corrected kernel constants in the published version of Buttcane--Zhou \cite[Section 4]{BZ}. With the $\mu$-coordinates and the contour measure $\mathrm{d}\mu=\mathrm{d}\mu_1\,\mathrm{d}\mu_2$ fixed above, the level-one diagonal term in that normalization is
\[
\dfrac{\delta_{m_1,n_1}\delta_{m_2,n_2}}{256\pi^5}
\int_{\Re\mu=0}
h(\mu)\operatorname{spec}(\mu)\,d\mu.
\]
We multiply the level-one Kuznetsov formula by $256\pi^5$, so that the coefficient of the diagonal integral becomes one. The same scalar is incorporated into the spectral constants $c_{\mathrm{cusp}},c_{\min},c_{\max}$ and into the kernels $K_{w_4}$ and $K_{w_6}^{\mathrm{sym}}$. The kernel formulas in \S\,\ref{subsec:kuz-transforms} are written in this rescaled normalization, and the integral transforms $\Phi_{w_4,h}$, $\Phi_{w_5,h}$ and $\Phi_{w_6,h}$ are defined using these kernels without any additional rescaling.

Fix positive integers $m_1,m_2,n_1,n_2$. With the coefficients and normalizations above, define the five spectral terms by
\begin{equation*}
\begin{aligned}
\mathcal C_{\rm new}\colon \!\!= &~c_{\rm cusp}\sum_{\pi\in\mathcal B(p)}
\dfrac{h(\mu_\pi)\overline{A_\pi(m_1,m_2)}A_\pi(n_1,n_2)}{L(1,\pi,\mathrm{Ad})},\\
\mathcal C_{\rm old}\colon \!\!= &~c_{\rm cusp}\sum_{\pi\in\mathcal B(1)}\sum_{j=0}^2
\dfrac{h(\mu_\pi)\overline{B_\pi^{(j)}(m_1,m_2)}B_\pi^{(j)}(n_1,n_2)} {L(1,\pi,\mathrm{Ad})},\\
\mathcal E_{\min}:= &~c_{\min}\sum_{j=0}^2\int_{\Re(s_1)=0}\int_{\Re(s_2)=0}
\dfrac{h(\mathbf s)\overline{B_{\mathbf s}^{(j)}(m_1,m_2)}  B_{\mathbf s}^{(j)}(n_1,n_2)}{\prod_{1\leq i<k\leq3}|\zeta(1+s_i-s_k)|^2}\,\dd s_1\,\dd s_2,\\
\mathcal E_{\max,\rm new}\colon \!\!= &~c_{\max}\sum_{\phi\in\mathcal B_2(p)} \int_{\Re(s)=0}\dfrac{h(s+\tau_\phi,s-\tau_\phi,-2s) \overline{A_{\phi,s}(m_1,m_2)}A_{\phi,s}(n_1,n_2)}
 {L(1,\phi,\mathrm{Ad})|L(1+3s,\phi)|^2}\,\dd s,\\
\mathcal E_{\max,\rm old}\colon \!\!= &~c_{\max}\sum_{\phi\in\mathcal B_2(1)}\sum_{j=0}^2 \int_{\Re(s)=0}\dfrac{h(s+\tau_\phi,s-\tau_\phi,-2s)
 \overline{B_{\phi,s}^{(j)}(m_1,m_2)}B_{\phi,s}^{(j)}(n_1,n_2)} {L(1,\phi,\mathrm{Ad})|L(1+3s,\phi)|^2}\,\dd s.
\end{aligned}
\end{equation*}
Here $\dd s$ and $\dd s_1\dd s_2$ denote positive spectral measures, with the fixed normalizing factors included in $c_{\max}$ and $c_{\min}$ respectively.

\begin{theorem}
\label{thm:kuznetsov}
For $h\in\widetilde{\mathcal H}(10)$ and prime $p$,
the Kuznetsov formula is
\begin{equation}
\textcolor{black}{\chi_p}\mathcal C_{\rm new}
+\mathcal C_{\rm old}+\mathcal E_{\min}+\textcolor{black}{\chi_p}\mathcal E_{\max,\rm new}
+\mathcal E_{\max,\rm old}= \Delta+S_4+S_5+S_6,
\label{eq:Kuz-full}
\end{equation}
where
\begin{gather}
\Delta= P\,\delta_{m_1=n_1}\delta_{m_2=n_2}
 \int_{\Re(\mu)=0}h(\mu)\operatorname{spec}(\mu)\,\dd\mu,
\label{eq:Delta-exact}\\
S_4= P\sum_{\epsilon=\pm1}
\sum_{\substack{D_1,D_2\geq1\\pD_2\mid D_1\\m_2D_1=n_1D_2^2}}
\frac{\widetilde S(-\epsilon n_2,m_2,m_1,D_2,D_1)}{D_1D_2}
\Phi_{w_4}\!\left(\frac{\epsilon m_1m_2n_2}{D_1D_2}\right),\notag\\
S_5= P\sum_{\epsilon=\pm1}
\sum_{\substack{D_1,D_2\geq1\\p\mid D_1\mid D_2\\m_1D_2=n_2D_1^2}}
\frac{\widetilde S(\epsilon n_1,m_1,m_2,D_1,D_2)}{D_1D_2}
\Phi_{w_5}\!\left(\frac{\epsilon n_1m_1m_2}{D_1D_2}\right),\notag\\
S_6= P\sum_{\epsilon_1,\epsilon_2=\pm1}
\sum_{\substack{D_1,D_2\geq1\\p\mid D_1,\ p\mid D_2}}
\frac{S^{(p)}(\epsilon_2n_2,\epsilon_1n_1,m_1,m_2,D_1,D_2)}{D_1D_2}\Phi_{w_6}\!\left(
-\frac{\epsilon_2m_1n_2D_2}{D_1^2},
-\frac{\epsilon_1m_2n_1D_1}{D_2^2}\right).\notag
\end{gather}
\end{theorem}
\begin{proof}
Combining the discussion in \S\,\ref{subsec:kuz-oldforms}, \S\,\ref{subsec:kuz-maximal} and \cite[Lemma 16]{BC}, we can get the trace formula. 
\end{proof}

\begin{remark}
Since $\chi_p=1+O(p^{-1})$, conversion to full finite adjoint values changes the moment main terms within the error terms stated below.
\end{remark}

\begin{remark}
In the moment calculations below, we insert the appropriate approximate functional equation from \S\,\ref{subsec:afe} and combine the coefficient products by the Hecke relations before applying Theorem \ref{thm:kuznetsov}. If a dual sum is reindexed by $\pi\longmapsto\widetilde\pi$, the family and the full adjoint value are unchanged, and the full spectral test, including its weight function, is evaluated at $-\mu$. The diagonal contribution is determined by equality of the two resulting Fourier index pairs as in \eqref{eq:Delta-exact}.
\end{remark}

\subsection{Bounds for the Kuznetsov transforms}
\label{subsec:kuznetsov-kernel-bounds}

We recall the following kernel bounds from \cite[Lemmas~2.4, 2.5 and 2.6(a)]{BL}.

\begin{lemma}[Kuznetsov kernel bounds]
\label{lem:kuz-kernel-bounds}
Let $A,\varepsilon>0$ and $j,j_1,j_2\in\mathbb Z_{\geq0}$. For sufficiently large $a_0$ and $h\in\widetilde{\mathcal H}(a_0)$, the following bounds hold for $y,y_1,y_2\in\mathbb R\setminus\{0\}$. \par\noindent\textup{(1)} Short Weyl kernels:
\begin{equation}
|y|^j\left|\dfrac{\dd^j}{\dd y^j}\Phi_{w_i}(y)\right| \ll_{A,j,\varepsilon,h} |y|^{j/3} \min\{|y|^A,|y|^{-1/6+\varepsilon}\},\quad i=4,5.
\label{eq:kuz-short-bound}
\end{equation}
\par\noindent\textup{(2)} Derivatives of the long Weyl kernel:
\begin{equation*}
\begin{split}
&\left|y_1^{j_1}y_2^{j_2}
\dfrac{\partial^{j_1+j_2}\Phi_{w_6}(y_1,y_2)}
{\partial y_1^{j_1}\partial y_2^{j_2}}\right|\\
&\qquad\ll_{j_1,j_2,h}
\left(1+|y_1|^{1/2}+|y_1|^{1/3}|y_2|^{1/6}\right)^{j_1}
\left(1+|y_2|^{1/2}+|y_2|^{1/3}|y_1|^{1/6}\right)^{j_2}.\end{split}
\end{equation*}
\par\noindent\textup{(3)} Long Weyl kernel at small arguments:
\begin{equation}
\Phi_{w_6}(y_1,y_2)\ll_{A,\varepsilon,h} \min\{|y_1y_2^2|,|y_1^2y_2|\}^{A} (1+|y_1|+|y_2|)^\varepsilon.
\label{eq:w6-small}
\end{equation}
\par\noindent\textup{(4)} Global bound for the long Weyl kernel:
\begin{equation*}
\Phi_{w_6}(y_1,y_2)\ll_h  (2+|\log|y_1||+|\log|y_2||)^3 (1+|y_1|)^{-1/4}(1+|y_2|)^{-1/4}.
\end{equation*}
\end{lemma}

\begin{lemma}[Mean squares of Fourier coefficients]
\label{lem:coefficient-square-average}
Fix $0<\delta<1$ and $C\geq0$, and let $h$ be fixed and rapidly decreasing on the unitary spectrum. Uniformly for positive integers $a,b$ with $(ab,p)=1$ and $ab\leq p^{1-\delta}$,
\begin{equation}
\sum_{\pi\in\mathcal B(p)}
 \dfrac{|h(\mu_\pi)|(1+\|\mu_\pi\|)^C|A_\pi(a,b)|^2}
 {L(1,\pi,\mathrm{Ad})}
\ll_{h,C,\delta} p^2,
\label{eq:coefficient-square-average}
\end{equation}
and
\begin{equation}
\sum_{\pi\in\mathcal B(p)}
 \dfrac{|h(\mu_\pi)|(1+\|\mu_\pi\|)^C|A_\pi(a,b)|}
 {L(1,\pi,\mathrm{Ad})}
\ll_{h,C,\delta} p^2.
\label{eq:coefficient-absolute-average}
\end{equation}
\end{lemma}

\begin{proof}
Choose $A>\max(4/3,2/(3\delta))$ and take $K$ sufficiently large for the corresponding kernel estimates. For $R\geq1$, we set
\begin{equation*}
g_R(\mu)\colon \!\!= \exp\!\left(R^{-2}\sum_{j=1}^3\mu_j^2\right )  \prod_{1\leq i<j\leq3}\prod_{k=0}^K  \big((\mu_i-\mu_j)^2-(2k+1)^2\big)^2.
\end{equation*}
For $g_R$, the bounds \eqref{eq:kuz-short-bound} and \eqref{eq:w6-small} used below have implied constants $O_\delta(R^N)$, for some $N$ independent of $R$, $p$. On the unitary spectrum, $\mu$ and
$-\overline\mu$ lie in the same Weyl orbit, so $g_R\geq0$. The bound $|\Re(\mu_j)|\leq5/14$ also gives $g_R(\mu)\gg_K1$ when $\|\mu\|\leq R$.
Apply Theorem~\ref{thm:kuznetsov} with indices $(a,b,a,b)$ and the test function $g_R$.
Put $g=(a,b)$, $a^\prime=a/g$ and $b^\prime=b/g$. The short modulus pairs are
\begin{equation*}
(D_1,D_2)= (p^2a^\prime b^\prime k^2,pb^\prime k)
 \quad\text{or}\quad(pa^\prime k,p^2a^\prime b^\prime k^2),\qquad k\geq1.
\end{equation*}
Both kernel arguments have modulus $g^3/(p^3k^3)$. Since $g^2\leq ab<p$, the trivial Kloosterman bound and \eqref{eq:kuz-short-bound} give
\begin{equation*}
|S_4|+|S_5|\ll_\delta  p^2R^N(g/p)^{3A}\sum_{k\geq1}k^{-3A}\ll_\delta p^2R^N.
\end{equation*}
For the long term,
\begin{equation*}
\min\{|y_1^2y_2|,|y_1y_2^2|\}
= \dfrac{(ab)^3}{\max(D_1,D_2)^3},\qquad\text{and}\qquad
1+|y_1|+|y_2|\leq \dfrac{3\max(D_1,D_2)}p.
\end{equation*}
Apply \eqref{eq:w6-small} with decay exponent $A+1$ and growth exponent one. Using $ab\leq p\leq\max(D_1,D_2)$ gives
\begin{equation*}
\Phi_{w_6}(y_1,y_2)\ll_\delta
 R^N(ab)^{3A}\max(D_1,D_2)^{-3A}.
\end{equation*}
Thus
\begin{equation*}
\begin{aligned}
|S_6|\ll_\delta &~p^2R^N(ab)^{3A}
 \sum_{p\mid D_1,D_2}D_1D_2\max(D_1,D_2)^{-3A}\\
\ll_\delta &~R^Np^{4-3A}(ab)^{3A}
\leq R^Np^{4-3A\delta}\ll_\delta p^2R^N.
\end{aligned}
\end{equation*}
The diagonal term is $O_\delta(p^2R^N)$. Positivity therefore yields
\begin{equation*}
\sum_{\substack{\pi\in\mathcal B(p)\\\|\mu_\pi\|\leq R}}
 \dfrac{|A_\pi(a,b)|^2}{L(1,\pi,\mathrm{Ad})}
\ll_\delta p^2R^N.
\end{equation*}
Dyadic decomposition and rapid decrease of $h$ give \eqref{eq:coefficient-square-average}. The Cauchy--Schwarz inequality and the case $a=b=1$ of \eqref{eq:coefficient-square-average} give \eqref{eq:coefficient-absolute-average}.
\end{proof}

\subsection{Compatibility with the Kuznetsov test class}
\label{subsec:kuz-afe-compatibility}
We first work with a general test function. Fix a tube width $a\geq10$, an integer $J\geq\lceil a\rceil+4$ and $h_0\in\widetilde{\mathcal H}(a+1)$. We set
\begin{gather}
h(\mu)= P_J(0,\mu)h_0(\mu),
\label{eq:moment-test-class}\\
F_u^\pm(\mu)= h_0(\pm\mu)e^{u^2}P_J(u,\mu)
 \dfrac{\gamma_\mu(1/2+u)}{\gamma_{\pm\mu}(1/2)}.\notag
\end{gather}
Here $P_J$ is given by \eqref{eq:afe-polynomial}, and $F_u^-$ is the Mellin multiplier of $h(-\mu)W_{-\mu}$.
\begin{remark}
The fixed test function in \eqref{eq:positive-moment-test} is a special case of \eqref{eq:moment-test-class} obtained by taking $h_0$ as in \eqref{eq:moment-test-quotient}.
In this specialization, $h$ remains fixed and $h_0$ depends on $J$.
\end{remark}

\begin{lemma}[Uniform admissibility]
\label{lem:test-function-uniform}
For every fixed $C>0$ and $-2\leq\Re(u)\leq C$, the functions $F_u^\pm$ belong to $\widetilde{\mathcal H}(a)$. For all $M,k,R\geq0$,
\begin{equation}
\sup_{\substack{\mu\in\Lambda_a\\|\alpha|\leq k}}
(1+\|\mu\|)^M|\partial_\mu^\alpha F_u^\pm(\mu)|
\ll_{h_0,a,J,C,M,k,R} (1+|\Im(u)|)^{-R}.
\label{eq:hu-seminorm}
\end{equation}
\end{lemma}

\begin{proof}
The gamma recurrence gives
\begin{equation*}
P_J(u,\mu)\gamma_\mu(1/2+u) = (2\pi)^{3(J+1)} \dfrac{P_J(u,\mu)}
{\prod_{j=1}^3\prod_{\ell=0}^{J}(1/2+u-\mu_j+2\ell)}
 \gamma_\mu(1/2+u+2(J+1)).
\end{equation*}
By \eqref{eq:afe-polynomial}, the quotient is a polynomial. For $\mu\in\Lambda_{a+1}$ and $\Re(u)\geq-2$, the shifted gamma arguments have real parts at least $2J-a-1/2>0$. Thus $F_u^\pm$ are holomorphic on a neighborhood of $\Lambda_a$. And they are Weyl invariant and keep the structural zeros of $h_0$.

For $u=\sigma+\i v$ with $-2\leq\sigma\leq C$, \eqref{eq:gamma-ratio-bound} applied to the shifted product gives
\begin{equation*}
|F_u^\pm(\mu)|\ll_{a,J,C} |h_0(\pm\mu)|(1+\|\mu\|+|v|)^N e^{-v^2+3\pi|v|/4},
\end{equation*}
for some fixed $N$. Rapid decrease of $h_0$ on $\Lambda_{a+1}$, the Gaussian decay in $v$ and Cauchy's integral formula on the larger tube give \eqref{eq:hu-seminorm}.
\end{proof}

For $n\geq1$ and $c>0$, define the weighted test functions by \begin{equation}
H_n^\pm(\mu)\colon \!\!=h(\pm\mu)U_{\pm\mu}^\pm(n/\sqrt p)
= \dfrac1{2\pi \i}\int_{(c)}  \left(\dfrac{\sqrt p}{n}\right)^u F_u^\pm(\mu)\dfrac{\dd u}{u}.
\label{eq:weighted-tests-Mellin}
\end{equation}
Here $U_\mu^\pm$ are as in \S~\ref{subsec:afe}, and the integral gives the holomorphic continuation of the products. Lemma \ref{lem:test-function-uniform} allows any positive Mellin line and gives $H_n^\pm\in\widetilde{\mathcal H}(a)$ with
\begin{equation}
\sup_{\substack{\mu\in\Lambda_a\\|\alpha|\leq k}}
(1+\|\mu\|)^M|\partial_\mu^\alpha H_n^\pm(\mu)|
\ll_c (\sqrt p/n)^c,
\label{eq:weighted-test-bound}
\end{equation}
for every fixed $M,k\geq0$. 

\begin{remark}
\label{rem:kuz-kernel-uniformity}
For fixed exponents in \eqref{eq:kuz-short-bound} and \eqref{eq:w6-small}, their implied constants are bounded by a fixed finite sum of weighted suprema of the test function and its derivatives on a sufficiently wide fixed tube. The tube width and derivative orders depend only on the prescribed exponents. Indeed, the proofs of \cite[Lemmas~2.4 and 2.5]{BL} use finitely many fixed contour shifts and residues. Stirling's formula and Cauchy estimates on a slightly larger tube control the resulting integrals and residues by these suprema.

Consequently, \eqref{eq:weighted-test-bound} and \eqref{eq:product-test-seminorm} give the kernel bounds for $H_n^\pm$ and $H_t$ with respective additional factors $(\sqrt p/n)^c$ and $(p/t)^c$, for every fixed $c>0$. Similarly, \eqref{eq:holomorphic-square-test-seminorm} gives the bounds for $K_\mu^\pm(y)$ with a factor $y^{-c}$. Taking $c$ sufficiently large gives arbitrary decay in $p$ when $n>p^{1/2+\varepsilon}$ for fixed $\varepsilon>0$.
\end{remark}

\section{Harmonic twisted moments}
\label{sec:harmonic-twisted-moments}
In this section we evaluate the harmonic twisted
first and second moments, uniformly in the twist ranges needed
for mollification.
\subsection{The harmonic twisted first moment}
\label{sec:harmonic-first-moment}
For a general test $h$ as in \eqref{eq:moment-test-class},
define the harmonic twisted first moment and the corresponding
spectral integral by
\begin{equation}
M_h(r,s,p):= c_{\rm cusp}\sum_{\pi\in\mathcal B(p)}
 \frac{\overline{A_\pi(r,s)}L(1/2,\pi)}{L(1,\pi,\mathrm{Ad})}
 h(\mu_\pi),
\label{eq:def-Mrs}
\end{equation}
 and
 \begin{equation}
\mathcal H_h:= \int_{\Re(\mu)=0}h(\mu)
 \operatorname{spec}(\mu)\,\dd\mu.
\label{eq:Hh}
\end{equation}

\begin{theorem}
\label{thm:harmonic-first-moment}
Fix a positive integer $Q$ and $0<\eta<1/2$. Let $a_0\geq10$
be sufficiently large in terms of $\eta$, and fix an integer
$J\geq\lceil a_0\rceil+4$ and
$h_0\in\widetilde{\mathcal H}(a_0+1)$, independently of $p$.
Set $h(\mu)=P_J(0,\mu)h_0(\mu)$.
Then, for every $\varepsilon>0$, as $p\to\infty$ through primes
and uniformly for
\begin{equation*}
1\leq r\leq Qp^\eta,\qquad 1\leq s\leq Q,\qquad (rs,p)=1,
\end{equation*}
we have
\begin{equation}
M_h(r,s,p)= \frac{\delta_{r=1}}{\sqrt s}
 (p^2+p+1)\mathcal H_h
 +O_{Q,\eta,h_0,J,\varepsilon}(p^{5/4+\varepsilon}).
\label{eq:first-moment-main}
\end{equation}
Note that the fixed test \eqref{eq:positive-moment-test} satisfies these
hypotheses for every $0<\eta<1/2$, by
\textup{\S}\,\ref{subsec:kuz-afe-compatibility}.
\end{theorem}

Substituting \eqref{eq:AFE-final} into \eqref{eq:def-Mrs}, we have
\begin{equation}
\begin{split}
M_h(r,s,p)= &~c_{\rm cusp}\sum_{\pi\in\mathcal B(p)}
 \frac{h(\mu_\pi)\overline{A_\pi(r,s)}}{L(1,\pi,\mathrm{Ad})}
 \sum_{n\geq1}\frac{A_\pi(1,n)}{\sqrt n}
 V_{\mu_\pi}(n/\sqrt p)\\
&~ -\sqrt p\,c_{\rm cusp}\sum_{\pi\in\mathcal B(p)}
 \frac{h(\mu_\pi)\overline{A_\pi(r,s)}}{L(1,\pi,\mathrm{Ad})}
 \sum_{n\geq1}\frac{\overline{A_\pi(p,n)}}{\sqrt n}
 W_{\mu_\pi}(n/\sqrt p)\\
 = &~M_{h,1}(r,s,p)+M_{h,2}(r,s,p),\quad \text{say}.
\end{split}
\label{eq:first-moment-afe}
\end{equation}
The family $\mathcal B(p)$ and the full adjoint value are invariant
under contragredience. Reindexing $M_{h,2}(r,s,p)$ in
\eqref{eq:first-moment-afe} by $\pi\mapsto\widetilde\pi$, and using
\eqref{eq:contragredient-coeff},
\eqref{eq:unramified-unitarity},
\eqref{eq:functional-equation} and $(rs,p)=1$, we obtain
\begin{equation*}
M_{h,2}(r,s,p)= -\sqrt p\,c_{\rm cusp}
 \sum_{\pi\in\mathcal B(p)}
 \frac{h(-\mu_\pi)\overline{A_\pi(s,r)}}{L(1,\pi,\mathrm{Ad})}
 \sum_{n\geq1}\frac{A_\pi(p,n)}{\sqrt n}
 W_{-\mu_\pi}(n/\sqrt p).
\end{equation*}
With $H_n^\pm$ as in \eqref{eq:weighted-tests-Mellin}, for any fixed
$B>0$, truncation at
$n\leq p^{1/2+\varepsilon}$ yields
\begin{equation*}
\begin{aligned}
M_h(r,s,p)= &~c_{\rm cusp}
 \sum_{n\leq p^{1/2+\varepsilon}}\frac1{\sqrt n}
 \sum_{\pi\in\mathcal B(p)}
 \frac{\overline{A_\pi(r,s)}A_\pi(1,n)H_n^+(\mu_\pi)}
 {L(1,\pi,\mathrm{Ad})}\\
&~ -\sqrt p\,c_{\rm cusp}
 \sum_{n\leq p^{1/2+\varepsilon}}\frac1{\sqrt n}
 \sum_{\pi\in\mathcal B(p)}
 \frac{\overline{A_\pi(s,r)}A_\pi(p,n)H_n^-(\mu_\pi)}
 {L(1,\pi,\mathrm{Ad})}
 +O_{Q,\eta,h_0,J,\varepsilon,B}(p^{-B}).
\end{aligned}
\end{equation*}

Apply \eqref{eq:Kuz-full} with all four indices equal to $1$ and
the positive test $g_R$ from the proof of
Lemma~\ref{lem:coefficient-square-average}.
All five spectra satisfy $|\Re(\mu_j)|\leq5/14$, so $g_R\geq0$
and $g_R\gg1$ on $\|\mu\|\leq R$ in every sector.
The new coefficients at $(1,1)$ equal $1$, while the three oldvectors
associated with each level one representation satisfy
$\sum_{j=0}^2|B_\Pi^{(j)}(1,1)|^2
\geq|B_\Pi^{(0)}(1,1)|^2=1$.
Positivity and the same geometric estimates therefore bound the
arithmetic mass of each sector on $\|\mu\|\leq R$ by $O(p^2R^N)$,
counting each old representation once.
Dyadic decomposition and rapid decrease of $h$ give harmonic mass
$O_{h,C}(p^2)$ in each sector with weight
$(|h(\mu)|+|h(-\mu)|)(1+\|\mu\|)^C$, for every fixed $C\geq0$.
By \eqref{eq:kuz-old-coeff}--\eqref{eq:kuz-orth-coeff},
the unramified coefficient bounds and \eqref{eq:ramified-Hecke},
\begin{equation}
|B_\Pi^{(j)}(1,n)|+|B_\Pi^{(j)}(p,n)|
\ll_\eps p^{1+\eps}n^{\theta_3+\eps},
\qquad 0\leq j\leq2,
\label{eq:uniform-coefficient-tail-bound}
\end{equation}
in the three level one sectors.
The same bound holds for the new coefficients $A_\pi$ and $A_{\phi,z}$.
At $(r,s)$ and $(s,r)$, the twist coefficient is bounded by
$(rs)^{\theta_3+\eps}\ll_Q p^{\eta(\theta_3+\eps)}$ uniformly
in the stated range.
Taking the decay exponent in \eqref{eq:V-decay} sufficiently large
makes the contribution of $n>p^{1/2+\varepsilon}$, including the dual
factor $\sqrt p$, of size
$O_{Q,\eta,h_0,J,\varepsilon,B}(p^{-B})$ for any $B>0$.
These bounds also give absolute convergence and justify interchanging
the sums and spectral integrals.

Apply Theorem~\ref{thm:kuznetsov} with indices $(r,s,1,n)$ and
$(s,r,p,n)$ and tests $H_n^+$ and $H_n^-$, respectively.
Admissibility follows from Lemma~\ref{lem:test-function-uniform}.
The diagonal contribution is
\begin{equation}
\mathscr D:= P\sum_{n\leq p^{1/2+\varepsilon}}\frac1{\sqrt n}
 \int_{\Re(\mu)=0}\Bigl(
 \delta_{r=1}\delta_{s=n}H_n^+(\mu)
 -\sqrt p\,\delta_{s=p}\delta_{r=n}H_n^-(\mu)\Bigr)
 \operatorname{spec}(\mu)\,\dd\mu.
\label{eq:first-moment-diagonal-def}
\end{equation}
With the transforms in \eqref{eq:kuz-transforms}, define the short
Weyl contributions by
\begin{gather*}
\begin{aligned}
\mathscr S_4^+:= &~P
 \sum_{n\leq p^{1/2+\varepsilon}}\frac1{\sqrt n}
 \sum_{\epsilon=\pm1}
 \sum_{\substack{D_1,D_2\geq1\\pD_2\mid D_1\\sD_1=D_2^2}}
 \frac{\widetilde S(-\epsilon n,s,r,D_2,D_1)}{D_1D_2}
 \Phi_{w_4,H_n^+}\!\left(\frac{\epsilon rsn}{D_1D_2}\right),\\
\mathscr S_4^-:= &~P\sqrt p
 \sum_{n\leq p^{1/2+\varepsilon}}\frac1{\sqrt n}
 \sum_{\epsilon=\pm1}
 \sum_{\substack{D_1,D_2\geq1\\pD_2\mid D_1\\rD_1=pD_2^2}}
 \frac{\widetilde S(-\epsilon n,r,s,D_2,D_1)}{D_1D_2}
 \Phi_{w_4,H_n^-}\!\left(\frac{\epsilon rsn}{D_1D_2}\right),
\end{aligned}
\\
\begin{aligned}
\mathscr S_5^+:= &~P
 \sum_{n\leq p^{1/2+\varepsilon}}\frac1{\sqrt n}
 \sum_{\epsilon=\pm1}
 \sum_{\substack{D_1,D_2\geq1\\p\mid D_1\mid D_2\\rD_2=nD_1^2}}
 \frac{\widetilde S(\epsilon,r,s,D_1,D_2)}{D_1D_2}
 \Phi_{w_5,H_n^+}\!\left(\frac{\epsilon rs}{D_1D_2}\right),\\
\mathscr S_5^-:= &~P\sqrt p
 \sum_{n\leq p^{1/2+\varepsilon}}\frac1{\sqrt n}
 \sum_{\epsilon=\pm1}
 \sum_{\substack{D_1,D_2\geq1\\p\mid D_1\mid D_2\\sD_2=nD_1^2}}
 \frac{\widetilde S(\epsilon p,s,r,D_1,D_2)}{D_1D_2}
 \Phi_{w_5,H_n^-}\!\left(\frac{\epsilon prs}{D_1D_2}\right).
\end{aligned}
\end{gather*}
The long Weyl contributions are
\begin{equation}
\begin{split}
\mathscr S_6^+:= &~P
 \sum_{n\leq p^{1/2+\varepsilon}}\frac1{\sqrt n}
 \sum_{\epsilon_1,\epsilon_2=\pm1}
 \sum_{\substack{D_1,D_2\geq1\\p\mid D_1,\ p\mid D_2}}
 \frac{S^{(p)}(\epsilon_2n,\epsilon_1,r,s,D_1,D_2)}{D_1D_2}\Phi_{w_6,H_n^+}\!\left(
 -\frac{\epsilon_2rnD_2}{D_1^2},
 -\frac{\epsilon_1sD_1}{D_2^2}\right),\\*
\mathscr S_6^-:= &~P\sqrt p
 \sum_{n\leq p^{1/2+\varepsilon}}\frac1{\sqrt n}
 \sum_{\epsilon_1,\epsilon_2=\pm1}
 \sum_{\substack{D_1,D_2\geq1\\p\mid D_1,\ p\mid D_2}}
 \frac{S^{(p)}(\epsilon_2n,\epsilon_1p,s,r,D_1,D_2)}{D_1D_2}\\*
&~\qquad\times\Phi_{w_6,H_n^-}\!\left(
 -\frac{\epsilon_2snD_2}{D_1^2},
 -\frac{\epsilon_1prD_1}{D_2^2}\right).
\end{split}
\label{eq:first-moment-w6-def}
\end{equation}
Thus $\mathscr S_i^-$ includes the factor $\sqrt p$ but not the
minus sign in the moment expansion.

The four remaining spectral contributions are
\begin{gather*}
\begin{aligned}
\mathscr R_{\rm old}:= &~c_{\rm cusp}
 \sum_{n\leq p^{1/2+\varepsilon}}\frac1{\sqrt n}
 \sum_{\pi\in\mathcal B(1)}\sum_{j=0}^2
 \frac1{L(1,\pi,\mathrm{Ad})}\\
&~\quad\times\Bigl(
 \overline{B_\pi^{(j)}(r,s)}B_\pi^{(j)}(1,n)H_n^+(\mu_\pi)
 -\sqrt p\,\overline{B_\pi^{(j)}(s,r)}
 B_\pi^{(j)}(p,n)H_n^-(\mu_\pi)\Bigr),
\end{aligned}
\\
\begin{aligned}
\mathscr R_{\min}:= &~c_{\min}
 \sum_{n\leq p^{1/2+\varepsilon}}\frac1{\sqrt n}
 \sum_{j=0}^2\int_{\Re(s_1)=0}\int_{\Re(s_2)=0}
 \frac1{\prod_{1\leq i<k\leq3}|\zeta(1+s_i-s_k)|^2}\\
&~\quad\times\Bigl(
 \overline{B_{\mathbf s}^{(j)}(r,s)}
 B_{\mathbf s}^{(j)}(1,n)H_n^+(\mathbf s)
 -\sqrt p\,\overline{B_{\mathbf s}^{(j)}(s,r)}
 B_{\mathbf s}^{(j)}(p,n)H_n^-(\mathbf s)\Bigr)
 \,\dd s_1\,\dd s_2,
\end{aligned}
\\
\begin{aligned}
\mathscr R_{\max,\rm new}:= &~\chi_p c_{\max}
 \sum_{n\leq p^{1/2+\varepsilon}}\frac1{\sqrt n}
 \sum_{\phi\in\mathcal B_2(p)}\int_{\Re(z)=0}
 \frac1{L(1,\phi,\mathrm{Ad})|L(1+3z,\phi)|^2}\\
&~\quad\times\Bigl(
 \overline{A_{\phi,z}(r,s)}A_{\phi,z}(1,n)
 H_n^+(z+\tau_\phi,z-\tau_\phi,-2z)\\
&~\qquad\quad-\sqrt p\,\overline{A_{\phi,z}(s,r)}A_{\phi,z}(p,n)
 H_n^-(z+\tau_\phi,z-\tau_\phi,-2z)\Bigr)\,\dd z,
\end{aligned}
\\
\begin{aligned}
\mathscr R_{\max,\rm old}:= &~c_{\max}
 \sum_{n\leq p^{1/2+\varepsilon}}\frac1{\sqrt n}
 \sum_{\phi\in\mathcal B_2(1)}\sum_{j=0}^2\int_{\Re(z)=0}
 \frac1{L(1,\phi,\mathrm{Ad})|L(1+3z,\phi)|^2}\\
&~\quad\times\Bigl(
 \overline{B_{\phi,z}^{(j)}(r,s)}B_{\phi,z}^{(j)}(1,n)
 H_n^+(z+\tau_\phi,z-\tau_\phi,-2z)\\
&~\qquad\quad-\sqrt p\,\overline{B_{\phi,z}^{(j)}(s,r)}
 B_{\phi,z}^{(j)}(p,n)H_n^-(z+\tau_\phi,z-\tau_\phi,-2z)\Bigr)\,\dd z.
\end{aligned}
\end{gather*}
The continuous integrals use the measures specified in \S\,\ref{subsec:kuznetsov-identity}.
Set
\begin{equation*}
\begin{aligned}
\mathscr S:= &~\sum_{i=4}^6(\mathscr S_i^+-\mathscr S_i^-),\\
\mathscr R:= &~\mathscr R_{\rm old}+\mathscr R_{\min}
 +\mathscr R_{\max,\rm new}+\mathscr R_{\max,\rm old}.
\end{aligned}
\end{equation*}
Then the trace expansion is
\begin{equation}
\chi_p M_h(r,s,p)= \mathscr D+\mathscr S-\mathscr R
 +O_{Q,\eta,h_0,J,\varepsilon,B}(p^{-B}).
\label{eq:moment-decomposition}
\end{equation}

\subsubsection{The diagonal term}

The first diagonal in \eqref{eq:first-moment-diagonal-def}
requires $r=1$ and
$n=s$, which is included for large $p$ since
$s\leq Q$. The dual diagonal requires $s=p$ and $r=n$, contrary to
$(s,p)=1$. Hence
\begin{equation*}
\mathscr D= \frac{\delta_{r=1}P}{\sqrt s}
 \int_{\Re(\mu)=0}h(\mu)V_\mu(s/\sqrt p)
 \operatorname{spec}(\mu)\,\dd\mu.
\end{equation*}
For sufficiently large $p$, $s/\sqrt p\leq1$
uniformly for $s\leq Q$, so \eqref{eq:V-small} gives
$V_\mu(s/\sqrt p)=1+O_J(s^2p^{-1}(1+\|\mu\|)^C)$.
Using \eqref{eq:integrable-spectral-density}, we obtain
\begin{equation}
\mathscr D= \frac{\delta_{r=1}P}{\sqrt s}\mathcal H_h
+O_{Q,h,J}(p).
\label{eq:diagonal-asymptotic}
\end{equation}

\subsubsection{The short Weyl elements}

\begin{lemma}
\label{lem:moment-short}
The four short Weyl sums in $\mathscr S$ are
$O_{Q,\eta,h_0,J}(p^{-1})$, uniformly in the range of
Theorem~\ref{thm:harmonic-first-moment}.
\end{lemma}

\begin{proof}
Choose a fixed decay exponent $A$ sufficiently
large in terms of $\eta$, and take $a_0$ large enough for the
corresponding kernel estimates.
Counting the summands in
\eqref{eq:kuz-short-def} gives $|\widetilde S|\leq D_1D_2$.
Combining \eqref{eq:weighted-test-bound} with $c=2\eps$ and
\eqref{eq:kuz-short-bound}, we obtain, for either sign,
\begin{equation*}
|\Phi_{w_4,H_n^\pm}(y)|+|\Phi_{w_5,H_n^\pm}(y)|
\ll_{h_0,J,A,\eps} p^\eps|y|^A.
\end{equation*}
For $(m_1,m_2,n_1,n_2)=(r,s,1,n)$, the $w_4$ conditions
$sD_1=D_2^2$ and $pD_2\mid D_1$ give
\begin{equation*}
D_2= spk,\qquad D_1=sp^2k^2,\qquad
|y|=\frac{rn}{sp^3k^3},\qquad k\geq1.
\end{equation*}
For $w_5$, put $g=(r,n)$, $r'=r/g$, $n'=n/g$.
Writing $D_2=D_1d$ gives $r'd=n'D_1$. Since $(r',n')=1$
and $(r',p)=1$, all solutions are
\begin{equation*}
D_1= pr'k,\qquad D_2=p^2r'n'k^2,\qquad
|y|=\frac{rs}{p^3r'^2n'k^3}.
\end{equation*}

For $(s,r,p,n)$, the $w_4$ conditions $rD_1=pD_2^2$ and
$pD_2\mid D_1$ give
\begin{equation*}
D_2= rk,\qquad D_1=prk^2,\qquad
|y|=\frac{sn}{prk^3}.
\end{equation*}
For its $w_5$ term, put $g=(s,n)$, $s'=s/g$, $n'=n/g$ to obtain
\begin{equation*}
D_1= ps'k,\qquad D_2=p^2s'n'k^2,\qquad
|y|=\frac{sr}{p^2s'^2n'k^3}.
\end{equation*}
Summing $k^{-3A}$, using $s\leq Q$ and $r'^2n',s'^2n'\geq1$,
and then summing over $n$ gives
\begin{equation*}
\begin{aligned}
|\mathscr S_4^+|
\ll_{Q,h_0,J,A,\eps} &~p^{2+\eps-3A}r^A
 p^{(1/2+\eps)(A+1/2)},\\
|\mathscr S_5^+|
\ll_{Q,h_0,J,A,\eps} &~p^{2+\eps-3A}r^Ap^{1/4+\eps/2},\\
|\mathscr S_4^-|
\ll_{Q,h_0,J,A,\eps} &~p^{5/2+\eps-A}
 p^{(1/2+\eps)(A+1/2)},\\
|\mathscr S_5^-|
\ll_{Q,h_0,J,A,\eps} &~p^{5/2+\eps-2A}r^Ap^{1/4+\eps/2}.
\end{aligned}
\end{equation*}
Since $r\leq Qp^\eta$ and $\eta<1/2$, a sufficiently large fixed
$A$ and the arbitrarily small $\eps$ make all four terms
$O_{Q,\eta,h_0,J}(p^{-1})$.
\end{proof}

\subsubsection{The long Weyl element}

\begin{lemma}
\label{lem:moment-long}
The two long Weyl sums in $\mathscr S$ have total
size $O_{Q,\eta,h_0,J}(p^{-1})$, uniformly in the range of
Theorem~\ref{thm:harmonic-first-moment}.
\end{lemma}

\begin{proof}
The kernel arguments in
\eqref{eq:first-moment-w6-def}, with
$(m_1,m_2,n_1,n_2)=(r,s,1,n)$ or $(s,r,p,n)$, satisfy
\begin{equation*}
\begin{aligned}
|y_1^2y_2|= &~\frac{m_1^2m_2n_1n_2^2}{D_1^3},\\
|y_1y_2^2|= &~\frac{m_1m_2^2n_1^2n_2}{D_2^3}.
\end{aligned}
\end{equation*}
Hence
\begin{equation*}
\min(|y_1^2y_2|,|y_1y_2^2|)^A
\leq \frac{(m_1m_2n_1n_2)^{3A/2}}
 {(D_1D_2)^{3A/2}}.
\end{equation*}
For sufficiently large $p$, we have $r,s,n<p$ uniformly in the
stated range. At $(r,s,1,n)$ and $(s,r,p,n)$, the conditions
$D_1,D_2\geq p$ therefore give
\begin{equation*}
1+|y_1|+|y_2|\ll D_1D_2.
\end{equation*}
By \eqref{eq:weighted-test-bound} with $c=2\eps$,
\eqref{eq:w6-small} applies to $H_n^\pm$ with an additional
factor $p^\eps$. Counting the variables in
\eqref{eq:kuz-long-def} gives $|S^{(p)}|\leq(D_1D_2)^2$.
For fixed $n,\epsilon_1,\epsilon_2$, \eqref{eq:w6-small} gives
\begin{equation*}
\begin{aligned}
P\frac{|S^{(p)}(\epsilon_2n_2,\epsilon_1n_1,m_1,m_2,D_1,D_2)|}
 {D_1D_2}\,&~|\Phi_{w_6,H_n^\pm}(y_1,y_2)|\\
&~\ll_{h_0,J,A,\eps}
 p^{2+\varepsilon}(m_1m_2n_1n_2)^{3A/2}
 (D_1D_2)^{1-3A/2+\varepsilon}.
\end{aligned}
\end{equation*}
We first sum over $D_1,D_2$ with $p\mid D_1,D_2$.
For sufficiently large $A$, we have
\begin{equation*}
\sum_{p\mid D_1,D_2}(D_1D_2)^{1-3A/2+\varepsilon}
= p^{2-3A+2\varepsilon}\zeta(3A/2-1-\varepsilon)^2.
\end{equation*}
Since $m_1m_2n_1n_2$ equals $rsn$ in the first piece and $rspn$
in the dual piece, summing over $n\leq p^{1/2+\eps}$ with weight
$n^{-1/2}$ and over $\epsilon_1,\epsilon_2=\pm1$, and including
the additional factor $\sqrt p$ in the dual term, gives
\begin{equation*}
\begin{aligned}
\mathscr S_6^+\ll_{Q,h_0,J,A,\eps}
 &~p^{4+3\eps-3A}r^{3A/2}
 p^{(1/2+\eps)(3A/2+1/2)},\\
\mathscr S_6^-\ll_{Q,h_0,J,A,\eps}
 &~p^{9/2+3\eps-3A/2}r^{3A/2}
 p^{(1/2+\eps)(3A/2+1/2)}.
\end{aligned}
\end{equation*}
Since $r\leq Qp^\eta$ and $\eta<1/2$, taking $A$ sufficiently
large in terms of $\eta$ gives the stated bound.
\end{proof}

\subsubsection{The old and minimal Eisenstein spectra}
The following lemma records the cancellation in
the dual contribution obtained by summing over the three oldvectors.

\begin{lemma}
\label{lem:oldspace-cancellation}
Let $\Pi$ be a cuspidal, minimal Eisenstein, or maximal Eisenstein representation of level one. Suppose its unramified coefficients satisfy
$|A_\Pi(m,n)|\ll_\varepsilon(mn)^{\theta+\varepsilon}$,
with respectively $\theta=5/14$, $0$, or $7/64$.
For $(m_1m_2,p)=1$ and $n<p$,
\begin{equation*}
\begin{aligned}
\sum_{j=0}^2\overline{B_\Pi^{(j)}(m_1,m_2)}B_\Pi^{(j)}(1,n)
= &~\overline{A_\Pi(m_1,m_2)}A_\Pi(1,n)(Q_\Pi^{-1})_{00},\\
\sum_{j=0}^2\overline{B_\Pi^{(j)}(m_1,m_2)}B_\Pi^{(j)}(p,n)
\ll &~|A_\Pi(m_1,m_2)A_\Pi(1,n)|p^{\theta-1}.
\end{aligned}
\end{equation*}
Also $(Q_\Pi^{-1})_{00}\ll1$, uniformly as $p\to\infty$.
\end{lemma}

\begin{proof}

Put $\lambda=A_\Pi(1,p)$ in
\eqref{eq:kuz-gram}. Since $\Pi$ is unramified at $p$,
$|\lambda|\ll p^\theta$ and $\theta<1/2$ give
$Q_\Pi=I+O(p^{\theta-1/2})$, including
$b_p|\lambda|^2=O(p^{2\theta-1})$. Thus, uniformly in $\Pi$,
\begin{equation*}
Q_\Pi^{-1}= I+O(p^{\theta-1/2}),\qquad
(Q_\Pi^{-1})_{00}= O(1),\qquad
(Q_\Pi^{-1})_{01}=O(p^{\theta-1/2}).
\end{equation*}

Let $\mathfrak{C}$ be the change of basis from the original oldvectors to the
orthogonalized oldvectors in \eqref{eq:kuz-orth-coeff}.
Since the inner product is linear in its first variable,
$\mathfrak{C}Q_\Pi \mathfrak{C}^*=I$ and $\mathfrak{C}^*\mathfrak{C}=Q_\Pi^{-1}$. For coefficient columns
$a,b$, the sum over the three oldvectors is therefore
\begin{equation*}
\sum_{j=0}^2\overline{(\mathfrak{C}a)_j}(\mathfrak{C}b)_j= a^*Q_\Pi^{-1}b.
\end{equation*}
For $(m_1m_2,p)=1$ and $n<p$, the coefficient columns for the original
oldvectors at $(m_1,m_2)$ and
$(1,n)$ are $A_\Pi(m_1,m_2)e_0$ and $A_\Pi(1,n)e_0$, respectively,
where $e_0=(1,0,0)^{\mathsf t}$. This proves the first identity.
At $(p,n)$, multiplicativity, unitarity and $p\nmid n$ give
\begin{equation*}
\begin{aligned}
A_\Pi(p,n)= &~\overline\lambda A_\Pi(1,n),\\
A_\Pi^{(1)}(p,n)= &~\sqrt p\lambda A_\Pi(1,n),\\
A_\Pi^{(2)}(p,n)= &~pA_\Pi(1,n).
\end{aligned}
\end{equation*}
Thus the coefficient column is $A_\Pi(1,n)v$, with
$v=(\overline\lambda,\sqrt p\lambda,p)^{\mathsf t}$.

The definitions of $a_p,b_p$ give
\begin{equation*}
1-pb_p= p^{-1/2}a_p,\qquad
\sqrt p-pa_p=p^{-1/2}b_p.
\end{equation*}
Consequently
\begin{equation*}
v-pQ_\Pi e_2= p^{-1/2}(a_p\overline\lambda,b_p\lambda,0)^{\mathsf t},
\qquad e_2=(0,0,1)^{\mathsf t}.
\end{equation*}
Since $e_0^{\mathsf t}e_2=0$, this yields
\begin{equation*}
e_0^{\mathsf t}Q_\Pi^{-1}v
= p^{-1/2}\bigl(a_p\overline\lambda(Q_\Pi^{-1})_{00}
              +b_p\lambda(Q_\Pi^{-1})_{01}\bigr)
\ll p^{\theta-1}+p^{2\theta-2}
\ll p^{\theta-1}.
\end{equation*}
Taking absolute values before summing over the three oldvectors would lose this saving.
This proves the second assertion. The cancellation holds in the
weighted trace because $H_n^\pm$ is common to all three oldvectors.
\end{proof}

\begin{proposition}
\label{prop:old-remainder}
Uniformly in the range of Theorem~\ref{thm:harmonic-first-moment},
\begin{equation*}
|\mathscr R_{\rm old}|+|\mathscr R_{\min}|+|\mathscr R_{\max,\rm old}|
\ll_{h_0,J,\eps} (rs)^{5/14+\eps}p^{3/7+\eps}.
\end{equation*}
\end{proposition}

\begin{proof}

Lemma~\ref{lem:afe-weights} with $j=0$ gives, uniformly in the
three level one sectors and $n\leq p^{1/2+\eps}$,
\begin{equation*}
|H_n^+(\mu)|+|H_n^-(\mu)|
\ll_{J} \bigl(|h(\mu)|+|h(-\mu)|\bigr)(1+\|\mu\|)^C.
\end{equation*}
Apply Lemma~\ref{lem:oldspace-cancellation} to $(r,s)$ and $(s,r)$
in the first and dual pieces. The twist coefficient is
$O_\eps((rs)^{\theta+\eps})$, and the first piece, apart from
the rapidly decreasing spectral factor, is bounded by
\begin{equation*}
(rs)^{\theta+\eps}\sum_{n\leq p^{1/2+\eps}}n^{-1/2+\theta+\eps}
\ll_\eps (rs)^{\theta+\eps}
 p^{(1/2+\eps)(1/2+\theta+\eps)}.
\end{equation*}
By Lemma~\ref{lem:oldspace-cancellation}, the additional factor in the
dual contribution is $\sqrt p\,p^{\theta-1}=p^{\theta-1/2}\leq1$.

Brumley's lower bounds for Rankin--Selberg residues
and values \cite[Theorems~3 and 5]{Br}, together with an estimate for \(1/\zeta(s)\) \cite[Theorem~8.29]{IK}, give, at level one and for
$\Re(s_i)=\Re(z)=0$,
\begin{equation}
\begin{split}
L(1,\Pi,\mathrm{Ad})^{-1}\ll &~(1+\|\mu_\Pi\|)^C,
\qquad \Pi\in\mathcal B(1),\\
\prod_{i<j}|\zeta(1+s_i-s_j)|^{-2}
\ll &~(1+\|\mathbf s\|)^C,\\
\frac1{L(1,\phi,\mathrm{Ad})|L(1+3z,\phi)|^2}
\ll &~(1+|z|+|\tau_\phi|)^C,
\qquad \phi\in\mathcal B_2(1).
\label{eq:level-one-normalizer-bounds}
\end{split}
\end{equation}
Here
$L(1,\Pi,\mathrm{Ad})=\operatorname*{Res}_{w=1}
L(w,\Pi\times\widetilde\Pi)$, and similarly for $\phi$.
The standard $L$-value is covered by \cite[Theorem~5]{Br}
with degrees $2$ and $1$.
These factors admit continuous extensions, with value zero when the parameters coincide.
M\"uller's Weyl law \cite[Theorem~0.1]{Mu},
with $n=3$ and $n=2$, implies, for some fixed $C$ and all $R\geq1$,
\begin{equation}
\#\{\Pi\in\mathcal B(1):\|\mu_\Pi\|\leq R\}
+\#\{\phi\in\mathcal B_2(1):|\tau_\phi|\leq R\}
\ll R^C.
\label{eq:level-one-counting-bound}
\end{equation}
These bounds and rapid decrease of $h$ make the remaining sums
and integrals finite, independently of $p$.
Taking $\theta=5/14,0,7/64$, respectively,
we obtain
\begin{equation*}
\begin{aligned}
\mathscr R_{\rm old}\ll_{h_0,J,\eps}
 &~(rs)^{5/14+\eps}p^{3/7+\eps},\\
\mathscr R_{\min}\ll_{h_0,J,\eps}
 &~(rs)^\eps p^{1/4+\eps},
\end{aligned}
\end{equation*} and \begin{equation*}
\mathscr R_{\max,\rm old}\ll_{h_0,J,\eps}
 (rs)^{7/64+\eps}p^{39/128+\eps}.
\end{equation*}
The stated bound follows from these estimates.
\end{proof}

\subsubsection{The maximal Eisenstein new spectrum}
\label{subsec:first-moment-maximal-new}
We first bound the weighted Dirichlet series
arising from the maximal Eisenstein new spectrum.

\begin{lemma}
\label{lem:maximal-series}
Let $\phi\in\mathcal B_2(p)$, $z\in \i\mathbb R$ and
$\mu=(z+\tau_\phi,z-\tau_\phi,-2z)$. For either sign, and for some $C=C(J,\varepsilon)$, we have
\begin{equation*}
\sum_{n\geq1}\frac{A_{\phi,z}(1,n)}{\sqrt n}U_{\pm\mu}^\pm(n/\sqrt p)
\ll_{J,\varepsilon} p^{1/4+\varepsilon}
 (1+|z|+|\tau_\phi|)^C.
\end{equation*}
\end{lemma}

\begin{proof}

Here $\tau_\phi\in\i\mathbb R\cup[-7/64,7/64]$.
The convolution formula gives, in a right half plane,
\begin{equation*}
\sum_{n\geq1}\frac{A_{\phi,z}(1,n)}{n^w}
= \sum_{a,b\geq1}\frac{\lambda_\phi(a)a^zb^{-2z}}{(ab)^w}
= L(w-z,\phi)\zeta(w+2z).
\end{equation*}
On the line $\Re(u)=c>1$, absolute convergence justifies interchanging the integral and the $n$-sum. This gives
\begin{equation*}
\frac1{2\pi\i}\int_{(c)}p^{u/2}G_J(u,\mu)
 \frac{\gamma_\mu(1/2+u)}{\gamma_{\pm\mu}(1/2)}
 L(1/2+u-z,\phi)\zeta(1/2+u+2z)\frac{\dd u}{u}.
\end{equation*}
The denominator signs correspond to $V_\mu$ and $W_{-\mu}$,
respectively.

Shifting the integral to $\Re(u)=1/8$, we find that the pole of $\zeta$ at
$u=1/2-2z=1/2+\mu_3$ is cancelled by the factor
$(1/2+\mu_3)^2-u^2$ in $P_J$, defined by \eqref{eq:afe-polynomial}.
$P_J(0,\mu)$ is non-zero since $|\Re(\mu_j)|<1/2$, and
\begin{equation*}
\Re(1/2+u-\mu_j)\geq 1/2+1/8-7/64>0.
\end{equation*}
Thus no gamma pole or pole at $u=0$ is crossed.
By \eqref{eq:zeta-pole-removed} and \eqref{eq:afe-polynomial},
\begin{equation*}
|P_J(u,\mu)\zeta(1/2+u+2z)|
\ll_J (1+\|\mu\|+|u|)^{C_J},
\qquad \Re(u)\geq1/8.
\end{equation*}
The $\rm GL_2$ $L$-function is entire. Its bound
\cite[Theorem~5.41]{IK}, together with \eqref{eq:gamma-ratio-bound}
and the Gaussian $e^{u^2}$, makes the horizontal integrals tend
to $0$, so the shift crosses no residue.

For $u=1/8+\i v$, $\rm GL_2$ convexity \cite[Theorem~5.41]{IK},
\eqref{eq:zeta-five-eighths-bound} and \eqref{eq:gamma-ratio-bound}
give the following bounds, respectively, valid also for exceptional parameters,
\begin{equation*}
\begin{aligned}
L(5/8+\i v-z,\phi)
\ll_\varepsilon &~p^{3/16+\varepsilon}
 (1+|v|+|z|+|\tau_\phi|)^C,\\
|\zeta(5/8+\i v+2z)|\ll_\eps &~
(1+|v|+|z|)^{3/16+\eps},
\end{aligned}
\end{equation*}
and 
\begin{equation*}
\left|G_J(u,\mu)
 \frac{\gamma_\mu(1/2+u)}{\gamma_{\pm\mu}(1/2)}\right|
\ll_{J} (1+|z|+|\tau_\phi|+|v|)^{C_J}
 e^{-v^2+3\pi|v|/4}.
\end{equation*}
Integration over $v$ and
$p^{1/16}p^{3/16+\varepsilon}=p^{1/4+\varepsilon}$
prove the estimate. Unitarity gives
$\gamma_{-\mu}(1/2)=\overline{\gamma_\mu(1/2)}$, so the change of the
central denominator has modulus $1$.
\end{proof}
We next bound the harmonic averages of the $\mathrm{GL}_2$ Hecke coefficients.
\begin{lemma}
\label{lem:gl2-coefficient-average}
For $R\geq1$ and $(k,p)=1$, we have
\begin{equation}
\sum_{\substack{\phi\in\mathcal B_2(p)\\|\tau_\phi|\leq R}}
 \frac{|\lambda_\phi(k)|^2}{L(1,\phi,\mathrm{Ad})}
\ll pR^6\left(1+\frac{k^2}{p^2}\right).
\label{eq:gl2-coefficient-average}
\end{equation}
Let $h$ be fixed and rapidly decreasing on the unitary spectrum.
For every fixed $C\geq0$ and $\varepsilon>0$, uniformly for
$(abmn,p)=1$ and with $\mu=(z+\tau_\phi,z-\tau_\phi,-2z)$, we have
\begin{equation}
\begin{split}
&~\sum_{\phi\in\mathcal B_2(p)}
 \frac1{L(1,\phi,\mathrm{Ad})}
 \int_{\Re(z)=0}
 \frac{(|h(\mu)|+|h(-\mu)|)(1+\|\mu\|)^C}
 {|L(1+3z,\phi)|^2}
 |A_{\phi,z}(a,b)A_{\phi,z}(m,n)|\,\dd z\\
&~\qquad\qquad
\ll_{h,C,\varepsilon}
 p^{1+\varepsilon}(abmn)^\varepsilon
 \left(1+\frac{ab}{p}\right)
 \left(1+\frac{mn}{p}\right).
\end{split}
\label{eq:gl2-induced-coefficient-average}
\end{equation}
\end{lemma}

\begin{proof}
For this proof only, rescale $\phi$ to have norm $1$ for the hyperbolic measure $\dd x\,\dd y/y^2$ on $\Gamma_0(p)\backslash\mathbb H$. This leaves its Hecke eigenvalues unchanged. Write
$\rho_\phi(n)$ for its coefficients with Whittaker factor
$\sqrt yK_{\i t_\phi}(2\pi|n|y)$, where $\tau_\phi=\i t_\phi$.
Then $\rho_\phi(k)=\rho_\phi(1)\lambda_\phi(k)$ for $(k,p)=1$.
The Hecke relations, including $\lambda_\phi(p)=\pm p^{-1/2}$, give
\begin{equation*}
\sum_{n\geq1}\frac{|\lambda_\phi(n)|^2}{n^w}
= \frac{\zeta(w)L(w,\phi,\mathrm{Ad})}
{\zeta(2w)(1+p^{-w})},
\qquad \Re(w)>1.
\end{equation*}
The Eisenstein series at infinity has residue $3/(\pi(p+1))$
in \cite[\S\,15.4]{IK}, and the Bessel Mellin formula
\cite[(10.32.19)]{DLMF} yields
\begin{equation*}
\int_0^\infty K_{\i t_\phi}(y)^2\dd y
= \frac{\pi^2}{4\cosh(\pi t_\phi)}.
\end{equation*}
Since $|\rho_\phi(-n)|=|\rho_\phi(n)|$, Rankin--Selberg unfolding
and taking the residue at $w=1$ give
\begin{equation*}
\frac{3}{\pi(p+1)}
= 2|\rho_\phi(1)|^2(2\pi)^{-1}
 \frac{\pi^2}{4\cosh(\pi t_\phi)}
 \frac{pL(1,\phi,\mathrm{Ad})}{(p+1)\zeta(2)}.
\end{equation*}
Here the factor $2$ comes from the positive and negative Fourier
frequencies. Substituting $\zeta(2)=\pi^2/6$, we obtain
\begin{equation}
\frac{|\rho_\phi(1)|^2}{\cosh(\pi t_\phi)}
= \frac{2}{pL(1,\phi,\mathrm{Ad})},
\label{eq:gl2-fourier-normalization}
\end{equation}
and hence the conversion to the Kuznetsov weight is
\begin{equation*}
\frac{|\lambda_\phi(k)|^2}{L(1,\phi,\mathrm{Ad})}
= \frac p2\,\frac{|\rho_\phi(k)|^2}{\cosh(\pi t_\phi)}
\qquad((k,p)=1).
\end{equation*}
The Bessel integral converges for $|\Im(t_\phi)|<1/2$, so these
identities also hold for exceptional $t_\phi\in\i[-7/64,7/64]$.

For $R\geq 1$, we take the even entire test function
\begin{equation*}
b_R(t)\colon \!\!= e^{-t^2/R^2}(t^2+1/4)^2.
\end{equation*}
This is non-negative on the real and exceptional spectra and bounded below by a positive constant for $|\tau_\phi|\leq R$. Its transforms \cite[Eqs.\,(16.35) and (16.36)]{IK} are
\begin{equation*}
g_{0,R}= \pi^{-2}\int_{\mathbb R}t\,b_R(t)\tanh(\pi t)\dd t
 \ll R^6,
\end{equation*}
 \and\begin{equation*}
g_R^{(2)}(x)= \frac{2\i}{\pi}\int_{\mathbb R}
 J_{2\i t}(x)\frac{t\,b_R(t)}{\cosh(\pi t)}\dd t.
\end{equation*}
Moving the second contour to $t=v-\i$,
the only crossed pole is $-\i/2$, which can be cancelled by $b_R$.
Poisson's integral representation \cite[Eq.\,(10.9.4)]{DLMF} gives, for $x>0$,
\begin{equation*}
J_{2+2\i v}(x)= 
\frac{(x/2)^{2+2\i v}}{\sqrt\pi\,\Gamma(5/2+2\i v)}
 \int_{-1}^{1}e^{\i xt}(1-t^2)^{3/2+2\i v}\dd t,
\end{equation*}
 and 
 \begin{equation*}
|J_{2+2\i v}(x)|\ll 
x^2e^{\pi|v|}(1+|v|)^{-2}.
\end{equation*}
Stirling's formula and the bound
\begin{equation*}
|e^{-t^2/R^2}|\leq e^{1/R^2}e^{-(\Re(t))^2/R^2}
\qquad(-1\leq\Im(t)\leq0)
\end{equation*}
justify the contour shift.
Since $\cosh(\pi(v-\i))=-\cosh(\pi v)$,
\begin{equation*}
|g_R^{(2)}(x)|
\ll x^2\int_{\mathbb R}
 (1+|v|)^3e^{-v^2/R^2}\dd v
\ll x^2R^6\qquad(x>0).
\end{equation*}

Applying the $\rm GL_2$  Kuznetsov trace formula \cite[Theorem~16.3, Eq.\,(16.34)]{IK}
with both indices equal to $k$,
the old and continuous spectra are non-negative.
Using the trivial bound $|S(k,k;c)|\leq c$, the geometric contribution is bounded by
\begin{equation*}
\sum_{p\mid c}\frac{|S(k,k;c)|}{c}
 |g_R^{(2)}(4\pi k/c)|
\ll R^6k^2\sum_{p\mid c}c^{-2}
\ll R^6\frac{k^2}{p^2}.
\end{equation*}
The identity \eqref{eq:gl2-fourier-normalization} and the lower bound
for $b_R$ prove \eqref{eq:gl2-coefficient-average}.
In particular,
\begin{equation*}
\sum_{\substack{\phi\in\mathcal B_2(p)\\|\tau_\phi|\leq R}}
 \frac1{L(1,\phi,\mathrm{Ad})}
\ll pR^6.
\end{equation*}

The unramified Hecke eigenvalues of $\phi$ are real, since its
central character is trivial. The Hecke relations give
\begin{equation*}
\begin{aligned}
A_{\phi,z}(a,b)= &~
\sum_{d\mid(a,b)}\mu(d)
 \overline{A_{\phi,z}(1,a/d)}A_{\phi,z}(1,b/d),\\
A_{\phi,z}(1,n)= &~
\sum_{uv=n}\lambda_\phi(u)u^zv^{-2z},\\
\lambda_\phi(u)\lambda_\phi(v)= &~
\sum_{d\mid(u,v)}\lambda_\phi(uv/d^2).
\end{aligned}
\end{equation*}
For $z\in\i\mathbb R$, these express $A_{\phi,z}(a,b)$
as a sum of at most $\tau(ab)^4$ terms of the form
$c\,\lambda_\phi(k)$, with $|c|\leq1$ and $k\leq ab$.
The coefficients may depend on $z$, but the bound is uniform in $z$.
The Cauchy--Schwarz inequality,
the divisor bound $\tau(n)\ll_\varepsilon n^\varepsilon$
(see \cite[Eq.\,(1.81)]{IK}) and
\eqref{eq:gl2-coefficient-average} therefore give
\begin{equation*}
\sum_{\substack{\phi\in\mathcal B_2(p)\\|\tau_\phi|\leq R}}
 \frac{|A_{\phi,z}(a,b)|^2}{L(1,\phi,\mathrm{Ad})}
\ll_\varepsilon 
pR^6(ab)^\varepsilon\left(1+\frac{(ab)^2}{p^2}\right).
\end{equation*}
Applying the Cauchy--Schwarz inequality again yields the corresponding
product bound with factors
$(abmn)^\varepsilon(1+ab/p)(1+mn/p)$.
Applying \cite[Theorem~C(2)]{HR} to the unitary twist
$\phi\otimes|\det|^{3z}$, we obtain
\begin{equation}
|L(1+3z,\phi)|^{-2}
\ll_\varepsilon 
p^\varepsilon(1+|z|+|\tau_\phi|)^{C_\varepsilon}.
\label{eq:gl2-cross-reciprocal}
\end{equation}
Finally, $\|\mu\|^2=6|z|^2+2|\tau_\phi|^2$.
Rapid decrease of $h$, integration in $z$, and dyadic
decomposition in $|\tau_\phi|$ prove
\eqref{eq:gl2-induced-coefficient-average}.
Omitting the factor $|L(1+3z,\phi)|^{-2}$, the case $k=1$ of
\eqref{eq:gl2-coefficient-average} similarly gives,
for every fixed $C\geq0$,
\begin{equation}
\sum_{\phi\in\mathcal B_2(p)}\frac1{L(1,\phi,\mathrm{Ad})}
\int_{\Re(z)=0}
(|h(\mu)|+|h(-\mu)|)(1+\|\mu\|)^C\,\dd z
\ll_{h,C} p.
\label{eq:gl2-weighted-harmonic-mass}
\end{equation}
This completes the proof.
\end{proof}

\begin{proposition}
\label{prop:maximal-remainder}
Uniformly in the range of Theorem~\ref{thm:harmonic-first-moment},
\begin{equation*}
\mathscr R_{\max,\rm new}\ll_{Q,\eta,h_0,J,\varepsilon}
 p^{5/4+\varepsilon}.
\end{equation*}
\end{proposition}
\begin{proof}
By \S\,\ref{subsec:kuz-maximal}, the two non-zero standard parameters at $p$ are
$p^z\lambda_\phi(p)$ and $p^{-2z}$, with
$\lambda_\phi(p)=\pm p^{-1/2}$.
Their product is $p^{-z}\lambda_\phi(p)$.
The ramified Schur identity \eqref{eq:local-Schur} gives, for every $n\geq1$,
\begin{equation*}
A_{\phi,z}(p,n)= p^{-z}\lambda_\phi(p)A_{\phi,z}(1,n),\quad\text{and}\quad
\sqrt p\,|p^{-z}\lambda_\phi(p)|= 1.
\end{equation*}
Thus Lemma~\ref{lem:maximal-series} bounds both full sums, including
the dual factor $\sqrt p$, after extracting $h(\mu)$ or $h(-\mu)$.
The truncation argument following \eqref{eq:first-moment-afe},
using Lemma~\ref{lem:afe-weights} and
\eqref{eq:uniform-coefficient-tail-bound}, replaces the truncated
sums by full sums with error
$O_{Q,\eta,h_0,J,\eps,B}(p^{-B})$, uniformly in $r,s$.

Apply \eqref{eq:gl2-induced-coefficient-average} with one
coefficient at $(1,1)$ and the other at $(r,s)$ or $(s,r)$.
Together with Lemma~\ref{lem:maximal-series}, this gives
\begin{equation*}
|\mathscr R_{\max,\rm new}|
\ll_{h_0,J,\eps} \chi_p p^{5/4+\eps}(rs)^\eps
 \left(1+\frac{rs}{p}\right).
\end{equation*}
Since $rs\leq Q^2p^\eta<p$ for sufficiently large $p$ and
$\chi_p\ll1$, we obtain the stated bound.
\end{proof}

\subsubsection{Proof of Theorem~\ref{thm:harmonic-first-moment}}
Combining \eqref{eq:diagonal-asymptotic},
Lemmas~\ref{lem:moment-short} and \ref{lem:moment-long}, and
Propositions~\ref{prop:old-remainder} and \ref{prop:maximal-remainder}
in \eqref{eq:moment-decomposition}, and using $rs\leq Q^2p^\eta$,
we obtain
\begin{equation*}
\chi_p M_h(r,s,p)= \frac{\delta_{r=1}P}{\sqrt s}\mathcal H_h
 +O_{Q,\eta,h_0,J,\eps}(p^{5/4+\eps}).
\end{equation*}
Indeed, the level one terms are
$O_{Q,\eta,h_0,J,\eps}(p^{3/7+5\eta/14+\eps})$ and the
diagonal error is $O_{Q,h,J}(p)$.
Dividing by $\chi_p=1+O(p^{-1})$ changes the main term by
$O_{Q,h}(p)$ and proves \eqref{eq:first-moment-main}.

\subsection{The harmonic twisted second moment}
\label{sec:harmonic-second-moment}

From this subsection onward, we use the fixed test
\eqref{eq:positive-moment-test}.
For positive integers $a,b$, define the harmonic twisted second moment
\begin{equation*}
\mathcal N_2(a,b,p):= c_{\rm cusp}\sum_{\pi\in\mathcal B(p)}
 \frac{h(\mu_\pi)\overline{A_\pi(a,b)}|L(1/2,\pi)|^2}
 {L(1,\pi,\mathrm{Ad})}.
\end{equation*}
We write $\mathcal N_2(p):=\mathcal N_2(1,1,p)$ for the untwisted moment.
With $\mathcal H_h$ as in \eqref{eq:Hh}, we also set
\begin{equation*}
\kappa(\mu):= \frac{\gamma_\mu'}{\gamma_\mu}(1/2)
 +\frac{\gamma_{-\mu}'}{\gamma_{-\mu}}(1/2),
\end{equation*}
 and
 \begin{equation*}
\mathcal K_h:= \int_{\Re(\mu)=0}h(\mu)\kappa(\mu)
 \operatorname{spec}(\mu)\,\dd\mu.
\end{equation*}

The positivity and rapid decrease of $h$ give
$\mathcal H_h>0$ and absolute convergence of $\mathcal K_h$.

\begin{theorem} 
\label{thm:harmonic-second-moment}

Fix a positive integer $Q$ and $0<\eta<1/2$, independently of $p$.
Choose $a_0$ sufficiently large in terms of $\eta$
and $J\geq\lceil a_0\rceil+4$, independently of $p$.
Uniformly for positive integers $a,b$ with
$(ab,p)=1$ and $ab\leq Q^2p^{2\eta}$,
\begin{equation}
\begin{split}
\mathcal N_2(a,b,p)= &~\frac{P}{\sqrt{ab}}
 \left((1-p^{-1})
  \left\{\mathcal H_h\left(\log\frac p{ab}+2\gamma\right)
   +\mathcal K_h\right\}
  +\frac{2\log p}{p}\mathcal H_h\right)\\
&~{}+O_{h,Q,\eta,\varepsilon}\big(
 p^{7/4}(ab)^{-1/4}+p^{3/2+\varepsilon}\big).
\end{split}
\label{eq:twisted-second-moment}
\end{equation}
Here $\gamma$ is Euler's constant, and the small exponents in the
error can be chosen in terms of $\eta$. No condition $(a,b)=1$ is
imposed. In particular, uniformly in this range,
\begin{equation*}
\mathcal N_2(a,b,p)= \frac{P\mathcal H_h}{\sqrt{ab}}
 \log\frac p{ab}\left(1+O_{h,Q,\eta}\left(\frac1{\log p}\right)\right).
\end{equation*}
For $a=b=1$, taking $Q=1$ gives
\begin{equation}
\mathcal N_2(p)= P\mathcal H_h\log p
 +P(2\gamma\mathcal H_h+\mathcal K_h)
 +O_{h,\eta,\varepsilon}(p^{7/4+\varepsilon}).
\label{eq:ordinary-second}
\end{equation}
\end{theorem}

Substituting \eqref{eq:product-afe}
into the definition of
$\mathcal N_2(a,b,p)$ gives
\begin{equation*}
\mathcal N_2(a,b,p)= 2c_{\rm cusp}\sum_{\pi\in\mathcal B(p)}
 \frac{h(\mu_\pi)\overline{A_\pi(a,b)}}{L(1,\pi,\mathrm{Ad})}
 \sum_{m,n\geq1}
 \frac{\overline{A_\pi(1,m)}A_\pi(1,n)}{\sqrt{mn}}
 \mathcal U_{\mu_\pi}(mn/p).
\end{equation*}
By \eqref{eq:ramified-Hecke},
$|A_\pi(1,p^k)|\leq(1-p^{-1/2})^{-1}$ for $k\geq0$.

Writing $m=p^km'$ and $n=p^\ell n'$ with $(m'n',p)=1$, the terms
with $k+\ell\geq1$ satisfy, by \eqref{eq:product-weight-decay},
\begin{equation*}
\begin{aligned}
&~\sum_{\substack{m,n\geq1\\p\mid mn}}
 \frac{|A_\pi(1,m)A_\pi(1,n)|}{\sqrt{mn}}
 |\mathcal U_{\mu_\pi}(mn/p)|\\
&~\ll_{B,J,\varepsilon} (1+\|\mu_\pi\|)^C
 \sum_{k+\ell\geq1}p^{-(k+\ell)/2}
 \sum_{m',n'\geq1}
 (m'n')^{-1/2+\theta_3+\varepsilon}
 (1+p^{k+\ell-1}m'n')^{-B}\\
&~\ll_{J,\varepsilon} p^{-1/2}(1+\|\mu_\pi\|)^C.
\end{aligned}
\end{equation*}
Here $B>2+2\theta_3+2\varepsilon$ is fixed.
For sufficiently large $p$, we have
$ab\leq Q^2p^{2\eta}\leq p^{1/2+\eta}$.
Apply \eqref{eq:coefficient-absolute-average} with
$\delta=1/2-\eta$. The total ramified contribution is therefore
$O_{h,Q,\eta,\varepsilon}(p^{3/2+\varepsilon})$, uniformly in $a,b$.
For $(mn,p)=1$, the Hecke relation gives
\begin{equation*}
\overline{A_\pi(1,m)}A_\pi(1,n)
= \sum_{r\mid(m,n)}A_\pi(m/r,n/r).
\end{equation*}
Replacing the indices by $rm,rn$, we obtain
\begin{equation*}
\mathcal N_2(a,b,p)= 2c_{\rm cusp}
 \sum_{\substack{r,m,n\geq1\\(rmn,p)=1}}\frac1{r\sqrt{mn}}
 \sum_{\pi\in\mathcal B(p)}
 \frac{\overline{A_\pi(a,b)}A_\pi(m,n)H_{r^2mn}(\mu_\pi)}
 {L(1,\pi,\mathrm{Ad})}
 +O_{h,Q,\eta,\varepsilon}(p^{3/2+\varepsilon}).
\end{equation*}

The weighted mass bound after \eqref{eq:first-moment-afe} and
\eqref{eq:product-weight-decay} justify the interchanges and allow
truncation at $r^2mn\leq p^{1+\varepsilon}$ with error
$O_{h,Q,\eta,\varepsilon,B}(p^{-B})$,
uniformly for $ab\leq Q^2p^{2\eta}$.
Apply Theorem~\ref{thm:kuznetsov} with indices $(a,b,m,n)$ and test
$H_{r^2mn}$. Admissibility follows from Lemma~\ref{lem:product-weights}.
Removing the truncation $r^2mn\leq p^{1+\varepsilon}$ from the diagonal
and the other spectral terms introduces an error of the same size. Define
\begin{equation}
\mathcal D_2(a,b,p):= \frac{2P}{\sqrt{ab}}
 \sum_{\substack{r\geq1\\(r,p)=1}}\frac1r
 \int_{\Re(\mu)=0}H_{r^2ab}(\mu)
 \operatorname{spec}(\mu)\,\dd\mu.
\label{eq:second-diagonal-def}
\end{equation}
The three Weyl contributions are
\begin{gather}
\mathcal S_4(a,b):= 2P
 \sum_{\substack{r,m,n\geq1\\(rmn,p)=1\\r^2mn\leq p^{1+\varepsilon}}}
 \frac1{r\sqrt{mn}}\sum_{\epsilon=\pm1}
 \sum_{\substack{D_1,D_2\geq1\\pD_2\mid D_1\\bD_1=mD_2^2}}
 \frac{\widetilde S(-\epsilon n,b,a,D_2,D_1)}{D_1D_2}\Phi_{w_4,H_{r^2mn}}\!
 \left(\frac{\epsilon abn}{D_1D_2}\right),
\label{eq:second-w4-def}\\
\mathcal S_5(a,b):= 2P
 \sum_{\substack{r,m,n\geq1\\(rmn,p)=1\\r^2mn\leq p^{1+\varepsilon}}}
 \frac1{r\sqrt{mn}}\sum_{\epsilon=\pm1}
 \sum_{\substack{D_1,D_2\geq1\\p\mid D_1\mid D_2\\aD_2=nD_1^2}}
 \frac{\widetilde S(\epsilon m,a,b,D_1,D_2)}{D_1D_2}\Phi_{w_5,H_{r^2mn}}\!
 \left(\frac{\epsilon mab}{D_1D_2}\right),
\label{eq:second-w5-def}\\
\begin{aligned}
\mathcal S_6(a,b):= &~2P
 \sum_{\substack{r,m,n\geq1\\(rmn,p)=1\\r^2mn\leq p^{1+\varepsilon}}}
 \frac1{r\sqrt{mn}}\sum_{\epsilon_1,\epsilon_2=\pm1}
 \sum_{\substack{D_1,D_2\geq1\\p\mid D_1,\ p\mid D_2}}
 \frac{S^{(p)}(\epsilon_2n,\epsilon_1m,a,b,D_1,D_2)}{D_1D_2}
 \\
&~\qquad\times\Phi_{w_6,H_{r^2mn}}\!
 \left(-\frac{\epsilon_2anD_2}{D_1^2},
       -\frac{\epsilon_1bmD_1}{D_2^2}\right).
\end{aligned}
\label{eq:second-w6-def}
\end{gather}
The four remaining spectral contributions, with full index sums, are
\begin{gather}
\mathcal R_{\rm old}(a,b):= 2c_{\rm cusp}
 \sum_{\substack{r,m,n\geq1\\(rmn,p)=1}}\frac1{r\sqrt{mn}}
 \sum_{\pi\in\mathcal B(1)}\sum_{j=0}^2
 \frac{\overline{B_\pi^{(j)}(a,b)}B_\pi^{(j)}(m,n)
 H_{r^2mn}(\mu_\pi)}{L(1,\pi,\mathrm{Ad})},\notag
\\
\begin{aligned}
\mathcal R_{\min}(a,b):= &~2c_{\min}
 \sum_{\substack{r,m,n\geq1\\(rmn,p)=1}}\frac1{r\sqrt{mn}}
 \sum_{j=0}^2\int_{\Re(s_1)=0}\int_{\Re(s_2)=0}\\
&~\qquad\times
 \frac{\overline{B_{\mathbf s}^{(j)}(a,b)}
 B_{\mathbf s}^{(j)}(m,n)H_{r^2mn}(\mathbf s)}
 {\prod_{1\leq i<k\leq3}|\zeta(1+s_i-s_k)|^2}
 \,\dd s_1\,\dd s_2,
\end{aligned}\notag
\\
\begin{aligned}
\mathcal R_{\max,\rm new}(a,b):= &~2\chi_p c_{\max}
 \sum_{\substack{r,m,n\geq1\\(rmn,p)=1}}\frac1{r\sqrt{mn}}
 \sum_{\phi\in\mathcal B_2(p)}\int_{\Re(z)=0}\\
&~\qquad\times
 \frac{\overline{A_{\phi,z}(a,b)}A_{\phi,z}(m,n)
 H_{r^2mn}(z+\tau_\phi,z-\tau_\phi,-2z)}
 {L(1,\phi,\mathrm{Ad})|L(1+3z,\phi)|^2}\,\dd z,
\end{aligned}
\label{eq:second-max-new-def}\\
\begin{aligned}
\mathcal R_{\max,\rm old}(a,b):= &~2c_{\max}
 \sum_{\substack{r,m,n\geq1\\(rmn,p)=1}}\frac1{r\sqrt{mn}}
 \sum_{\phi\in\mathcal B_2(1)}\sum_{j=0}^2
 \int_{\Re(z)=0}\\
&~\qquad\times
 \frac{\overline{B_{\phi,z}^{(j)}(a,b)}B_{\phi,z}^{(j)}(m,n)
 H_{r^2mn}(z+\tau_\phi,z-\tau_\phi,-2z)}
 {L(1,\phi,\mathrm{Ad})|L(1+3z,\phi)|^2}\,\dd z.
\end{aligned}\notag
\end{gather}
The transforms and continuous measures are normalized as in \S\,\ref{subsec:kuznetsov-identity}. The trace expansion is
\begin{equation}
\begin{split}
\chi_p\mathcal N_2(a,b,p)= &~\mathcal D_2(a,b,p)
 +\sum_{i=4}^6\mathcal S_i(a,b)
 -\mathcal R_{\rm old}(a,b)-\mathcal R_{\min}(a,b)\\
&~-\mathcal R_{\max,\rm new}(a,b)
 -\mathcal R_{\max,\rm old}(a,b)
 +O_{h,Q,\eta,\varepsilon}(p^{3/2+\varepsilon}).
\end{split}
\label{eq:second-moment-decomposition}
\end{equation}

\subsubsection{The diagonal term}
\label{subsec:ordinary-diagonal}

The diagonal conditions are $m=a,n=b$.
With $R_\mu$ as defined in the proof of
Lemma~\ref{lem:product-weights},
\eqref{eq:second-diagonal-def} and \eqref{eq:product-weight-mellin}
give, by the Mellin inversion formula on $\Re(u)=c>1$,
\begin{equation*}
2\sum_{\substack{r\geq1\\(r,p)=1}}\frac1r
 \mathcal U_\mu(r^2ab/p)
= \frac{2}{2\pi\i}\int_{(c)}
 \left(\frac{ab}{p}\right)^{-u}G_J(u,\mu)R_\mu(u)
 \zeta^{(p)}(1+2u)\frac{\dd u}{u},
\end{equation*}
where $\zeta^{(p)}(w)=(1-p^{-w})\zeta(w)$.
For $\Re(\mu)=0$, move the contour to $\Re(u)=-1/4$, crossing
only the double pole at $u=0$. To compute its residue, first use
the Laurent expansion of the zeta factor,
\begin{equation*}
\zeta^{(p)}(1+2u)=\frac{1-p^{-1}}{2u}
 +(1-p^{-1})\gamma+\frac{\log p}{p}+O_p(u).
\end{equation*}
The remaining factor is holomorphic and equals one at zero.
Since $G_J$ is even and $R_\mu'(0)=\kappa(\mu)$, its Taylor expansion is
\begin{equation*}
\left(\frac{p}{ab}\right)^uG_J(u,\mu)R_\mu(u)
=1+u\left(\log\frac{p}{ab}+\kappa(\mu)\right)
 +O_{p,a,b,\mu,J}(u^2).
\end{equation*}
Multiplying these expansions by $2/u$ and applying the residue theorem gives
\begin{equation*}
2\sum_{\substack{r\geq1\\(r,p)=1}}\frac1r
 \mathcal U_\mu(r^2ab/p)
= (1-p^{-1})\left(\log\frac{p}{ab}+2\gamma+\kappa(\mu)\right)
 +\frac{2\log p}{p}
+O_J\big((ab/p)^{1/4}(1+\|\mu\|)^C\big).
\end{equation*}
The error comes from the integral on $\Re(u)=-1/4$, where
$|1-p^{-1-2u}|\leq1+p^{-1/2}$.
By \eqref{eq:zeta-pole-removed} and \eqref{eq:gamma-ratio-bound},
for $\Re(\mu)=0$,
\begin{equation*}
|G_J(-1/4+\i v,\mu)R_\mu(-1/4+\i v)
\zeta^{(p)}(1/2+2\i v)|
\ll_J (1+\|\mu\|+|v|)^C e^{-v^2+3\pi|v|/2}.
\end{equation*}
This justifies the integral estimate and the horizontal limits.
Integrating against $h$ using \eqref{eq:integrable-spectral-density} gives
\begin{equation}
\mathcal D_2(a,b,p)= \frac{P}{\sqrt{ab}}
 \left((1-p^{-1})
 \left\{\mathcal H_h\left(\log\frac p{ab}+2\gamma\right)
 +\mathcal K_h\right\}
 +\frac{2\log p}{p}\mathcal H_h\right)+O_{h,\eta}\big(p^{7/4}(ab)^{-1/4}\big).
\label{eq:ordinary-diagonal}
\end{equation}

\subsubsection{The short Weyl elements}
\label{sec:geometry}

\begin{lemma} 
\label{lem:second-moment-short}
Uniformly in the twist range of Theorem~\ref{thm:harmonic-second-moment},
$\mathcal S_4(a,b)$ and $\mathcal S_5(a,b)$ are $O_{h,Q,\eta}(p^{-1})$.
\end{lemma}

\begin{proof}
Choose $A$ sufficiently large in terms of $\eta$ and take the tube
width $a_0$ large enough for the kernel estimates. By
\eqref{eq:product-test-seminorm} and \eqref{eq:kuz-short-bound},
\begin{equation*}
|\Phi_{w_4,H_t}(y)|+|\Phi_{w_5,H_t}(y)|
\ll_{h,\eta,\eps} p^\eps|y|^A,\qquad t\geq1.
\end{equation*}
The modulus conditions in \eqref{eq:second-w4-def} and
\eqref{eq:second-w5-def} give, respectively,
\begin{equation*}
\begin{aligned}
D_2= &~\frac{pbk}{(b,m)},\qquad
D_1=\frac{p^2mbk^2}{(b,m)^2},\\
D_1= &~\frac{pak}{(a,n)},\qquad
D_2=\frac{p^2nak^2}{(a,n)^2},\qquad k\geq1.
\end{aligned}
\end{equation*}
In both cases $D_1D_2\geq p^3k^3$. Hence, on
$r^2mn\leq p^{1+\eps}$, the kernel argument satisfies
\begin{equation*}
|y|\leq \frac{ab\,p^{-2+\eps}}{k^3}.
\end{equation*}
A direct computation gives
\begin{equation}
\sum_{r^2mn\leq p^{1+\eps}}\frac1{r\sqrt{mn}}
\ll_\eps p^{1/2+\eps}.
\label{eq:second-outside-weights}
\end{equation}
Using $|\widetilde S|\leq D_1D_2$, $P\ll p^2$ and
$\sum_{k\geq1}k^{-3A}\ll_A1$, we obtain
\begin{equation*}
|\mathcal S_4(a,b)|+|\mathcal S_5(a,b)|
\ll_{h,\eta,\eps} 
 p^{5/2+(A+2)\eps}(ab)^A p^{-2A}.
\end{equation*}
Since $ab\leq Q^2p^{2\eta}$, $Q$ is fixed and $\eps$ is arbitrarily
small, the claim follows by our choice of $A$.
\end{proof}

\subsubsection{The long Weyl element}
\label{subsec:second-long}

\begin{lemma} 
\label{lem:second-moment-long}
Uniformly in the twist range of Theorem~\ref{thm:harmonic-second-moment},
$\mathcal S_6(a,b)=O_{h,Q,\eta}(p^{-1})$.
\end{lemma}

\begin{proof}
Use the same choice of $A$ as in the preceding proof.
\eqref{eq:w6-small} and
\eqref{eq:product-test-seminorm} give
\begin{equation*}
|\Phi_{w_6,H_t}(y_1,y_2)|\ll_{h,\eta,\eps} p^\eps
 \min(|y_1^2y_2|,|y_1y_2^2|)^A
 (1+|y_1|+|y_2|)^\eps.
\end{equation*}
The arguments in \eqref{eq:second-w6-def} satisfy
\begin{equation*}
|y_1^2y_2|= \frac{a^2bmn^2}{D_1^3},\qquad
|y_1y_2^2|=\frac{ab^2m^2n}{D_2^3}.
\end{equation*}
Further, we have
\begin{equation*}
\min(|y_1^2y_2|,|y_1y_2^2|)^A
\leq \frac{(abmn)^{3A/2}}{(D_1D_2)^{3A/2}}.
\end{equation*}
For sufficiently large $p$,
we have $an,bm\leq abmn\leq Q^2p^{2\eta+1+\varepsilon}<p^3$.
Since $D_1,D_2\geq p$, this gives
$1+|y_1|+|y_2|\leq3D_1D_2$.
Using $|S^{(p)}|\leq(D_1D_2)^2$ and
\begin{equation*}
\sum_{p\mid D_1,D_2}(D_1D_2)^{1-3A/2+\eps}
= p^{2-3A+2\eps}\zeta(3A/2-1-\eps)^2,
\end{equation*}
and then \eqref{eq:second-outside-weights}, gives
\begin{equation*}
|\mathcal S_6(a,b)|\ll_{h,\eta,\eps} 
 p^{9/2+(3A/2+4)\eps}
 (ab)^{3A/2}p^{-3A/2}.
\end{equation*}

Since $ab\leq Q^2p^{2\eta}$, $Q$ is fixed and $\eta<1/2$,
our choice of $A$ and the arbitrarily small $\eps$ give the claimed bound.
The sum over $D_1,D_2$ converges absolutely, completing the proof.
\end{proof}

\subsubsection{The old and minimal Eisenstein spectra}
\label{sec:spectrum}

By \eqref{eq:kuz-old-coeff}--\eqref{eq:kuz-orth-coeff}, for $(abmn,p)=1$,
\begin{equation*}
\sum_{j=0}^2\overline{B_\Pi^{(j)}(a,b)}B_\Pi^{(j)}(m,n)
= \overline{A_\Pi(a,b)}A_\Pi(m,n)(Q_\Pi^{-1})_{00},
\end{equation*}
where $(Q_\Pi^{-1})_{00}\ll1$ by Lemma~\ref{lem:oldspace-cancellation}.
For the coefficient sums below, fix $\theta\geq0$ and
$B>2+\theta+\eps$.
Applying the divisor bound separately to
$r\leq\sqrt p$ and $r>\sqrt p$ gives
\begin{equation}
\begin{split}
\sum_{r,m,n\geq1}&~\frac{(mn)^\theta}{r\sqrt{mn}}
 \left(1+\frac{r^2mn}{p}\right)^{-B}\\
\ll_{B,\theta,\eps} &~p^{1/2+\theta+\eps}
 \sum_{r\leq\sqrt p}r^{-2-2\theta-2\eps}
 +p^B\sum_{r>\sqrt p}r^{-1-2B}\\
\ll_{B,\theta,\eps} &~p^{1/2+\theta+\eps}.
\label{eq:second-coefficient-sum}
\end{split}
\end{equation}

\begin{proposition}
\label{prop:second-level-one-remainders}
Uniformly in the twist range, we have
\begin{equation*}
\begin{aligned}
\mathcal R_{\rm old}(a,b)\ll_{h,\eta,\eps}
 &~(ab)^{5/14+\eps}p^{6/7+\eps},\\
\mathcal R_{\min}(a,b)\ll_{h,\eta,\eps}
 &~(ab)^\eps p^{1/2+\eps},
\end{aligned}
\end{equation*}
 and
 \begin{equation*}
\mathcal R_{\max,\rm old}(a,b)\ll_{h,\eta,\eps}
 (ab)^{7/64+\eps}p^{39/64+\eps}.
\end{equation*}
\end{proposition}

\begin{proof}
By \eqref{eq:level-one-normalizer-bounds} and
\eqref{eq:level-one-counting-bound}, rapid decrease of $h$ makes
the level one spectral sums and integrals with any fixed polynomial
weight bounded independently of $p$.
The twist contributes $(ab)^{\theta+\eps}$, while
\eqref{eq:product-weight-decay} and \eqref{eq:second-coefficient-sum}
bound the remaining index sums by $p^{1/2+\theta+\eps}$. Taking $\theta=5/14,0,7/64$ for the old
cuspidal, minimal and maximal old spectra proves the three estimates.
\end{proof}

\subsubsection{The maximal Eisenstein new spectrum}
\label{subsec:second-maximal-new}

\begin{proposition}
\label{prop:second-maximal-remainder}
Uniformly in the twist range,
\begin{equation*}
\mathcal R_{\max,\rm new}(a,b)\ll_{h,Q,\eta,\eps}
 p^{3/2+\eps}.
\end{equation*}
\end{proposition}

\begin{proof}
Apply \eqref{eq:gl2-induced-coefficient-average} in
\eqref{eq:second-max-new-def}. For sufficiently large $p$, we have
$ab\leq Q^2p^{2\eta}<p$. Together with $\chi_p\ll1$ and
$1+mn/p\leq1+r^2mn/p$, \eqref{eq:product-weight-decay} gives
\begin{equation*}
\begin{aligned}
|\mathcal R_{\max,\rm new}(a,b)|
&~\ll_{h,\eta,B,\eps} p^{1+\eps}(ab)^\eps
 \sum_{r,m,n\geq1}\frac{(mn)^\eps}{r\sqrt{mn}}
 \left(1+\frac{r^2mn}{p}\right)^{-B+1}\\
&~\ll_{h,Q,\eta,\eps} p^{3/2+\eps}.
\end{aligned}
\end{equation*}
The last step follows from \eqref{eq:second-coefficient-sum} with
$\theta=\eps$ and $B$ sufficiently large.
\end{proof}

\subsubsection{Proof of Theorem~\ref{thm:harmonic-second-moment}}
Combine \eqref{eq:ordinary-diagonal},
Lemmas~\ref{lem:second-moment-short} and \ref{lem:second-moment-long},
and Propositions~\ref{prop:second-level-one-remainders} and
\ref{prop:second-maximal-remainder} in
\eqref{eq:second-moment-decomposition}.
For sufficiently large $p$, $ab\leq Q^2p^{2\eta}<p$.
The level one remainders are
$O_{h,Q,\eta,\eps}(p^{17/14+\eps})$ and are absorbed by
$p^{3/2+\eps}$. Dividing by $\chi_p=1+O(p^{-1})$ changes the
diagonal main term by $O_{h,Q,\eta}(p\log p/\sqrt{ab})$, which is
absorbed by $p^{7/4}(ab)^{-1/4}$. This proves
\eqref{eq:twisted-second-moment}.

For the uniform asymptotic, $ab\leq Q^2p^{2\eta}$ gives
\begin{equation*}
\frac{\sqrt{ab}}{p^2}
\big(p^{7/4}(ab)^{-1/4}+p^{3/2+\eps}\big)
\leq Q^{1/2}p^{-1/4+\eta/2}+Qp^{-1/2+\eta+\eps}
=o_{Q,\eta}(1).
\end{equation*}
Here $\eta<1/2$ and $\eps$ is arbitrarily small.
Since $\mathcal H_h>0$ and
\begin{equation*}
\log\frac p{ab}\geq (1-2\eta)\log p-2\log Q
\gg_{Q,\eta}\log p,
\end{equation*}
the constant term and the finite Euler corrections have relative
size $O_{h,Q,\eta}(1/\log p)$. The stated asymptotic follows.
Taking $Q=1$ and $a=b=1$ gives \eqref{eq:ordinary-second}.

\section{Mollified harmonic moments}
\label{sec:mollified-moments}

We evaluate the mollified harmonic first and second moments using
the uniform twisted moment formulas of
\S\,\ref{sec:harmonic-twisted-moments}. We also compute the holomorphic
mixed moment associated with $\widehat M_\pi$ in
\eqref{eq:two-piece-mollifier}. These estimates will be combined in
\S\,\ref{sec:harmonic-nonvanishing} to prove the non-vanishing theorem.

\subsection{Mollified harmonic first moment}
\label{subsec:mollified-first}

Recall the mollified harmonic first moment
$\mathcal M_1(p)$ defined in \eqref{eq:intro-mollified-first-moment}.
\begin{theorem}
\label{prop:mollified-first}
Fix a square-free integer $Q\geq1$ and $0<\eta<1/2$, independently
of $p$, and put $L=p^\eta$. As $p\to\infty$ through primes
$p\nmid Q$, we have
\begin{equation}
\begin{split}
\mathcal M_1(p)= &~P\mathcal H_h
 +O_{h,Q,\eta,\eps}\!\left(p^{5/4+\eps}L^{1/2}\right)\\
= &~P\mathcal H_h+O_{h,Q,\eta}(p^{3/2}).
\end{split}
\label{eq:mollified-first-result}
\end{equation}
\end{theorem}
\begin{proof}
The finite product in \eqref{eq:intro-mollifier} expands as
\begin{equation*}
\prod_{q\mid Q}
\left(1-\frac{A_\pi(1,q)}{\sqrt q}
+\frac{A_\pi(q,1)}q-\frac1{q^{3/2}}\right)
= \sum_{uvw\mid Q}
 \frac{\mu(u)\mu(w)}{\sqrt u\,v\,w^{3/2}}A_\pi(v,u),
\end{equation*}
where $u,v,w$ are positive and pairwise coprime.
Since $p\nmid Q$, $(\ell,Q)=1$ on the support of $x_\ell$, and
$\ell\leq L<p$, multiplicativity and
\eqref{eq:unramified-unitarity} give
$A_\pi(v,u)A_\pi(1,\ell)=A_\pi(v,u\ell)
=\overline{A_\pi(u\ell,v)}$.
Recalling the definition of $M_h$ in \eqref{eq:def-Mrs}, we obtain
\begin{equation}
\mathcal M_1(p)= \sum_{uvw\mid Q}
 \frac{\mu(u)\mu(w)}{\sqrt u\,v\,w^{3/2}}
 \sum_{\ell\leq L}\frac{x_\ell}{\sqrt\ell}M_h(u\ell,v,p).
\label{eq:mollified-first-afe}
\end{equation}
Here $u\ell\leq Qp^\eta$ and $v\leq Q$. So Theorem \ref{thm:harmonic-first-moment} applies uniformly to every term in \eqref{eq:mollified-first-afe}. Its main term requires $u\ell=1$, hence $u=\ell=1$.
Since $x_1=1$, it follows that
\begin{equation*}
\mathcal M_1(p)= P\mathcal H_h
 \sum_{vw\mid Q}\frac{\mu(w)}{(vw)^{3/2}}+O_{h,Q,\eta,\eps}\!\left(
 p^{5/4+\eps}\sum_{\ell\leq L}\frac{|x_\ell|}{\sqrt\ell}\right).
\end{equation*}
Finally,
\begin{equation*}
\sum_{vw\mid Q}\frac{\mu(w)}{(vw)^{3/2}}
= \prod_{q\mid Q}(1+q^{-3/2}-q^{-3/2})=1,
\end{equation*}
and $|x_\ell|\leq1$ gives
$\sum_{\ell\leq L}|x_\ell|/\sqrt\ell\ll L^{1/2}$.
This proves the first estimate in \eqref{eq:mollified-first-result}.
Since $5/4+\eta/2<3/2$, the arbitrarily small $\eps$ gives the
second estimate.
\end{proof}

\subsection{Mollified harmonic second moment}
\label{subsec:mollified-second}
Recall the mollified harmonic second moment
$\mathcal M_2(p)$ defined in \eqref{eq:intro-mollified-second-moment}.
\begin{theorem}
\label{prop:mollified-second}
Fix a square-free integer $Q\geq1$ and $0<\eta<1/2$, and put
$L=p^\eta$. As $p\to\infty$ through primes $p\nmid Q$, we have
\begin{equation}
\mathcal M_2(p)= \zeta(2)\prod_{q\mid Q}(1-q^{-2})
 (1+\eta^{-1})P\mathcal H_h+O_{h,Q,\eta}\left(
 \frac{p^2\log(2+\log p)}{\log p}\right).
\label{eq:mollified-second-asymptotic}
\end{equation}
\end{theorem}

\subsubsection{Hecke expansion and the twisted second moment}

The finite Euler factors in \eqref{eq:intro-mollifier} have the exact
Hecke expansion
\begin{equation}
\prod_{q\mid Q}
 \left|1-\frac{A_\pi(1,q)}{\sqrt q}
 +\frac{A_\pi(q,1)}q-\frac1{q^{3/2}}\right|^2
=\sum_{uv\mid Q^2}c_Q(u,v)\overline{A_\pi(u,v)}.
\label{eq:finite-euler-square}
\end{equation}
Here $c_Q(u,v)=\prod_{q\mid Q}c_q(\operatorname{ord}_q u,
\operatorname{ord}_q v)$ for $uv\mid Q^2$, where
\begin{equation*}
c_q(i,j):= 
\begin{cases}
(1-q^{-3/2})^2+q^{-1}+q^{-2},&(i,j)=(0,0),\\
(1-q^{-3/2})(q^{-1}-q^{-1/2})-q^{-3/2},
 &(i,j)=(1,0),(0,1),\\
q^{-1}+q^{-2},&(i,j)=(1,1),\\
-q^{-3/2},&(i,j)=(2,0),(0,2),\\
0,&\text{otherwise}.
\end{cases}
\end{equation*}
These coefficients follow from the unramified Hecke relations
$A_\pi(1,q)A_\pi(q,1)=A_\pi(q,q)+1$ and
$A_\pi(1,q)^2=A_\pi(1,q^2)+A_\pi(q,1)$, together with their conjugates.
In particular,
\begin{equation}
\sum_{uv\mid Q^2}\frac{c_Q(u,v)}{\sqrt{uv}}
=\prod_{q\mid Q}\left(
 c_q(0,0)+\frac{2c_q(1,0)}{\sqrt q}
 +\frac{c_q(1,1)}q+\frac{2c_q(2,0)}q\right)
=\prod_{q\mid Q}(1-q^{-1}).
\label{eq:finite-euler-diagonal-sum}
\end{equation}
Recall that the coefficients $x_\ell$ are real and supported on
$\ell\leq L$ with $(\ell,Q)=1$. The unramified Hecke relation
and \eqref{eq:unramified-unitarity} give
\begin{equation*}
\left|\sum_{\ell\leq L}
 \frac{x_\ell A_\pi(1,\ell)}{\sqrt\ell}\right|^2
= \sum_{d,\ell_1,\ell_2\geq1}
 \frac{x_{d\ell_1}x_{d\ell_2}}{d\sqrt{\ell_1\ell_2}}
 \overline{A_\pi(\ell_1,\ell_2)}.
\end{equation*}
Here $d\ell_1,d\ell_2\leq L$ are the original mollifier indices.
Since $(\ell_1\ell_2,uv)=1$ in every non-zero term, multiplicativity and
\eqref{eq:finite-euler-square} yield
\begin{equation}
\mathcal M_2(p)= \sum_{d,\ell_1,\ell_2\geq1}
 \frac{x_{d\ell_1}x_{d\ell_2}}{d\sqrt{\ell_1\ell_2}}
 \sum_{uv\mid Q^2}c_Q(u,v)\mathcal N_2(u\ell_1,v\ell_2,p).
\label{eq:trace-expansion}
\end{equation}
Every non-zero term in \eqref{eq:trace-expansion} satisfies
$uv\ell_1\ell_2\leq Q^2L^2=Q^2p^{2\eta}$ and
$(uv\ell_1\ell_2,p)=1$.
Thus Theorem~\ref{thm:harmonic-second-moment} applies uniformly,
with the same fixed $Q,\eta$ and $h$.
Substituting \eqref{eq:twisted-second-moment} into
\eqref{eq:trace-expansion}, we obtain
\begin{equation}
\begin{split}
\label{M2}
\mathcal M_2(p)= &~P\sum_{d,\ell_1,\ell_2\geq1}
 \frac{x_{d\ell_1}x_{d\ell_2}}{d\ell_1\ell_2}
 \sum_{uv\mid Q^2}\frac{c_Q(u,v)}{\sqrt{uv}}\\*
&~\quad\times\Biggl((1-p^{-1})
 \left\{\mathcal H_h\left(\log\frac p{uv\ell_1\ell_2}+2\gamma\right)
 +\mathcal K_h\right\}
 +\frac{2\log p}{p}\mathcal H_h\Biggr)\\
&~\qquad{}+O_{h,Q,\eta,\eps}\Biggl(
 p^{7/4}\sum_{d,\ell_1,\ell_2\geq1}
 \frac{|x_{d\ell_1}x_{d\ell_2}|}{d(\ell_1\ell_2)^{3/4}}+
 p^{3/2+\eps}\sum_{d,\ell_1,\ell_2\geq1}
 \frac{|x_{d\ell_1}x_{d\ell_2}|}{d\sqrt{\ell_1\ell_2}}\Biggr).
 \end{split}
\end{equation}
In the error terms, the $u,v$-sums are absorbed into the
implied constant, since $Q$ is fixed.

For $0<\beta<1$, the coefficients multiplying the twisted moments satisfy
\begin{equation}
\sum_{d,\ell_1,\ell_2\geq1}
 \frac{|x_{d\ell_1}x_{d\ell_2}|}{d(\ell_1\ell_2)^\beta}
\leq  \sum_{d\leq L}\frac1d
 \left(\sum_{\ell_1\leq L/d}\ell_1^{-\beta}\right)^2
\ll_\beta  L^{2(1-\beta)}
 \sum_{d\geq1}d^{-3+2\beta}
\ll_\beta L^{2(1-\beta)}.
\label{eq:mollified-second-outside-sum}
\end{equation}
Taking $\beta=3/4$ and $1/2$ in
\eqref{eq:mollified-second-outside-sum} bounds the total error in \eqref{M2} by
\begin{equation*}
O_{h,Q,\eta,\eps}\!\left(p^{7/4}L^{1/2}+p^{3/2+\eps}L\right).
\end{equation*}

\subsubsection{Arithmetic quadratic forms}

It remains to evaluate the arithmetic main term, which is expressed through
\begin{equation}
\begin{split}
Q_0(L):= &~\sum_{d,\ell_1,\ell_2\geq1}
 \frac{x_{d\ell_1}x_{d\ell_2}}{d\ell_1\ell_2},\\
Q_1(L):= &~-\sum_{d,\ell_1,\ell_2\geq1}
 \frac{x_{d\ell_1}x_{d\ell_2}\log(\ell_1\ell_2)}{d\ell_1\ell_2}.
\end{split}
\label{eq:quadratic-forms}
\end{equation}
All non-zero terms have $(d\ell_1\ell_2,Q)=1$.
Using \eqref{eq:finite-euler-diagonal-sum} and \eqref{eq:quadratic-forms}, we therefore obtain
\begin{equation}
\begin{split}
\mathcal M_2(p)= &~P(1-p^{-1})\prod_{q\mid Q}(1-q^{-1})
\left\{\mathcal H_h\big((\log p+2\gamma)Q_0(L)+Q_1(L)\big)
 +\mathcal K_hQ_0(L)\right\}\\
&~-P(1-p^{-1})\mathcal H_hQ_0(L)
 \sum_{uv\mid Q^2}\frac{c_Q(u,v)\log(uv)}{\sqrt{uv}}\\
&~+\frac{2P\log p}{p}\mathcal H_hQ_0(L)
 \prod_{q\mid Q}(1-q^{-1})\\
&~+O_{h,Q,\eta,\varepsilon}
 \left(p^{7/4}L^{1/2}+p^{3/2+\varepsilon}L\right).
\end{split}
\label{eq:second-moment-result}
\end{equation}

\begin{lemma}
\label{lem:mobius-quadratic-forms}
Let $x_\ell$ be in \eqref{eq:mollifier-coefficients} and the quadratic forms
in \eqref{eq:quadratic-forms}.
For fixed square-free $Q$ and uniformly for $L\geq2$, we have
\begin{equation}
\begin{split}
Q_0(L)= &~\frac{\zeta(2)}{\log L}
 \prod_{q\mid Q}(1+q^{-1})
 +O_Q\left(\frac{\log(2+\log L)}{(\log L)^2}\right),\\
Q_1(L)= &~\zeta(2)\prod_{q\mid Q}(1+q^{-1})
 +O_Q\left(\frac{\log(2+\log L)}{\log L}\right).
\end{split}
\label{eq:quadratic-form-asymptotics}
\end{equation}
In particular, $Q_0(L)\geq0$, $Q_0(L)\ll_Q1/\log L$, and
$|Q_1(L)|\ll_Q1$.
\end{lemma}

\begin{proof}
Let
\begin{equation*}
F(X):= \sum_{n\leq X}\frac{\mu(n)}n\log\frac Xn.
\end{equation*}
Perron's formula \cite[Proposition~5.54]{IK} expresses $F(X)$ as the integral of
$X^s/(s^2\zeta(1+s))$. Shifting the contour using the zero-free
region and the bound of $1/\zeta(s)$ in \cite[Theorem~8.29]{IK}, and taking
the residue at $s=0$, gives, for some $c>0$,
\begin{equation}
F(X)= 1+O\big(e^{-c\sqrt{\log X}}\big),\qquad
F(X)\ll1\quad(X\geq1).
\label{eq:mobius-log-asymptotic}
\end{equation}
For $d\geq1$ put
\begin{equation}
\begin{split}
F_d(X)\colon \!\!= &~\sum_{\substack{\ell_1\leq X\\(\ell_1,dQ)=1}}
 \dfrac{\mu(\ell_1)}{\ell_1}
 \prod_{q\mid\ell_1}\frac q{q+1}\log\frac X{\ell_1},\\
T_d(X)\colon \!\!= &~\sum_{\substack{\ell_1\leq X\\(\ell_1,dQ)=1}}
 \dfrac{\mu(\ell_1)\log\ell_1}{\ell_1}
 \prod_{q\mid\ell_1}\frac q{q+1}\log\frac X{\ell_1}.
\end{split}
\label{eq:coprime-mobius-definitions}
\end{equation}
Define the non-negative multiplicative function $k_d$ by $k_d(1)=1$
and, for primes $q$ and $j\geq1$,
\begin{equation*}
k_d(q^j):= 
\begin{cases}
1,&q\mid dQ,\\
(q+1)^{-1},&q\nmid dQ.
\end{cases}
\end{equation*}
By the Euler product, a direct calculation gives
\begin{equation}
F_d(X)=\sum_{e\leq X}\frac{k_d(e)}eF(X/e),
\label{eq:coprime-mobius-convolution}
\end{equation}
and
\begin{equation}
\sum_{e\geq1}\frac{k_d(e)}e=\zeta(2)\prod_{q\mid dQ}(1+q^{-1}).
\label{eq:coprime-euler-mass}
\end{equation}
The Euler product for $\sum_{e\geq1}k_d(e)e^{-s}$ converges
absolutely for $\Re(s)>0$. Logarithmic differentiation at $s=1$
yields
\begin{equation}
\begin{split}
\sum_{e\geq1}\frac{k_d(e)\log e}e
= &~\zeta(2)\prod_{q\mid dQ}(1+q^{-1})
 \left(\sum_{q\mid dQ}\frac{\log q}{q-1}
 +\sum_{q\nmid dQ}\frac{\log q}{q(q-1)}\right)\\
\ll &~\prod_{q\mid dQ}(1+q^{-1})
 \left(1+\sum_{q\mid dQ}\frac{\log q}{q-1}\right).
\end{split}
\label{eq:coprime-euler-log-sum}
\end{equation}
Since $k_d(e)\geq0$, \eqref{eq:mobius-log-asymptotic}--\eqref{eq:coprime-mobius-convolution}
also imply
$F_d(X)\ll\prod_{q\mid dQ}(1+q^{-1})$.

For $X\geq2$, we subtract the sum in \eqref{eq:coprime-euler-mass}
from \eqref{eq:coprime-mobius-convolution}, split at $e=\sqrt X$, and obtain
\begin{equation*}
F_d(X)-\sum_{e\geq1}\frac{k_d(e)}e
= \sum_{e\leq\sqrt X}\frac{k_d(e)}e\big(F(X/e)-1\big)+\sum_{\sqrt X<e\leq X}\frac{k_d(e)}eF(X/e)
 -\sum_{e>\sqrt X}\frac{k_d(e)}e.
\end{equation*}
By \eqref{eq:mobius-log-asymptotic},
$|F(X/e)-1|\ll\exp(-c\sqrt{(\log X)/2})$ for $e\leq\sqrt X$
and $|F(X/e)|\ll1$ for $\sqrt X<e\leq X$.
Since $k_d(e)\geq0$, it follows that
\begin{equation*}
\left|F_d(X)-\sum_{e\geq1}\frac{k_d(e)}e\right|
\ll \exp(-c\sqrt{(\log X)/2})
 \sum_{e\geq1}\frac{k_d(e)}e
 +\frac2{\log X}\sum_{e\geq1}\frac{k_d(e)\log e}e.
\end{equation*}
Inserting \eqref{eq:coprime-euler-mass} and
\eqref{eq:coprime-euler-log-sum}, and using
$\exp(-c\sqrt{(\log X)/2})\ll(1+\log X)^{-1}$, we obtain
\begin{equation}
F_d(X)= \zeta(2)\prod_{q\mid dQ}(1+q^{-1})+O\left(\frac{\prod_{q\mid dQ}(1+q^{-1})}{1+\log X}
 \left(1+\sum_{q\mid dQ}\frac{\log q}{q-1}\right)\right).
\label{eq:coprime-mobius-asymptotic}
\end{equation}
For $1\leq X<2$, we have $F_d(X)=\log X$, so
\eqref{eq:coprime-mobius-asymptotic} holds for all $X\geq1$.
By the definitions of $F_d$ and $T_d$ in
\eqref{eq:coprime-mobius-definitions}, 
and then substituting \eqref{eq:coprime-mobius-asymptotic} for $F_d$, we have
\begin{equation}
\begin{split}
T_d(X)&=(\log X)\,F_d(X)-2\int_1^X F_d(t)\frac{\dd t}{t}\\
&=-\zeta(2)(\log X)\prod_{q\mid dQ}(1+q^{-1})+O\left(\log(2+\log X)\prod_{q\mid dQ}(1+q^{-1})
 \left(1+\sum_{q\mid dQ}\frac{\log q}{q-1}\right)\right).
\label{eq:coprime-log-insertion-asymptotic}
\end{split}
\end{equation}

Since $x_\ell$ is supported on square-free integers,
$x_{d\ell_1}x_{d\ell_2}=0$ unless $(d,\ell_1\ell_2)=1$.
Substituting \eqref{eq:mollifier-coefficients} into
\eqref{eq:quadratic-forms} and using
\eqref{eq:coprime-mobius-definitions}, we obtain
\begin{equation}
\begin{split}
Q_0(L)= &~\frac1{(\log L)^2}
 \sum_{\substack{d\leq L\\(d,Q)=1}}\frac{\mu(d)^2}{d}
 \prod_{q\mid d}\left(\frac q{q+1}\right)^2F_d(L/d)^2,\\
Q_1(L)= &~-\frac2{(\log L)^2}
 \sum_{\substack{d\leq L\\(d,Q)=1}}\frac{\mu(d)^2}{d}
 \prod_{q\mid d}\left(\frac q{q+1}\right)^2F_d(L/d)T_d(L/d).
\end{split}
\label{eq:quadratic-forms-expanded}
\end{equation}
The factor $2$ in the second identity follows from
$\log(\ell_1\ell_2)=\log\ell_1+\log\ell_2$ and symmetry in
$\ell_1,\ell_2$. The first identity also shows that $Q_0(L)\geq0$.

It remains to evaluate the main terms and bound the errors in
\eqref{eq:quadratic-forms-expanded}. For $(d,Q)=1$,
\begin{equation*}
\prod_{q\mid dQ}(1+q^{-1})\prod_{q\mid d}\frac q{q+1}
= \prod_{q\mid Q}(1+q^{-1}).
\end{equation*}
Substituting \eqref{eq:coprime-mobius-asymptotic} and
\eqref{eq:coprime-log-insertion-asymptotic} into
\eqref{eq:quadratic-forms-expanded}, and using
$|F_d(X)|\ll\prod_{q\mid dQ}(1+q^{-1})$, gives
\begin{equation}
Q_0(L)= \frac{\zeta(2)^2}{(\log L)^2}
 \prod_{q\mid Q}(1+q^{-1})^2
 \sum_{\substack{d\leq L\\(d,Q)=1}}\frac{\mu(d)^2}d+O_Q\!\left(\frac1{(\log L)^2}
 \sum_{d\leq L}\frac{1+\sum_{q\mid d}\frac{\log q}{q-1}}
 {d(1+\log(L/d))}\right),
\label{eq:quadratic-form-q0-main-error}
\end{equation}
 and
 \begin{equation}
\begin{split}
Q_1(L)= &~\frac{2\zeta(2)^2}{(\log L)^2}
 \prod_{q\mid Q}(1+q^{-1})^2
 \sum_{\substack{d\leq L\\(d,Q)=1}}\frac{\mu(d)^2}d\log(L/d)\\
&~+O_Q\!\left(\frac1{(\log L)^2}
 \sum_{d\leq L}\frac{1+\sum_{q\mid d}\frac{\log q}{q-1}}d
 \log(2+\log(L/d))\right).
\end{split}
\label{eq:quadratic-forms-main-error}
\end{equation}
In these error terms, we used $\mu(d)^2\leq1$ and
\begin{equation*}
1+\sum_{q\mid dQ}\frac{\log q}{q-1}
\ll_Q 1+\sum_{q\mid d}\frac{\log q}{q-1}.
\end{equation*}

We first evaluate the sums in the main terms of
\eqref{eq:quadratic-form-q0-main-error} and
\eqref{eq:quadratic-forms-main-error}. Using
$\mathbf1_{(m,Q)=1}=\sum_{e\mid(m,Q)}\mu(e)$, we have, for $y\geq1$,
\begin{equation*}
\sum_{\substack{m\leq y\\(m,Q)=1}}\frac1m
= \sum_{e\mid Q}\frac{\mu(e)}e
 \sum_{n\leq y/e}\frac1n
=\prod_{q\mid Q}(1-q^{-1})\log y+O_Q(1).
\end{equation*}
Since $\mu(d)^2=\sum_{e^2\mid d}\mu(e)$, the above gives
\begin{equation}
\sum_{\substack{d\leq x\\(d,Q)=1}}\frac{\mu(d)^2}d
=\sum_{\substack{e\leq\sqrt x\\(e,Q)=1}}\frac{\mu(e)}{e^2}
 \sum_{\substack{m\leq x/e^2\\(m,Q)=1}}\frac1m
= \frac{\log x}{\zeta(2)}
 \prod_{q\mid Q}(1+q^{-1})^{-1}+O_Q(1).
\label{eq:mollifier-divisor-asymptotic}
\end{equation}
This evaluates the sum in the main term of $Q_0(L)$.
For $Q_1(L)$, the additional factor $\log(L/d)$ is handled by
$\log(L/d)=\int_d^L\dd t/t$. Thus
\begin{equation}
\sum_{\substack{d\leq L\\(d,Q)=1}}\frac{\mu(d)^2}d\log(L/d)
= \int_1^L
 \left(\sum_{\substack{d\leq t\\(d,Q)=1}}\frac{\mu(d)^2}d\right)
 \frac{\dd t}{t}
= \frac{(\log L)^2}{2\zeta(2)}
 \prod_{q\mid Q}(1+q^{-1})^{-1}+O_Q(\log L).
\label{eq:mollifier-log-divisor-asymptotic}
\end{equation}
Substituting \eqref{eq:mollifier-divisor-asymptotic} and
\eqref{eq:mollifier-log-divisor-asymptotic} into the main terms of
\eqref{eq:quadratic-form-q0-main-error} and
\eqref{eq:quadratic-forms-main-error} gives
$\zeta(2)\prod_{q\mid Q}(1+q^{-1})/\log L$ for $Q_0(L)$ and
$\zeta(2)\prod_{q\mid Q}(1+q^{-1})$ for $Q_1(L)$, with respective
errors $O_Q((\log L)^{-2})$ and $O_Q((\log L)^{-1})$.

We now bound the two explicit error sums in
\eqref{eq:quadratic-form-q0-main-error} and
\eqref{eq:quadratic-forms-main-error}. To control the factor
$1+\sum_{q\mid d}(\log q)/(q-1)$ in both errors,
we first sum it over $d\leq y$, where $y\geq1$:
\begin{equation*}
\sum_{d\leq y}\left(1+\sum_{q\mid d}\frac{\log q}{q-1}\right)
= \lfloor y\rfloor+
 \sum_{q\leq y}\frac{\log q}{q-1}\left\lfloor\frac yq\right\rfloor
\leq y\left(1+\sum_q\frac{\log q}{q(q-1)}\right)\ll y.
\end{equation*}
For the error in $Q_0(L)$, we apply Abel's summation formula and use
\begin{equation*}
\int_1^L\frac{\dd t}{t(1+\log(L/t))}
= \log(1+\log L).
\end{equation*}
For the error in $Q_1(L)$, we use
$\log(2+\log(L/d))\leq\log(2+\log L)$ and
\begin{equation*}
\sum_{d\leq L}\frac1d
 \left(1+\sum_{q\mid d}\frac{\log q}{q-1}\right)
\ll 1+\int_1^L\frac{\dd t}{t}\ll\log L.
\end{equation*}
Consequently, the errors in $Q_0(L)$ and $Q_1(L)$ satisfy,
respectively,
\begin{equation*}
\frac1{(\log L)^2}
 \sum_{d\leq L}\frac{1+\sum_{q\mid d}\frac{\log q}{q-1}}
 {d(1+\log(L/d))}
 \ll\frac{\log(2+\log L)}{(\log L)^2},
\end{equation*}
 and
 \begin{equation*}
\frac1{(\log L)^2}
 \sum_{d\leq L}\frac{1+\sum_{q\mid d}\frac{\log q}{q-1}}d
 \log(2+\log(L/d))
 \ll\frac{\log(2+\log L)}{\log L}.
\end{equation*}
Combining these bounds with the main terms above in
\eqref{eq:quadratic-form-q0-main-error} and
\eqref{eq:quadratic-forms-main-error} proves
\eqref{eq:quadratic-form-asymptotics}.
\end{proof}

For any real coefficients supported on $(\ell,Q)=1$ and
$\ell\leq L$, the definition gives
\begin{equation*}
Q_0(L)= \sum_{\substack{d\leq L\\(d,Q)=1}}d
 \left(\sum_{\substack{n\leq L\\d\mid n}}\frac{x_n}{n}\right)^2.
\end{equation*}
If $x_1=1$, the M\"obius inversion formula and the Cauchy--Schwarz inequality imply
\begin{equation*}
1=\left(\sum_{\substack{d\leq L\\(d,Q)=1}}\mu(d)
 \sum_{\substack{n\leq L\\d\mid n}}\frac{x_n}{n}\right)^2
\leq\left(\sum_{\substack{d\leq L\\(d,Q)=1}}
 \frac{\mu(d)^2}{d}\right)Q_0(L).
\end{equation*}
For fixed $Q$, \eqref{eq:mollifier-divisor-asymptotic} therefore gives
\begin{equation*}
Q_0(L)\geq \left(\sum_{\substack{d\leq L\\(d,Q)=1}}
 \frac{\mu(d)^2}{d}\right)^{-1}
=\frac{\zeta(2)}{\log L}\prod_{q\mid Q}(1+q^{-1})
 \left(1+O_Q\left(\frac1{\log L}\right)\right)
 \qquad(L\to\infty).
\end{equation*}
By \eqref{eq:quadratic-form-asymptotics}, the coefficients in
\eqref{eq:mollifier-coefficients} attain this lower bound asymptotically.
Thus they are asymptotically optimal for $Q_0(L)$ among all real
coefficients with $x_1=1$ supported on $\ell\leq L$ and $(\ell,Q)=1$.

\subsubsection{Proof of Theorem~\ref{prop:mollified-second}}

Since $\log L=\eta\log p$, Lemma~\ref{lem:mobius-quadratic-forms} gives
\begin{equation*}
\begin{aligned}
&~\prod_{q\mid Q}(1-q^{-1})
 \big((\log p+2\gamma)Q_0(L)+Q_1(L)\big)\\
= &~\zeta(2)\prod_{q\mid Q}(1-q^{-2})(1+\eta^{-1})
 +O_{Q,\eta}\left(\frac{\log(2+\log p)}{\log p}\right).
\end{aligned}
\end{equation*}
Since $Q$ is fixed and $Q_0(L)\ll_{Q,\eta}1/\log p$, the
terms involving $\log(uv)$ and $\mathcal K_h$ in
\eqref{eq:second-moment-result} are $O_{h,Q,\eta}(p^2/\log p)$.
Replacing $1-p^{-1}$ by $1$ introduces an error
$O_{h,Q,\eta}(p)$, and the term containing $2\log p/p$
satisfies the same bound.
Moreover, $L=p^\eta$ and $\eta<1/2$ give
$p^{7/4}L^{1/2}+p^{3/2+\eps}L\ll_\eta p^2/\log p$,
since $\eps$ is arbitrarily small.
Substituting these estimates into \eqref{eq:second-moment-result}
proves \eqref{eq:mollified-second-asymptotic}.

\subsection{The holomorphic mixed moment}
\label{subsec:mollified-holomorphic-square}
Recall the holomorphic mixed moment
$\mathcal M_{1,1}(p)$ defined in \eqref{eq:intro-mollified-mixed-moment}.
 We use the same $h$, $Q$, $\eta$, $L$ and the mollifier $M_\pi$ defined above.

\begin{theorem}
\label{thm:mollified-holomorphic-square}
For every fixed square-free $Q$ and $0<\eta<1/2$,
as $p\to\infty$ through primes $p\nmid Q$, we have
\begin{equation}
\mathcal M_{1,1}(p)= P\mathcal H_h
 \prod_{q\nmid Q}
 \left(1-\frac{2}{\sqrt q(q+1)^2}
             +\frac{1}{q(q+1)^2}\right)
 +O_{h,Q,\eta}\left(\frac{p^2}{\log p}\right).
\label{eq:mollified-holomorphic-square}
\end{equation}
The Euler product converges absolutely, lies in $(0,1]$, and tends
to $1$ as $Q$ ranges over products of all primes up to a bound
tending to infinity.
\end{theorem}

\subsubsection{Application of the Kuznetsov formula}
By \eqref{eq:asymmetric-holomorphic-square-afe}, the moment takes the form
\begin{equation}
\begin{split}
\mathcal M_{1,1}(p)=
&~c_{\rm cusp}\sum_{\pi\in\mathcal B(p)}
 \frac{M_\pi^2}{L(1,\pi,\mathrm{Ad})}
 \sum_{m,n\geq1}\frac{A_\pi(1,m)A_\pi(1,n)}{\sqrt{mn}}
 K_{\mu_\pi}^+(mn/(pL^2))\\*
&~+c_{\rm cusp}\sum_{\pi\in\mathcal B(p)}
 \frac{\epsilon(\pi)^2M_\pi^2}{L(1,\pi,\mathrm{Ad})}
 \sum_{m,n\geq1}\frac{\overline{A_\pi(1,m)A_\pi(1,n)}}{\sqrt{mn}}
 K_{\mu_\pi}^-(mnL^2/p).
\end{split}
\label{eq:holomorphic-moment-afe}
\end{equation}
After fixing $Q$ and $\eta$, choose the kernel decay exponents needed below, then take the auxiliary $\eps>0$ sufficiently small in terms of $\eta$ and those exponents, and finally fix $a_0,J$ as required by Lemma~\ref{lem:holomorphic-square-weights}. In particular, require $0<\eps<\min\{2\eta,1/2-\eta\}$. All these choices precede $p\to\infty$. This auxiliary $\eps$ is independent of the prescribed error in the final non-vanishing proportion.
We first justify truncating the two Dirichlet sums in
\eqref{eq:holomorphic-moment-afe}.
The finite Euler factors in \eqref{eq:intro-mollifier} are
$O_Q(1)$ by the coefficient bounds of \S\,\ref{prelim}.
The Cauchy--Schwarz inequality and
\eqref{eq:coefficient-square-average}, applied to $A_\pi(1,\ell)$
with $\ell\leq L<p^{1/2}$, give, for every fixed $C\geq0$,
\begin{equation*}
c_{\rm cusp}\sum_{\pi\in\mathcal B(p)}
 \frac{h(\mu_\pi)(1+\|\mu_\pi\|)^C|M_\pi|^2}
 {L(1,\pi,\mathrm{Ad})}
\ll_{h,Q,C} p^2
 \left(\sum_{\ell\leq L}\frac{|x_\ell|}{\sqrt\ell}\right)^2
\ll_{h,Q,C}p^2L.
\end{equation*}
Here we used $|x_\ell|\leq1$.
Combining \eqref{eq:holomorphic-square-weight-decay} with the
coefficient bounds of \S\,\ref{prelim} and the preceding average,
we may truncate the two Dirichlet sums in
\eqref{eq:holomorphic-moment-afe} at
\begin{equation*}
mn\leq p^{1+\eps}L^2,\qquad mn\leq p^{1+\eps}/L^2,
\end{equation*}
respectively, with total error $O_{h,Q,\eta,\eps,A}(p^{-A})$
for every fixed $A>0$. The same bounds give absolute convergence
and justify interchanging the Dirichlet sums with the cuspidal sum.
Since $L=p^\eta$ and $\eps<2\eta$, the truncated
dual range satisfies
\begin{equation*}
mn\leq p^{1+\eps}/L^2=p^{1-2\eta+\eps}<p.
\end{equation*}
Hence $(mn,p)=1$ throughout the truncated dual sum.
In the first sum, we retain the terms with $p\mid mn$.
For example, $mn=p$ gives the weight argument $L^{-2}<1$,
so \eqref{eq:holomorphic-square-weight-decay} gives no decay in $p$
for this term.

We next combine the Fourier coefficient products.
The classical Hecke relations \cite[\S\,2]{GL} imply
\begin{equation}
A_\pi(1,m)A_\pi(1,n)= 
 \sum_{r\mid(m,n)}A_\pi(r,mn/r^2).
\label{eq:holomorphic-pieri-two}
\end{equation}
This identity also holds at $p$, by the Schur polynomial identity
\eqref{eq:local-Schur} with $\gamma_p=0$.
For $(dbn,p)=1$, taking complex conjugates in
\cite[Eq.\,(1.6)]{BC} and using \eqref{eq:unramified-unitarity} gives
\begin{equation}
A_\pi(d,b)A_\pi(1,n)= 
 \sum_{\substack{n_1n_2n_3=n\\n_2\mid b,\ n_3\mid d}}
 A_\pi(dn_2/n_3,bn_1/n_2).
\label{eq:holomorphic-pieri-three}
\end{equation}
Expanding $M_\pi^2$ using \eqref{eq:intro-mollifier} introduces indices $\ell_1$, $\ell_2\leq L$ with weights $x_{\ell_1}x_{\ell_2}/\sqrt{\ell_1\ell_2}$.
For fixed $m$, $n$, applying \eqref{eq:holomorphic-pieri-two}, \eqref{eq:holomorphic-pieri-three} and their complex conjugates at unramified primes, we can express each product, including the finite
Euler factors, as a finite sum of single Fourier coefficients $A_\pi(u,v)$. After summing over $\ell_1$ and $\ell_2$, let $b_{m,n}(u,v)$ be the total coefficient of $A_\pi(u,v)$. Thus
\begin{equation}
M_\pi^2A_\pi(1,m)A_\pi(1,n)
= \sum_{u,v\geq1}b_{m,n}(u,v)A_\pi(u,v).
\label{eq:holomorphic-hecke-coefficients}
\end{equation}
Here $m$, $n$ are the Dirichlet indices in \eqref{eq:holomorphic-moment-afe}, and $u,v$ are the Fourier indices in this expansion. The coefficients $b_{m,n}(u,v)$ are real,
depend on $Q,L,m,n,u,v$, and are independent of $\pi$.
Since $A_\pi(1,1)=1$, setting $m=n=1$ and writing the Fourier
indices as $(d,b)$ gives
\begin{equation*}
M_\pi^2= \sum_{d,b\geq1}b_{1,1}(d,b)A_\pi(d,b).
\end{equation*}
We now bound the sum of the absolute coefficients in
\eqref{eq:holomorphic-hecke-coefficients} and the Fourier indices
that occur there. These bounds will be used to estimate the Weyl terms. For fixed $\ell_1$, $\ell_2$, $m$, $n$, the Hecke relations involve only divisor sums. The divisor bound shows that 
the number of resulting terms is $O_\eps((mn\ell_1\ell_2)^\eps)$ counted with multiplicity. The finite Euler factors contribute a further constant depending only on $Q$. Summing the absolute weights of these terms and 
using $|x_\ell|\leq1$ gives 
\begin{equation}
\sum_{u,v\geq1}|b_{m,n}(u,v)|
\ll_{Q,\eps}\sum_{\ell_1,\ell_2\leq L}
 \frac{|x_{\ell_1}x_{\ell_2}|}{\sqrt{\ell_1\ell_2}}
 (mn\ell_1\ell_2)^\eps
\ll_{Q,\eps}(mnL^2)^\eps L.
\label{eq:holomorphic-hecke-coefficient-mass}
\end{equation}

Next we bound $u^2v$, which controls the kernel arguments in the
Weyl terms below. In \eqref{eq:holomorphic-pieri-three}, multiplying
$A_\pi(d,b)$ by $A_\pi(1,n)$ produces Fourier indices
$(dn_2/n_3,bn_1/n_2)$. Since $n_1n_2n_3=n$, 
\begin{equation*}
\left(\frac{dn_2}{n_3}\right)^2\frac{bn_1}{n_2}
= \frac{d^2bn}{n_3^3}\leq d^2bn.
\end{equation*}
Iteration gives $u^2v\leq mn\ell_1\ell_2$ at unramified primes.
The same bound holds at $p$ by \eqref{eq:holomorphic-pieri-two},
since $(\ell_1\ell_2,p)=1$.
The finite Euler factors contribute only a constant depending on $Q$.
Thus, for $\ell_1,\ell_2\leq L$ and $mn\leq p^{1+\eps}L^2$,
every non-zero $b_{m,n}(u,v)$ satisfies
\begin{equation}
u^2v\ll_Q mnL^2\leq p^{1+\eps}L^4.
\label{eq:holomorphic-primal-indices}
\end{equation}
Since $A_\pi(1,1)=1$, the first sum has Kuznetsov indices $(1,1,u,v)$.
For the second sum in \eqref{eq:holomorphic-moment-afe}, reindex
by $\pi\mapsto\widetilde\pi$. The family and the full adjoint
value are unchanged. By \eqref{eq:contragredient-coeff},
$M_{\widetilde\pi}=\overline{M_\pi}$, while
\eqref{eq:global-root} and \eqref{eq:ramified-Hecke} give
$\epsilon(\widetilde\pi)^2=pA_\pi(p^2,1)$.
Thus the Fourier coefficient product becomes
\begin{equation}
p\,\overline{M_\pi^2}A_\pi(p^2,1)A_\pi(1,m)A_\pi(1,n),
\label{changed form}
\end{equation}
with test $K_{-\mu}^-(mnL^2/p)$, since
$\mu_{\widetilde\pi}=-\mu_\pi$ by \eqref{eq:functional-equation}
and $h$ is even.
Apply \eqref{eq:holomorphic-pieri-two} to the last two factors in \eqref{changed form}.
For each divisor $e\mid(m,n)$, write $m=em'$, $n=en'$ and
$v=m'n'$. The new coefficient is $A_\pi(e,v)$ and the weight
argument is $e^2vL^2/p$. Summing over the $\tau(v)$ ordered
factorizations $m'n'=v$, with $\sqrt{mn}=e\sqrt v$, gives the
scalar weight $\tau(v)/(e\sqrt v)$.
To expand $\overline{M_\pi^2}$, take the conjugate of the
$m=n=1$ case of \eqref{eq:holomorphic-hecke-coefficients}.
Since $b_{1,1}(d,b)$ is real and $(ev,p)=1$, multiplicativity gives
\begin{equation*}
p\,\overline{M_\pi^2}A_\pi(p^2,1)A_\pi(e,v)
= p\sum_{d,b\geq1}b_{1,1}(d,b)
 \overline{A_\pi(d,b)}A_\pi(p^2e,v).
\end{equation*}
The Fourier pair on the right has Kuznetsov indices $(d,b,p^2e,v)$.
All mollifier indices are prime to $p$. Together with
\eqref{eq:holomorphic-primal-indices} for $m=n=1$ and the dual
cutoff, this gives
\begin{equation}
d^2b\ll_Q L^2,\qquad
e^2v\leq p^{1+\eps}/L^2,\qquad (dbev,p)=1.
\label{eq:holomorphic-dual-indices}
\end{equation}
There is no dual diagonal, since the required equality $d=p^2e$
contradicts $(d,p)=1$.

To sum the short and long Weyl estimates below over all indices,
we bound the sums of the absolute values of their coefficients.
The first and dual sums have coefficients
$b_{m,n}(u,v)/\sqrt{mn}$ and
$p\tau(v)b_{1,1}(d,b)/(e\sqrt v)$, respectively.
Combining \eqref{eq:holomorphic-hecke-coefficient-mass} with the
divisor bound gives
\begin{equation}
\begin{split}
&~\sum_{\substack{m,n\geq1\\mn\leq p^{1+\eps}L^2}}
 \frac1{\sqrt{mn}}\sum_{u,v\geq1}|b_{m,n}(u,v)|
 \ll_{Q,\eps}p^{1/2+O(\eps)}L^{2+O(\eps)},\\
&~\sum_{\substack{e,v\geq1\\(ev,p)=1\\e^2v\leq p^{1+\eps}/L^2}}
 \frac{\tau(v)}{e\sqrt v}\sum_{d,b\geq1}|b_{1,1}(d,b)|
 \ll_{Q,\eps}p^{1/2+O(\eps)}.
\end{split}
\label{eq:holomorphic-outside-sums}
\end{equation}
The second bound excludes the factor $p$. Including it, the sum
of absolute coefficients in the dual contribution is
$O_{Q,\eps}(p^{3/2+O(\eps)})$.

Up to an error $O_{h,Q,\eta,\eps}(p^{-1})$,
applying the preceding coefficient identities to
\eqref{eq:holomorphic-moment-afe} with these cutoffs gives
\begin{equation}
\begin{split}
&\mathcal M_{1,1}(p)= c_{\rm cusp}
 \sum_{\substack{m,n,u,v\geq1\\mn\leq p^{1+\eps}L^2}}
 \sum_{\pi\in\mathcal B(p)}
 \frac{b_{m,n}(u,v)A_\pi(u,v)K_{\mu_\pi}^+(mn/(pL^2))}
 {\sqrt{mn}L(1,\pi,\mathrm{Ad})}\\*
&~\qquad+pc_{\rm cusp}
 \sum_{\substack{e,v,d,b\geq1\\(ev,p)=1\\e^2v\leq p^{1+\eps}/L^2}}
 \sum_{\pi\in\mathcal B(p)}
 \frac{\tau(v)b_{1,1}(d,b)\overline{A_\pi(d,b)}A_\pi(p^2e,v)
K_{-\mu_\pi}^-(e^2vL^2/p)}
 {e\sqrt v L(1,\pi,\mathrm{Ad})}.
\end{split}
\label{eq:holomorphic-cuspidal-expansion}
\end{equation}

The tests in \eqref{eq:holomorphic-cuspidal-expansion} are
admissible by Lemma~\ref{lem:holomorphic-square-weights} and
\S\,\ref{subsec:kuz-afe-compatibility}.
Since the truncated Dirichlet sums and their Hecke expansions
are finite, Theorem~\ref{thm:kuznetsov} applies term by term.
We retain the length cutoffs in the Weyl terms. In the diagonal
and the four remaining spectral contributions, we remove the
length restrictions, retaining $(ev,p)=1$ in the dual sums.
The additional terms have weight arguments
$mn/(pL^2)>p^\eps$ or $e^2vL^2/p>p^\eps$.

The truncation argument following \eqref{eq:first-moment-afe} applies with \eqref{eq:holomorphic-square-weight-decay} and \eqref{eq:holomorphic-hecke-coefficient-mass}. For the maximal new spectrum, we first separate the powers of $p$ in the local coefficients and use \eqref{eq:gl2-induced-coefficient-average}. The resulting polynomial losses are absorbed by choosing the weight decay sufficiently large. Thus each added contribution is $O_{h,Q,\eta,J,\eps,A}(p^{-A})$ for every fixed $A>0$, with the cutoff exponent $\eps$ and the desired $A$ fixed before choosing the decay exponent.

Taking $A=1$, define the diagonal contribution by
\begin{equation}
\mathcal D_{1,1}:= P\sum_{m,n\geq1}
 \frac{b_{m,n}(1,1)}{\sqrt{mn}}
 \int_{\Re(\mu)=0}K_\mu^+(mn/(pL^2))
 \operatorname{spec}(\mu)\,\dd\mu.
\label{eq:holomorphic-diagonal-def}
\end{equation}
For $i=4,5,6$, the Weyl contributions are
\begin{equation*}
\begin{split}
\mathcal S_{i,1,1}:=
&~\sum_{\substack{m,n\geq1\\mn\leq p^{1+\eps}L^2}}\frac1{\sqrt{mn}}
 \sum_{u,v\geq1}b_{m,n}(u,v) S_i\!\left(1,1,u,v\,;
 \mu\mapsto K_\mu^+(mn/(pL^2))\right)\\*
&~+p\sum_{\substack{e,v\geq1\\(ev,p)=1\\e^2v\leq p^{1+\eps}/L^2}}
 \frac{\tau(v)}{e\sqrt v}
 \sum_{d,b\geq1}b_{1,1}(d,b) S_i\!\left(d,b,p^2e,v\,;
 \mu\mapsto K_{-\mu}^-(e^2vL^2/p)\right).
\end{split}
\end{equation*}
Here $S_i$ uses the displayed Fourier indices and
spectral test, with the conventions of
\S\,\ref{subsec:kuznetsov-identity} and the transforms in
\eqref{eq:kuz-transforms}.
The four remaining spectral contributions are
\begin{equation*}
\begin{split}
&\mathcal R_{{\rm old},1,1}:=
c_{\rm cusp}\sum_{\Pi\in\mathcal B(1)}\sum_{j=0}^2
 \frac1{L(1,\Pi,\mathrm{Ad})}\\*
&~\quad\times\Biggl(
 \sum_{m,n\geq1}\frac1{\sqrt{mn}}
 \sum_{u,v\geq1}b_{m,n}(u,v)
 \overline{B_\Pi^{(j)}(1,1)}B_\Pi^{(j)}(u,v)K_{\mu_\Pi}^+(mn/(pL^2))\\
&~\qquad+p\sum_{\substack{e,v\geq1\\(ev,p)=1}}
 \frac{\tau(v)}{e\sqrt v}
 \sum_{d,b\geq1}b_{1,1}(d,b)
 \overline{B_\Pi^{(j)}(d,b)}B_\Pi^{(j)}(p^2e,v)K_{-\mu_\Pi}^-(e^2vL^2/p)\Biggr),
\end{split}
\end{equation*}
\begin{equation*}
\begin{split}
&\mathcal R_{\min,1,1}:=
c_{\min}\sum_{j=0}^2
 \int_{\Re(s_1)=0}\int_{\Re(s_2)=0}
 \frac1{\prod_{1\leq i<k\leq3}|\zeta(1+s_i-s_k)|^2}\\*
&~\quad\times\Biggl(
 \sum_{m,n\geq1}\frac1{\sqrt{mn}}
 \sum_{u,v\geq1}b_{m,n}(u,v)
 \overline{B_{\mathbf s}^{(j)}(1,1)}B_{\mathbf s}^{(j)}(u,v)K_{\mathbf s}^+(mn/(pL^2))\\
&~\qquad+p\sum_{\substack{e,v\geq1\\(ev,p)=1}}
 \frac{\tau(v)}{e\sqrt v}
 \sum_{d,b\geq1}b_{1,1}(d,b)
 \overline{B_{\mathbf s}^{(j)}(d,b)}B_{\mathbf s}^{(j)}(p^2e,v)K_{-\mathbf s}^-(e^2vL^2/p)\Biggr)
 \,\dd s_1\,\dd s_2,
\end{split}
\end{equation*}
and, with $\mu=(z+\tau_\phi,z-\tau_\phi,-2z)$,
\begin{equation*}
\begin{split}
&\mathcal R_{\max,{\rm new},1,1}:=
\chi_p c_{\max}\sum_{\phi\in\mathcal B_2(p)}
 \int_{\Re(z)=0}
 \frac1{L(1,\phi,\mathrm{Ad})|L(1+3z,\phi)|^2}\\*
&~\quad\times\Biggl(
 \sum_{m,n\geq1}\frac1{\sqrt{mn}}
 \sum_{u,v\geq1}b_{m,n}(u,v)A_{\phi,z}(u,v)
K_\mu^+(mn/(pL^2))\\
&~\qquad+p\sum_{\substack{e,v\geq1\\(ev,p)=1}}
 \frac{\tau(v)}{e\sqrt v}
 \sum_{d,b\geq1}b_{1,1}(d,b)
 \overline{A_{\phi,z}(d,b)}A_{\phi,z}(p^2e,v) K_{-\mu}^-(e^2vL^2/p)\Biggr)\,\dd z,
\end{split}
\end{equation*}
\begin{equation*}
\begin{split}
&\mathcal R_{\max,{\rm old},1,1}:=
c_{\max}\sum_{\phi\in\mathcal B_2(1)}\sum_{j=0}^2
 \int_{\Re(z)=0}
 \frac1{L(1,\phi,\mathrm{Ad})|L(1+3z,\phi)|^2}\\*
&~\quad\times\Biggl(
 \sum_{m,n\geq1}\frac1{\sqrt{mn}}
 \sum_{u,v\geq1}b_{m,n}(u,v)
 \overline{B_{\phi,z}^{(j)}(1,1)}B_{\phi,z}^{(j)}(u,v)K_\mu^+(mn/(pL^2))\\
&~\qquad+p\sum_{\substack{e,v\geq1\\(ev,p)=1}}
 \frac{\tau(v)}{e\sqrt v}
 \sum_{d,b\geq1}b_{1,1}(d,b)
 \overline{B_{\phi,z}^{(j)}(d,b)}B_{\phi,z}^{(j)}(p^2e,v)K_{-\mu}^-(e^2vL^2/p)\Biggr)\,\dd z.
\end{split}
\end{equation*}
The continuous measures retain the conventions of
\S\,\ref{subsec:kuznetsov-identity}. Thus
\begin{equation}
\begin{split}
\chi_p\mathcal M_{1,1}(p)=
&~\mathcal D_{1,1}+\sum_{i=4}^6\mathcal S_{i,1,1}
 -\mathcal R_{{\rm old},1,1}-\mathcal R_{\min,1,1}\\
&~-\mathcal R_{\max,{\rm new},1,1}
 -\mathcal R_{\max,{\rm old},1,1}
 +O_{h,Q,\eta,\eps}(p^{-1}).
\end{split}
\label{eq:holomorphic-moment-decomposition}
\end{equation}
We estimate these contributions in the following subsections.

\subsubsection{The diagonal term}

By \eqref{eq:holomorphic-diagonal-def}, the first trace diagonal
selects the coefficient of $A_\pi(1,1)$ in
\eqref{eq:holomorphic-hecke-coefficients}.
We first compute the sum of these coefficients with
$\mathcal V_\mu^+$ replaced by $1$, and then estimate the error
in this replacement.

For $q\nmid pQ$, write
\begin{equation*}
L_q=L_q(1/2,\pi)=\sum_{j\geq0}\frac{A_\pi(1,q^j)}{q^{j/2}}.
\end{equation*}
By \eqref{eq:holomorphic-pieri-two},
\begin{equation*}
L_q^2=\sum_{j,k\geq0}q^{-(j+k)/2}
 \sum_{r=0}^{\min(j,k)}A_\pi(q^r,q^{j+k-2r}).
\end{equation*}
The coefficients of $A_\pi(1,1)$, $A_\pi(q,1)$,
$A_\pi(1,q)$ and $A_\pi(q^2,1)$ in this expansion are
$1$, $q^{-1}$, $2q^{-1/2}$ and $q^{-2}$, respectively.
By \eqref{eq:holomorphic-pieri-three} with $n=q$,
multiplication by $A_\pi(1,q)$ produces $A_\pi(1,1)$ only
from $A_\pi(q,1)$, and produces $A_\pi(q,1)$ only from
$A_\pi(1,q)$ and $A_\pi(q^2,1)$.
Thus the coefficients of $A_\pi(1,1)$ in
$L_q^2$, $L_q^2A_\pi(1,q)$ and $L_q^2A_\pi(1,q)^2$ are
\begin{equation*}
1,\qquad q^{-1},\qquad 2q^{-1/2}+q^{-2},
\end{equation*}
respectively. The second value comes from $m=n=q$.
For the third, the contributing Dirichlet pairs are
$(m,n)=(q,1),(1,q),(q^2,q^2)$, with weights
$q^{-1/2},q^{-1/2},q^{-2}$.

We first compute the local diagonal with the logarithmic mollifier cutoffs removed, keeping $(\ell,pQ)=1$. The absolutely convergent majorant below justifies this formal computation.
The local mollifier factor at $q\nmid pQ$ is then
$1-\frac q{q+1}A_\pi(1,q)/\sqrt q$.
Expanding its square and using the three coefficients above,
the coefficient of $A_\pi(1,1)$ in the local product is
\begin{equation*}
1-\frac{2q}{q+1}q^{-3/2}
 +\left(\frac q{q+1}\right)^2(2q^{-3/2}+q^{-3})
=1-\frac{2}{\sqrt q(q+1)^2}+\frac{1}{q(q+1)^2}.
\end{equation*}
For $q\mid Q$, the finite inverse Euler factors in
\eqref{eq:intro-mollifier} cancel $L_q^2$ exactly.
At $p$, the mollifier has local factor $1$.
By \eqref{eq:holomorphic-pieri-two},
$A_\pi(1,p^j)A_\pi(1,p^k)$ contributes to $A_\pi(1,1)$
only when $j=k=0$, so the local diagonal factor at $p$ is also $1$.

We now justify removing the mollifier cutoffs.
Taking absolute values before collecting the Hecke coefficients
gives local diagonal factors $1+O(q^{-3/2})$ for $q\nmid pQ$.
Since $\sum_q(\log q)q^{-3/2}<\infty$, the corresponding
absolute series weighted by $\log\ell_1+\log\ell_2$ also
converges. The elementary bound
\begin{equation*}
\left|1-\mathbf1_{\ell\leq L}
 \frac{\log(L/\ell)}{\log L}\right|
\leq\frac{\log\ell}{\log L},\qquad \ell\geq1,
\end{equation*}
applied to both mollifier indices therefore gives
\begin{equation*}
\sum_{m,n\geq1}\frac{b_{m,n}(1,1)}{\sqrt{mn}}
=\prod_{q\nmid pQ}
 \left(1-\frac{2}{\sqrt q(q+1)^2}
             +\frac{1}{q(q+1)^2}\right)
 +O_Q\!\left(\frac1{\log L}\right).
\end{equation*}

Next fix $0<c<1/4$ to estimate the replacement of
$\mathcal V_\mu^+$ by $1$.
Weighting the Dirichlet pairs by $(mn)^c$ gives local absolute
factors
\begin{equation*}
1+O\big(q^{-3/2+2c}+q^{-3/2+c}+q^{-3+4c}\big).
\end{equation*}
The three powers of $q^c$ correspond to
$mn=q^2,q,q^4$, respectively.
These factors have an absolutely convergent product for
$0<c<1/4$. Since the cutoff factors have absolute value at most
$1$, this proves, uniformly for $L<p$ and $p\nmid Q$,
\begin{equation*}
\sum_{m,n\geq1}\frac{|b_{m,n}(1,1)|}{\sqrt{mn}}(mn)^c
\ll_{Q,c}1.
\end{equation*}
\eqref{eq:holomorphic-square-weight-small} and
\eqref{eq:holomorphic-square-weight-decay} imply, for all $y>0$,
\begin{equation*}
|\mathcal V_\mu^+(y)-1|
\ll_{J,c}y^c(1+\|\mu\|)^C.
\end{equation*}
Combining these bounds with \eqref{eq:integrable-spectral-density}
in \eqref{eq:holomorphic-diagonal-def} yields
\begin{equation*}
\mathcal D_{1,1}
=P\mathcal H_h\sum_{m,n\geq1}\frac{b_{m,n}(1,1)}{\sqrt{mn}}
 +O_{h,Q,J,c}\big(p^2(pL^2)^{-c}\big).
\end{equation*}

Finally, the factor at $q=p$ in the unramified Euler product is
\begin{equation*}
1-\frac{2}{\sqrt p(p+1)^2}+\frac{1}{p(p+1)^2}=1+O(p^{-5/2}).
\end{equation*}
Including it changes the main term by $O_{h,Q}(p^{-1/2})$.
Since $L=p^\eta$, we have $1/\log L=1/(\eta\log p)$ and
$(pL^2)^{-c}=p^{-c(1+2\eta)}\ll_{\eta,c}1/\log p$.
The preceding estimates therefore give
\begin{equation}
\mathcal D_{1,1}=P\mathcal H_h\prod_{q\nmid Q}
 \left(1-\frac{2}{\sqrt q(q+1)^2}
             +\frac{1}{q(q+1)^2}\right)
 +O_{h,Q,\eta,J}\left(\frac{p^2}{\log p}\right).
\label{eq:holomorphic-square-diagonal}
\end{equation}

\subsubsection{The short Weyl elements}

\begin{lemma} 
\label{lem:holomorphic-moment-short}
For $i=4,5$, the first and dual contributions to
$\mathcal S_{i,1,1}$ are each $O_{h,Q,\eta,J}(p^{-1})$.
\end{lemma}

\begin{proof}
Choose $A$ sufficiently large in terms of $\eta$, and take
$a_0,J$ large enough for the corresponding kernel estimates.
Since $L=p^\eta<p^{1/2}$, the weight arguments
$mn/(pL^2)$ and $e^2vL^2/p$ are at least $p^{-2}$.
By Remark~\ref{rem:kuz-kernel-uniformity}, \eqref{eq:holomorphic-square-test-seminorm} with $c=\eps/2$, and \eqref{eq:kuz-short-bound}, for either test functions in
\eqref{eq:holomorphic-cuspidal-expansion},
\begin{equation*}
|\Phi_{w_4}(y)|+|\Phi_{w_5}(y)|
\ll_{A,h,J,\eps}p^\eps|y|^A.
\end{equation*}

For the first indices $(1,1,u,v)$, the $w_4$ modulus conditions
are parametrized by $D_2=pk/(p,u)$ and $D_1=uD_2^2$.
For $w_5$, they give $D_1=pk$ and $D_2=vD_1^2$, with $k\geq1$.
Consequently,
\begin{equation*}
|y_{w_4}|=\frac{v(p,u)^3}{up^3k^3}
\leq\frac{u^2v}{p^3k^3},\qquad
|y_{w_5}|=\frac{u}{vp^3k^3}
\leq\frac{u^2v}{p^3k^3}.
\end{equation*}
By \eqref{eq:holomorphic-primal-indices}, both arguments are
$\ll_Q p^{-2+\eps}L^4k^{-3}$.
For the dual indices $(d,b,p^2e,v)$, using $(dbev,p)=1$,
the corresponding parametrizations are
\begin{equation*}
\begin{aligned}
w_4:\quad &D_2=\frac{bk}{(b,e)},\qquad
D_1=\frac{p^2bek^2}{(b,e)^2},\\
w_5:\quad &D_1=\frac{pdk}{(d,v)},\qquad
D_2=\frac{p^2dvk^2}{(d,v)^2},\qquad k\geq1.
\end{aligned}
\end{equation*}
Thus \eqref{eq:holomorphic-dual-indices} gives
\begin{equation*}
\begin{aligned}
|y_{w_4}|=\frac{dv(b,e)^3}{p^2bek^3}
\ll_Q p^{-1+\eps}k^{-3},\qquad\text{and}\qquad
|y_{w_5}|=\frac{eb(d,v)^3}{pdvk^3}
\ll_Q Lp^{-1/2+\eps}k^{-3}.
\end{aligned}
\end{equation*}

Apply $|\widetilde S|\leq D_1D_2$ and the kernel bound above.
In each of the four sums, the sum over $k$ is bounded by
$\sum_{k\geq1}k^{-3A}\ll_A1$.
Using $P\asymp p^2$ and \eqref{eq:holomorphic-outside-sums},
each first contribution is
\begin{equation*}
\ll_{A,h,Q,J,\eps}
p^{5/2+O(\eps)}L^{2+O(\eps)}
\left(\frac{p^\eps L^4}{p^2}\right)^A,
\end{equation*}
and each dual contribution, including its outside factor $p$, is
\begin{equation*}
\ll_{A,h,Q,J,\eps}
p^{7/2+O(\eps)}
\left(p^{-1+\eps}+Lp^{-1/2+\eps}\right)^A.
\end{equation*}
Since $L=p^\eta$ and $\eta<1/2$, both bounds are
$O_{h,Q,\eta,J}(p^{-1})$ for $A$ sufficiently large in terms
of $\eta$, as $\eps$ is arbitrarily small. This proves the lemma.
\end{proof}

\subsubsection{The long Weyl element}

\begin{lemma} 
\label{lem:holomorphic-moment-long}
The first and dual contributions to $\mathcal S_{6,1,1}$ are each
$O_{h,Q,\eta,J}(p^{-1})$.
\end{lemma}

\begin{proof}
Choose $A>1$ and the auxiliary parameters as in
Lemma~\ref{lem:holomorphic-moment-short}.
For the indices $(1,1,u,v)$ from the first sum, we have
\begin{equation*}
|y_1^2y_2|= \frac{uv^2}{D_1^3},\qquad
|y_1y_2^2|=\frac{u^2v}{D_2^3}.
\end{equation*}
By Remark~\ref{rem:kuz-kernel-uniformity} and \eqref{eq:holomorphic-square-test-seminorm} with $c=1$,
\eqref{eq:w6-small} with decay exponent $A+2$ and growth
exponent one, and the bounds $|S^{(p)}|\leq(D_1D_2)^2$ and
$1+|y_1|+|y_2|\leq3u^2vD_1D_2$, the sum over $D_1,D_2$
for fixed $m,n,u,v$ is bounded by
\begin{equation*}
\begin{aligned}
\ll_{A,h,J}P\frac{pL^2}{mn}u^4v^5(u^2v)^A
\sum_{\substack{D_1,D_2\geq1\\p\mid D_1,D_2}}
D_1^{-4}D_2^{2-3A}\ll_{A,h,J}\frac{pL^2}{mn}(u^2v)^5
\left(\frac{u^2v}{p^3}\right)^A.
\end{aligned}
\end{equation*}
By \eqref{eq:holomorphic-primal-indices} and
\eqref{eq:holomorphic-outside-sums}, the total first contribution is
\begin{equation*}
\ll_{A,h,Q,J,\eps}
p^C\left(\frac{p^\eps L^4}{p^2}\right)^A.
\end{equation*}

For the dual indices $(d,b,p^2e,v)$, twisted multiplicativity
reduces the case $\operatorname{ord}_p(D_2)=1$ to local moduli $(p^\alpha,p)$,
up to unit twists. The case $\alpha=1$ vanishes by the local
identities in \S\,\ref{subsec:kuz-kloosterman}, since $p\nmid d$.
For $\alpha\geq2$, the same substitution $B_1=p\beta$ in
\eqref{eq:kuz-long-def} as in that subsection gives
\begin{equation*}
S^{(p)}(\epsilon_2v,\epsilon_1p^2e,d,b,p^\alpha,p)
= r_p(\epsilon_1p^2e)S(d,\epsilon_2pv,p^\alpha)=0.
\end{equation*}
The last equality follows by summing over
$x\mapsto x+p^{\alpha-1}t$, $t\bmod p$, using $p\nmid d$.
Thus only $p^2\mid D_2$ contributes.
The dual arguments and \eqref{eq:holomorphic-dual-indices} give
\begin{equation*}
\begin{aligned}
|y_1^2y_2|=\frac{p^2d^2bev^2}{D_1^3},
\end{aligned}\end{equation*}
and 
\begin{equation*}
\begin{aligned}
|y_1y_2^2|&=\frac{p^4db^2e^2v}{D_2^3}
\ll_Q L^2p^{-1+\eps}\left(\frac{p^2}{D_2}\right)^3.
\end{aligned}
\end{equation*}
The same summation over $D_1,D_2$, now with $p^2\mid D_2$ and test
factor $p/(e^2vL^2)$ from
\eqref{eq:holomorphic-square-test-seminorm}, gives
\begin{equation*}
\ll_{A,h,J}\frac{p^9}{L^2}d^5b^3ev^4
\left(\frac{db^2e^2v}{p^2}\right)^A.
\end{equation*}
By \eqref{eq:holomorphic-dual-indices} and
\eqref{eq:holomorphic-outside-sums}, including the outside
factor $p$, the total dual contribution is
\begin{equation*}
\ll_{A,h,Q,J,\eps}
p^C\big(L^2p^{-1+\eps}\big)^A.
\end{equation*}
Here and above, $C$ is independent of $A$.
Since $L=p^\eta$ and $\eta<1/2$, taking $A$ sufficiently large
in terms of $\eta$ makes both contributions
$O_{h,Q,\eta,J}(p^{-1})$, as $\eps$ is arbitrarily small.
\end{proof}

\subsubsection{The old and minimal Eisenstein spectra}

\begin{proposition}
\label{prop:holomorphic-level-one-remainders}
For every $\eps>0$, we have
\begin{equation*}
\begin{aligned}
\mathcal R_{{\rm old},1,1}
&\ll_{h,Q,\eta,J,\eps}p^\eps L^{12/7+\eps},\\
\mathcal R_{\min,1,1}
&\ll_{h,Q,\eta,J,\eps}p^\eps L^{1+\eps},\\
\mathcal R_{\max,{\rm old},1,1}
&\ll_{h,Q,\eta,J,\eps}p^\eps L^{39/32+\eps}.
\end{aligned}
\end{equation*}
\end{proposition}

\begin{proof}
Let $\Pi$ lie in the old cuspidal, minimal or maximal old
spectrum, with $\theta=5/14,0,7/64$, respectively.
Define $M_\Pi$ by \eqref{eq:intro-mollifier} using the coefficients
of $\Pi$. The coefficient bounds of \S\,\ref{prelim} give
\begin{equation}
|M_\Pi|^2\ll_{Q,\eps}
\left(\sum_{\ell\leq L}\ell^{-1/2+\theta+\eps}\right)^2
\ll_{Q,\eps}L^{1+2\theta+2\eps}.
\label{eq:holomorphic-level-one-mollifier-bound}
\end{equation}

We compute the local series at $p$ to recover the full squared $L$-function after summing over the oldvectors. Insert \eqref{eq:holomorphic-square-weight} and put $s=1/2+u$ on the initial Mellin line $\Re(u)>1$.
Writing $m=p^km^\prime$ and $n=p^\ell n^\prime$ with $(m^\prime n^\prime,p)=1$, \eqref{eq:holomorphic-pieri-two} gives the local series
\begin{equation}
T_j(s)\colon \!\!=\sum_{k,\ell\geq0}p^{-(k+\ell)s}
 \sum_{r=0}^{\min(k,\ell)}
 A_\Pi^{(j)}(p^r,p^{k+\ell-2r}),\qquad j=0,1,2.
\label{eq:holomorphic-old-local-series}
\end{equation}
By \eqref{eq:kuz-old-coeff} and \eqref{eq:holomorphic-pieri-two}, \eqref{eq:holomorphic-old-local-series} gives
\begin{equation}
T_0(s)=L_p(s,\Pi)^2,\quad
T_1(s)=2p^{1/2-s}L_p(s,\Pi)^2,\quad
T_2(s)=p^{1-2s}L_p(s,\Pi)^2.
\label{eq:holomorphic-old-local-series1}
\end{equation}

We now sum the first Dirichlet series over the three oldvectors. By \eqref{eq:holomorphic-hecke-coefficients}, the primes other than $p$ contribute $M_\Pi^2L^{(p)}(s,\Pi)^2$.
The oldspace summation in the proof of Lemma \ref{lem:oldspace-cancellation} and \eqref{eq:holomorphic-old-local-series1} therefore give
\begin{equation}
\begin{aligned}
&\sum_{m,n\geq1}\dfrac1{(mn)^s}
 \sum_{d,b\geq1}b_{m,n}(d,b)
 \sum_{j=0}^2\overline{B_\Pi^{(j)}(1,1)}B_\Pi^{(j)}(d,b)\\
&\quad=M_\Pi^2L^{(p)}(s,\Pi)^2
 \sum_{j=0}^2(Q_\Pi^{-1})_{0j}T_j(s)\\
&\quad=M_\Pi^2L(s,\Pi)^2
 \Bigl((Q_\Pi^{-1})_{00}
 +2p^{1/2-s}(Q_\Pi^{-1})_{01}
 +p^{1-2s}(Q_\Pi^{-1})_{02}\Bigr).
\end{aligned}
\label{eq:holomorphic-level-one-first-series}
\end{equation}
Taking $s=1/2+u$ in \eqref{eq:holomorphic-level-one-first-series}, the factor in parentheses is $O(1)$ on $\Re(u)=\eps$ by Lemma \ref{lem:oldspace-cancellation}.
For fixed $\Pi$, the Mellin integral multiplying $h(\mu_\Pi)M_\Pi^2$ in the first contribution is
\begin{equation}
\begin{aligned}
&\dfrac1{2\pi\i}\int_{(c)}(pL^2)^uG_J(u,\mu_\Pi)^2
\left(\dfrac{\gamma_{\mu_\Pi}(1/2+u)}{\gamma_{\mu_\Pi}(1/2)}\right)^2L(1/2+u,\Pi)^2\\
&\qquad\times\Bigl(
(Q_\Pi^{-1})_{00}
+2p^{-u}(Q_\Pi^{-1})_{01}
+p^{-2u}(Q_\Pi^{-1})_{02}\Bigr)\frac{\dd u}{u}.
\end{aligned}
\label{eq:holomorphic-level-one-first-mellin}
\end{equation}
For each of the three level one spectra, move the $u$-contour in \eqref{eq:holomorphic-level-one-first-mellin} from $\Re(u)=c>1$ to $\Re(u)=\eps$ without a residue.
By \eqref{eq:afe-polynomial}, $P_J(u,\mu)^2$ cancels the Eisenstein poles with multiplicity.
The level one convexity bound \cite[Theorem~5.41]{IK}, \eqref{eq:zeta-pole-removed}, \eqref{eq:gamma-ratio-bound}, \eqref{eq:level-one-normalizer-bounds} and \eqref{eq:level-one-counting-bound}, together with the rapid decrease of $h$, bound the remaining spectral and Mellin integrals independently of $p$.
Combining \eqref{eq:holomorphic-level-one-mollifier-bound} with $|(pL^2)^u|=(pL^2)^\eps$ and rescaling $\eps$, the first contribution is
\begin{equation}\label{FB}
\ll_{h,Q,\eta,J,\eps}p^\eps L^{1+2\theta+\eps}.
\end{equation}

For the dual contribution, the Fourier pairs are $(d,b)$ and
$(p^2e,v)$, with $(dbev,p)=1$, by
\eqref{eq:holomorphic-dual-indices}.
We first combine the outside factor $p$ with the sum over
the three oldvectors.
Put $\lambda=A_\Pi(1,p)$ and $e_0=(1,0,0)^{\mathsf T}$.
The unramified Hecke relations give
$A_\Pi(p^2,1)=\overline\lambda^2-\lambda$ and
$A_\Pi(p,p)=|\lambda|^2-1$.
Hence \eqref{eq:kuz-old-coeff} gives the coefficient column
at $(p^2,1)$ as
\begin{equation*}
\begin{pmatrix}
A_\Pi^{(0)}(p^2,1)\\
A_\Pi^{(1)}(p^2,1)\\
A_\Pi^{(2)}(p^2,1)
\end{pmatrix}
=
\begin{pmatrix}
\overline\lambda^2-\lambda\\
\sqrt p(|\lambda|^2-1)\\
p\overline\lambda
\end{pmatrix}.
\end{equation*}
Since $(dbev,p)=1$, the coefficient columns at $(d,b)$ and
$(p^2e,v)$ are $A_\Pi(d,b)e_0$ and $A_\Pi(e,v)$ times
the displayed column, respectively.
The oldspace summation in the proof of
Lemma~\ref{lem:oldspace-cancellation} therefore gives
\begin{equation*}
\begin{aligned}
p\sum_{j=0}^2
 \overline{B_\Pi^{(j)}(d,b)}B_\Pi^{(j)}(p^2e,v)=\overline{A_\Pi(d,b)}A_\Pi(e,v)\,
 p\sum_{j=0}^2(Q_\Pi^{-1})_{0j}A_\Pi^{(j)}(p^2,1).
\end{aligned}
\end{equation*}
Using \eqref{eq:kuz-gram} and the identities
$a_p^2=b_p+b_p^2$,
$\sqrt p\,a_p=(p+1)b_p$ and $p(p+1)b_p=p-b_p$,
direct multiplication yields
\begin{equation}
p\sum_{j=0}^2(Q_\Pi^{-1})_{0j}A_\Pi^{(j)}(p^2,1)
=\frac{a_p^2\overline\lambda^2-b_p\lambda(1+b_p|\lambda|^2)}
       {\det(Q_\Pi)}
=(Q_\Pi^{-1})_{20}.
\label{eq:holomorphic-old-square-cancellation}
\end{equation}
The last equality is the cofactor formula.
Since $a_p^2,b_p\ll p^{-1}$, $|\lambda|\ll p^\theta$ and
$\det(Q_\Pi)\asymp1$ by
Lemma~\ref{lem:oldspace-cancellation}, the combined local
factor, including the outside $p$, satisfies
\begin{equation*}
|(Q_\Pi^{-1})_{20}|
\ll p^{2\theta-1}+p^{\theta-1}+p^{3\theta-2}
\ll p^{2\theta-1}.
\end{equation*}
The remaining dual indices are coprime to $p$. Thus the dual
Dirichlet series is
$\overline{M_\Pi^2}L^{(p)}(1/2+u,\Pi)^2$
multiplied by $(Q_\Pi^{-1})_{20}$.
The removed Euler factors are uniformly bounded on
$\Re(u)=\eps$, since $\theta<1/2$.
Using $\mathcal V_{-\mu_\Pi}^-$ in
\eqref{eq:holomorphic-square-weight}, the dual $u$-integral is
\begin{equation}
\frac1{2\pi\i}\int_{(c)}
\left(\frac p{L^2}\right)^uG_J(u,\mu_\Pi)^2
\left(\frac{\gamma_{\mu_\Pi}(1/2+u)}
{\gamma_{-\mu_\Pi}(1/2)}\right)^2
L^{(p)}(1/2+u,\Pi)^2\frac{\dd u}{u}.
\label{eq:holomorphic-level-one-dual-mellin}
\end{equation}
Move the $u$-contour in
\eqref{eq:holomorphic-level-one-dual-mellin} from $\Re(u)=c>1$
to $\Re(u)=\eps$ in the same way.
Combining this with \eqref{eq:holomorphic-level-one-mollifier-bound}
and \eqref{eq:holomorphic-old-square-cancellation} gives
\begin{equation*}
\ll_{h,Q,\eta,J,\eps}
p^{2\theta-1+\eps}L^{1+2\theta+\eps}.
\end{equation*}
Since $2\theta-1<0$, the first contribution \eqref{FB} dominates. Taking $\theta=5/14,0,7/64$ proves the three estimates.
\end{proof}

\subsubsection{The maximal Eisenstein new spectrum}

\begin{proposition}
\label{prop:holomorphic-maximal-new}
For every $\eps>0$, the first and dual contributions to
$\mathcal R_{\max,{\rm new},1,1}$ are each
\begin{equation*}
O_{h,Q,\eta,J,\eps}\big(p^{3/2+\eps}L^{1+\eps}\big).
\end{equation*}
\end{proposition}

\begin{proof}
Let $\Pi$ be the maximal Eisenstein representation associated with
$\phi\in\mathcal B_2(p)$ and $z\in\i\mathbb R$ in
\S\,\ref{subsec:kuz-maximal}, with
$\mu=(z+\tau_\phi,z-\tau_\phi,-2z)$.
Define $M_\Pi$ by \eqref{eq:intro-mollifier} using the coefficients
$A_{\phi,z}$. On the initial line $\Re(u)=c>1$,
\eqref{eq:holomorphic-hecke-coefficients} and
\S\,\ref{subsec:kuz-maximal} give
\begin{equation}
\begin{aligned}
&\sum_{m,n\geq1}\frac1{(mn)^{1/2+u}}
 \sum_{d,b\geq1}b_{m,n}(d,b)A_{\phi,z}(d,b)\\
&\qquad=M_\Pi^2L(1/2+u-z,\phi)^2
 \zeta(1/2+u+2z)^2.
\end{aligned}
\label{eq:holomorphic-maximal-new-first-series}
\end{equation}
Insert \eqref{eq:holomorphic-maximal-new-first-series} into the
Mellin representation \eqref{eq:holomorphic-square-weight}.
By \eqref{eq:afe-polynomial}, the factor
$\big((1/2+\mu_3)^2-u^2\big)^2$ in $P_J(u,\mu)^2$ cancels
the double zeta pole at $u=1/2-2z$.
The $\rm GL_2$ $L$-function is entire, and the numerator gamma factors
have no poles for $\Re(u)>0$, since
$|\Re(\mu_j)|\leq7/64$.
\eqref{eq:gamma-ratio-bound} and
\eqref{eq:zeta-pole-removed}, together with $e^{2u^2}$,
make the horizontal integrals tend to zero.
We may therefore move the $u$-contour to $\Re(u)=\eps$
without a residue.

On this line, the $\rm GL_2$ convexity bound
\cite[Theorem~5.41]{IK} gives
\begin{equation*}
|L(1/2+u-z,\phi)|^2
\ll_\eps p^{1/2+\eps}
 (1+|z|+|\tau_\phi|+|\Im(u)|)^C.
\end{equation*}
We next estimate the spectral average of $|M_\Pi|^2$.
The finite Euler product in \eqref{eq:intro-mollifier}
is $O_Q(1)$ by the coefficient bounds of \S\,\ref{prelim}
and the induced parameters of \S\,\ref{subsec:kuz-maximal}.
Thus
\begin{equation*}
|M_\Pi|^2\ll_Q
\sum_{\ell_1,\ell_2\leq L}
\frac{|x_{\ell_1}x_{\ell_2}|}{\sqrt{\ell_1\ell_2}}
|A_{\phi,z}(1,\ell_1)A_{\phi,z}(1,\ell_2)|.
\end{equation*}
Apply \eqref{eq:gl2-induced-coefficient-average} with
$(a,b,m,n)=(1,\ell_1,1,\ell_2)$.
Since $\ell_i\leq L<p$, all indices are coprime to $p$ and
$(1+\ell_1/p)(1+\ell_2/p)\ll1$.
For every fixed $C\geq0$, using $|x_\ell|\leq1$ gives
\begin{equation}
\begin{aligned}
&\sum_{\phi\in\mathcal B_2(p)}
 \frac1{L(1,\phi,\mathrm{Ad})}
 \int_{\Re(z)=0}
 \frac{|h(\mu)||M_\Pi|^2(1+\|\mu\|)^C}
 {|L(1+3z,\phi)|^2}\,|\dd z|\\
&\qquad\ll_{h,Q,C,\eps}p^{1+\eps}
 \left(\sum_{\ell\leq L}|x_\ell|\ell^{-1/2+\eps}\right)^2
 \ll_{h,Q,C,\eps}p^{1+\eps}L^{1+2\eps}.
\end{aligned}
\label{eq:holomorphic-maximal-new-mollifier-average}
\end{equation}
The estimates \eqref{eq:gamma-ratio-bound} and
\eqref{eq:zeta-pole-removed} and the factor $e^{2u^2}$
control the remaining integration in $\Im(u)$ by a fixed
polynomial in $1+\|\mu\|$.
Since $|(pL^2)^u|=(pL^2)^\eps$ and $\chi_p\ll1$,
the convexity bound and
\eqref{eq:holomorphic-maximal-new-mollifier-average}
bound the first contribution by
\begin{equation*}
\ll_{h,Q,\eta,J,\eps}
p^{1/2+\eps}(pL^2)^\eps
p^{1+\eps}L^{1+2\eps}
=p^{3/2+3\eps}L^{1+4\eps}.
\end{equation*}

For the dual contribution, \eqref{eq:local-Schur} and
\S\,\ref{subsec:kuz-maximal} give
\begin{equation*}
pA_{\phi,z}(p^2,1)
=p\lambda_\phi(p)^2p^{-2z}=p^{-2z},
\qquad |p^{-2z}|=1.
\end{equation*}
By \eqref{eq:holomorphic-dual-indices}, all remaining dual
indices are coprime to $p$.
Using \eqref{eq:holomorphic-hecke-coefficients} with $m=n=1$
and \eqref{eq:holomorphic-pieri-two}, the dual Dirichlet series is
\begin{equation*}
p^{-2z}\overline{M_\Pi^2}\,
L^{(p)}(1/2+u-z,\phi)^2
\zeta^{(p)}(1/2+u+2z)^2.
\end{equation*}
Compared with the full $L$-functions in
\eqref{eq:holomorphic-maximal-new-first-series}, the omitted
Euler factors contribute
\begin{equation*}
\begin{aligned}
&\left|
 (1-\lambda_\phi(p)p^{-1/2-u+z})^2
 (1-p^{-1/2-u-2z})^2
 \right|\\
&\qquad\leq(1+p^{-1-\eps})^2(1+p^{-1/2-\eps})^2\ll1
 \qquad(\Re(u)=\eps).
\end{aligned}
\end{equation*}
For the dual weight $\mathcal V_{-\mu}^-$, the numerator gamma
factor is again $\gamma_\mu(1/2+u)^2$, by
\eqref{eq:holomorphic-square-weight}.
The same pole cancellation therefore permits the contour shift
to $\Re(u)=\eps$.
Now $|\overline{M_\Pi^2}|=|M_\Pi|^2$, and the length factor is
$(p/L^2)^\eps\leq(pL^2)^\eps$.
The convexity bound and
\eqref{eq:holomorphic-maximal-new-mollifier-average},
together with the preceding uniform bound for the omitted
Euler factors, give
\begin{equation*}
\ll_{h,Q,\eta,J,\eps}
p^{1/2+\eps}(p/L^2)^\eps
p^{1+\eps}L^{1+2\eps}
=p^{3/2+3\eps}L.
\end{equation*}
This proves the proposition.
\end{proof}

\subsubsection{Proof of Theorem~\ref{thm:mollified-holomorphic-square}}
Insert \eqref{eq:holomorphic-square-diagonal},
Lemmas~\ref{lem:holomorphic-moment-short} and
\ref{lem:holomorphic-moment-long}, and
Propositions~\ref{prop:holomorphic-level-one-remainders} and
\ref{prop:holomorphic-maximal-new} into
\eqref{eq:holomorphic-moment-decomposition}.
Fix the auxiliary parameters in terms of $\eta$ as required
by the preceding lemmas. Since $L=p^\eta$, $\eta<1/2$ and
$\eps$ is arbitrarily small, the total spectral remainder is
\begin{equation*}
O_{h,Q,\eta,J,\eps}\big(p^{3/2+\eps}L^{1+\eps}\big)
=o_{h,Q,\eta,J}\!\left(\frac{p^2}{\log p}\right).
\end{equation*}
Division by $\chi_p=P/p^2$ replaces the diagonal factor $P$ by
$p^2$. Replacing $p^2$ by $P$ introduces an error $O_{h,Q}(p)$.
The fixed auxiliary parameters are absorbed into the dependence on $\eta$.
This proves \eqref{eq:mollified-holomorphic-square}.
Each Euler factor is positive and at most $1$, and differs from
$1$ by $O(q^{-5/2})$. The final assertions follow.

\section{Proof of Theorem~\ref{thm:intro-weighted-nonvanishing}}
\label{sec:harmonic-nonvanishing}

We now combine the mollified moment estimates of
\S\,\ref{sec:mollified-moments} with the total harmonic mass.
The functional equation relates the moments of
$L(1/2,\pi)\widehat M_\pi$ to those of $L(1/2,\pi)M_\pi$.
The Cauchy--Schwarz inequality then gives the
following lower bound for the harmonic non-vanishing proportion.

\begin{theorem}
\label{thm:harmonic-nonvanishing}
Let $h$ be the fixed non-negative test function in \eqref{eq:positive-moment-test}. For every $\eps>0$ and all sufficiently large primes $p$, we have
\begin{equation}
c_{\rm cusp}\sum_{\substack{\pi\in\mathcal B(p)\\L(1/2,\pi)\ne0}}
 \dfrac{h(\mu_\pi)}{L(1,\pi,\mathrm{Ad})}
\geq \left(\frac{1}{2}-\eps\right)c_{\rm cusp}\sum_{\pi\in\mathcal B(p)}
 \dfrac{h(\mu_\pi)}{L(1,\pi,\mathrm{Ad})}.
\label{eq:nonvanishing}
\end{equation}
In particular the non-vanishing harmonic mass is $\gg_h p^2$.
\end{theorem}

\subsection{The total harmonic mass}
\label{subsec:total-harmonic-mass}

We first evaluate the total harmonic mass.

\begin{lemma}
\label{lem:total-harmonic-mass}
For every $\eps>0$, we have
\begin{equation}
c_{\rm cusp}\sum_{\pi\in\mathcal{B}(p)}
 \dfrac{h(\mu_\pi)}{L(1,\pi,\mathrm{Ad})}
=P\mathcal H_h+O_{h,\eps}(p^{1+\eps}).
\label{eq:harmonic-mass}
\end{equation}
In particular, this sum is $\asymp_h p^2$ for sufficiently large $p$.
\end{lemma}

\begin{proof}
Apply Theorem~\ref{thm:kuznetsov} with Fourier indices $(m_1,m_2,n_1,n_2)=(1,1,1,1)$ and test function $h$. The diagonal term is $P\mathcal H_h$. For the short Weyl terms, the modulus pairs are $(d^2,d)$
and $(d,d^2)$ with $p| d$. The kernel arguments have modulus $d^{-3}$.
For the long Weyl term, we have $|y_1^2y_2|=D_1^{-3}$ and $|y_1y_2^2|=D_2^{-3}$. Using the trivial Kloosterman bounds $|\widetilde S|\leq D_1D_2$ and $|S^{(p)}|\leq(D_1D_2)^2$, \eqref{eq:kuz-short-bound} and \eqref{eq:w6-small} give 
\begin{equation*}
\begin{aligned}
|S_4|+|S_5|
&\ll_{h,A}P\sum_{\substack{d\geq1\\p\mid d}}d^{-3A}
 \ll_{h,A}p^{2-3A},\end{aligned}
 \end{equation*}
 and
\begin{equation*}
\begin{aligned}
|S_6|
&\ll_{h,A,\eps}P
 \sum_{\substack{D_1,D_2\geq1\\p\mid D_1,D_2}}
 (D_1D_2)^{1+\eps}\max(D_1,D_2)^{-3A}
 \ll_{h,A,\eps}p^{4-3A+2\eps} 
\end{aligned}
\end{equation*}
for sufficiently large $A$. 
Taking $A$ sufficiently large makes all three terms
$O_{h,\eps}(p^{-1})$.

For each level one sector, \eqref{eq:kuz-old-coeff} and the oldspace summation in Lemma \ref{lem:oldspace-cancellation} give
\begin{equation*}
\sum_{j=0}^2|B_\Pi^{(j)}(1,1)|^2
=(Q_\Pi^{-1})_{00}\ll1.
\end{equation*}
Together with \eqref{eq:level-one-normalizer-bounds}, \eqref{eq:level-one-counting-bound} and the rapid decrease of $h$, this yields
\begin{equation*}
|\mathcal{C}_{\rm old}|+|\mathcal{E}_{\min}|
+|\mathcal{E}_{\max,\rm old}|\ll_h 1.
\end{equation*}
The maximal new coefficient satisfies $A_{\phi,z}(1,1)=1$. 
Taking $(a,b,m,n)=(1,1,1,1)$ in \eqref{eq:gl2-induced-coefficient-average} and $\chi_p\ll1$, we have
\begin{equation*}
\chi_p|\mathcal{E}_{\max,\rm new}|\ll_{h,\eps}p^{1+\eps}.
\end{equation*}
Inserting these estimates into \eqref{eq:Kuz-full}, we obtain
\begin{equation*}
\chi_p c_{\rm cusp}\sum_{\pi\in\mathcal B(p)}
 \dfrac{h(\mu_\pi)}{L(1,\pi,\mathrm{Ad})}
=P\mathcal H_h+O_{h,\eps}(p^{1+\eps}).
\end{equation*}
Since $P/\chi_p=p^2$ and $P-p^2=p+1$, division by $\chi_p\asymp1$ proves \eqref{eq:harmonic-mass}.
The final assertion follows from $\mathcal H_h>0$.
\end{proof}

\subsection{Proof of Theorems~\ref{thm:harmonic-nonvanishing} and~\ref{thm:intro-weighted-nonvanishing}}

\begin{proof}
Fix a square-free $Q$ and $0<\eta<1/2$, independently of $p$, and put $L=p^\eta$. By the functional equation from \eqref{eq:functional-equation} and \eqref{eq:two-piece-mollifier},
\begin{equation*}
L(1/2,\pi)\widehat M_\pi =L(1/2,\pi)M_\pi+\overline{L(1/2,\pi)M_\pi}=2\Re\bigl(L(1/2,\pi)M_\pi\bigr).
\end{equation*}
Thus
\begin{equation*}
|L(1/2,\pi)\widehat M_\pi|^2 =2|L(1/2,\pi)M_\pi|^2  +2\Re\bigl((L(1/2,\pi)M_\pi)^2\bigr).
\end{equation*}
Summing these identities with the factors $c_{\rm cusp}h(\mu_\pi)/L(1,\pi,\mathrm{Ad})$ gives
\begin{equation}
\begin{aligned}
c_{\rm cusp}\sum_{\pi\in\mathcal B(p)}
 \dfrac{h(\mu_\pi)L(1/2,\pi)\widehat M_\pi}
 {L(1,\pi,\mathrm{Ad})}
&=2\Re\mathcal M_1(p),\\
c_{\rm cusp}\sum_{\pi\in\mathcal B(p)}
 \dfrac{h(\mu_\pi)|L(1/2,\pi)\widehat M_\pi|^2}
 {L(1,\pi,\mathrm{Ad})}
&=2\mathcal M_2(p)+2\Re\mathcal M_{1,1}(p).
\end{aligned}
\label{eq:two-piece-moment-identities}
\end{equation}
For sufficiently large $p$, \eqref{eq:mollified-first-result} gives $\Re\mathcal M_1(p)>0$, so the second moment in \eqref{eq:two-piece-moment-identities} is positive.
When $L(1/2,\pi)=0$, we have $L(1/2,\pi)\widehat M_\pi=0$. Then the Cauchy--Schwarz inequality with the non-negative factors $h(\mu_\pi)/L(1,\pi,\mathrm{Ad})$ yields
\begin{equation}
c_{\rm cusp}
\sum_{\substack{\pi\in\mathcal B(p)\\L(1/2,\pi)\ne0}}
 \dfrac{h(\mu_\pi)}{L(1,\pi,\mathrm{Ad})}
\geq
\dfrac{2\{\Re\mathcal M_1(p)\}^2} {\mathcal M_2(p)+\Re\mathcal M_{1,1}(p)}.
\label{eq:two-piece-cauchy}
\end{equation}

The Euler product in \eqref{eq:mollified-holomorphic-square} lies in $(0,1]$. Hence \eqref{eq:mollified-second-asymptotic} and \eqref{eq:mollified-holomorphic-square} give
\begin{equation*}
\mathcal M_2(p)+\Re\mathcal M_{1,1}(p)
\leq P\mathcal H_h
 \left(\zeta(2)(1+\eta^{-1})
       \prod_{q\mid Q}(1-q^{-2})+1\right)+O_{h,Q,\eta}\left(
 \dfrac{p^2\log(2+\log p)}{\log p}\right).
\end{equation*}
Inserting this bound and \eqref{eq:mollified-first-result} into \eqref{eq:two-piece-cauchy} and dividing by the total harmonic mass in \eqref{eq:harmonic-mass}, we then have
\begin{equation}
\dfrac{\displaystyle
 \sum_{\substack{\pi\in\mathcal B(p)\\L(1/2,\pi)\ne0}}
 \dfrac{h(\mu_\pi)}{L(1,\pi,\mathrm{Ad})}}
 {\displaystyle\sum_{\pi\in\mathcal B(p)}
 \dfrac{h(\mu_\pi)}{L(1,\pi,\mathrm{Ad})}}
\geq
\dfrac{2}{\zeta(2)(1+\eta^{-1})
 \prod_{q\mid Q}(1-q^{-2})+1}
-O_{h,Q,\eta}\left(\frac{\log(2+\log p)}{\log p}\right).
\label{eq:nonvanishing-fixed-parameters}
\end{equation}

Let $\eps>0$ be given. As $Q$ runs through products of all primes up to a bound tending to infinity, the Euler product for $\zeta(2)$ gives
\begin{equation*}
\zeta(2)\prod_{q\mid Q}(1-q^{-2})\longrightarrow1.
\end{equation*}
Since $1+\eta^{-1}\to3$ as $\eta\to(1/2)^-$, we may first choose
such a fixed $Q$, and then a fixed $0<\eta<1/2$, so that
\begin{equation*}
\dfrac{2}{\zeta(2)(1+\eta^{-1})
 \prod_{q\mid Q}(1-q^{-2})+1}
>\dfrac12-\dfrac\eps2.
\end{equation*}
With these choices fixed, take $p>Q$ sufficiently large. Then $p\nmid Q$ and the error in \eqref{eq:nonvanishing-fixed-parameters} is less than $\eps/2$. Multiplying \eqref{eq:nonvanishing-fixed-parameters} by the total harmonic mass gives \eqref{eq:nonvanishing}. The parameters $Q$, $\eta$ depend only on $\eps$. Thus the lower bound holds for all sufficiently large primes, with a threshold depending only on $h,\eps$. Lemma \ref{lem:total-harmonic-mass} also gives the asserted lower bound $\gg_h p^2$. This proves Theorem~\ref{thm:harmonic-nonvanishing}.
Cancelling the fixed positive constant $c_{\rm cusp}$ proves Theorem~\ref{thm:intro-weighted-nonvanishing}.
\end{proof}

\bigskip
\noindent{\bf Acknowledgements}
The authors would like to thank Professors Qingfeng Sun and Zhining Wei for their help and encouragement.

\end{document}